\documentclass[10pt]{article}
\usepackage{latexsym,amsmath,amssymb,graphics,amscd}

\addtocounter{footnote}{1}
\makeatletter
\@addtoreset{figure}{section}
\def\thefigure{\thesection.\@arabic\c@figure}
\def\fps@figure{h,t}
\@addtoreset{table}{bsection}
\def\thetable{\thesection.\@arabic\c@table}
\def\fps@table{h, t}
\@addtoreset{equation}{section}

\makeatother

\newtheorem{theorem}{Theorem}[section]
\newtheorem{definition}[theorem]{Definition}

\newtheorem{remark}[theorem]{Remark}
\newtheorem{proposition}[theorem]{Proposition}
\newtheorem{corollary}[theorem]{Corollary}
\newtheorem{example}[theorem]{Example}

\newcommand{\bfi}{\bfseries\itshape}
\newsavebox{\savepar}

\makeatletter

\begin{document}

\title{\textbf{Reduction, reconstruction, and skew-product decomposition of symmetric stochastic
differential equations}}
\author{Joan-Andreu L\'{a}zaro-Cam\'{\i}$^{1}$ and Juan-Pablo Ortega$^{2}$}
\date{}
\maketitle

\begin{abstract}
We present reduction and reconstruction procedures for the solutions of symmetric
stochastic differential equations, similar to those available for ordinary
differential equations. Additionally, we use the local tangent-normal decomposition, available when the symmetry group is proper, to construct local skew-product splittings in a neighborhood of any point in the open and dense principal orbit type. The general methods introduced in the first part of the
paper are then adapted to the Hamiltonian case, which is studied with special care and illustrated with several examples. The Hamiltonian category deserves a separate study since in that situation the presence of symmetries implies in most cases the existence of conservation laws, mathematically described via momentum maps,  that should be taken into account in the analysis.
\end{abstract}

\makeatletter
\addtocounter{footnote}{1} \footnotetext{%
Departamento de F\'{\i}sica Te\'{o}rica. Universidad de Zaragoza. Pedro
Cerbuna, 12. E-50009 Zaragoza. Spain. {\texttt{lazaro@unizar.es}}} 
\addtocounter{footnote}{1} \footnotetext{%
Centre National de la Recherche Scientifique, D\'{e}partement de Math\'{e}%
matiques de Besan\c{c}on, Universit\'{e} de Franche-Comt\'{e}, UFR des
Sciences et Techniques. 16, route de Gray. F-25030 Besan\c{c}on cedex.
France. {\texttt{Juan-Pablo.Ortega@univ-fcomte.fr} }}
\makeatother

\medskip

\normalsize

\section{Introduction}

Symmetries have historically played a role of paramount importance in the study of
dynamical systems in general (see~\cite{g2,golubitsky stewart 02, ChoLau}, and
references therein) and of physical, mechanical, and Hamiltonian systems in
particular (see for instance~\cite{fom, symmetry and mechanics, hsr} for general
presentations of the subject, historical overviews, and references). The presence of
symmetries in a system usually bring in its wake the occurrence of degeneracies,
conservation laws, and invariance properties that can be used to simplify or
\emph{reduce} the system and hence its analysis. In trying to pursue this strategy,
researchers have developed powerful mathematical tools that
optimize the benefit of this approach in specific situations.

The impressive
volume of work that has been done in this field over the centuries does not have a
counterpart in the context of stochastic dynamics, probably because most symmetry
based mathematical tools are formulated using global analysis and Lie theory in
an essential way, and this machinery has been adapted to the stochastic
context  relatively recently~\cite{meyer 81, meyer 82, schwartz 82, elworthy,
emery}. As we will show in this paper, \emph{most of the symmetry based techniques
available for  dynamical systems can be formulated and taken advantage of when
studying stochastic differential equations}.

In a first approach, symmetry based techniques can be roughly grouped into two
separate procedures, namely, \emph{reduction} and \emph{reconstruction}. Reduction
is explicitly implemented by combining the restriction of the
system to dynamically invariant submanifolds whose existence is implied by its
symmetries and by eliminating the remaining symmetry degeneracies through projection
to an appropriate orbit space. Even if the  space in which the system is originally
formulated is Euclidean, the resulting reduced space is most of the time a
non-Euclidean manifold hence showing the importance of global analysis in this
context. The reduction procedure yields a dimensionally smaller space in which the
symmetry degeneracies have been eliminated and that should, in principle, be easier
to study; in the stochastic context, reduction has the added value of being able in
some instances to isolate the non-stochastic part of the dynamics (see the example on
collective motion in Section~\ref{Stochastic collective Hamiltonian motion}).

If once the reduced system has been solved we want to come back to the original one,
we need to reconstruct the reduced solutions. In practice, this is obtained by
horizontally lifting the reduced motion using a connection and then correcting the
result with a curve in the group that satisfies a certain first order differential
equation. The strategy of combining reduction and reconstruction in the search for
the solutions of a symmetric dynamical system, splits the task into two parts, which
most of the time simplifies greatly the problem.

Another approach used to take advantage of the symmetries of a problem consists of using the Slice Theorem~\cite{pa} and the tangent-normal decomposition~\cite{krupa, field91} available for proper group actions to locally split the dynamics into a direction tangent to the group orbits and another one transversal to them. We will see that this tool, that is used in a standard fashion in the context of deterministic equivariant dynamics and equivariant bifurcation theory, yields in the stochastic case {\it skew-product splittings} that have already been extensively studied in the equivariant diffusions literature (see for instance \cite{Pauwels-Rogers, liao skew-product, taylor skew-product}, and references therein)
to construct decompositions of the associated second order differential operators. 

It must be noticed that the mathematical value of the results obtained with the two approaches that we just briefly discussed, that is, the one based on reduction-reconstruction and the one based on the tangent-normal decomposition, is morally the same. However, there are important technical conditions that make them different and preferable over one another in different specific situations: 
\begin{description}
\item [(i)]  The reduction-reconstruction technique uses very strongly the orbit space of the symmetry group in question; this space could be geometrically convoluted and we may need to use only its strata if we want to face regular quotient manifolds where the standard calculus on manifolds is valid. The main advantage of this technique is that it yields {\it global results}.
\item [(ii)] The use of the Slice Theorem and the tangent-normal decomposition makes unnecessary the use of quotient manifolds and the entire analysis takes place in the original manifold. However, the results obtained are {\it local and are limited to a tubular neighborhood of the orbits}.
\end{description}

In this paper we show how the symmetries of stochastic differential equations can be
used by implementing techniques similar to those available for
their deterministic counterparts. We start in Section~\ref{Symmetries and
conservation laws of stochastic differential equations} by introducing the notion of
group of symmetries of a stochastic differential equation and by studying the
associated invariant submanifolds as well as the implied degeneracies in the solutions. The reduction and reconstruction procedures are presented in Section~\ref{Reduction and reconstruction}; reconstruction is carried
out using the horizontal lifts for semimartingales introduced by~\cite{Shigekawa,
Catuogno}.

The skew-product decomposition of second order differential operators is a factorization technique that has been used in the stochastic processes literature in order to split the semielliptic and, in particular, the diffusion operators, associated to certain stochastic differential equations (see, for instance,~\cite{Pauwels-Rogers, liao skew-product, taylor skew-product}, and references therein). This splitting has important consequences as to the  properties of the solutions of these equations, like certain factorization properties of their probability laws and of the associated stochastic flows. In Section~\ref{Symmetries and skew-product decompositions} we show that symmetries are a natural way to obtain this kind of decompositions. Our work extends the existing results in two ways: first, we generalize the notion of skew-product to arbitrary stochastic differential equations by working with the notion of skew-product decomposition of the Stratonovich operator. Obviously, our approach coincides with the traditional one in the case of diffusions. Second, we use the Slice Theorem~\cite{pa} and the tangent-normal decomposition~\cite{krupa, field91}  to construct local skew-product decompositions in the presence of arbitrary proper symmetries (not necessarily free) in a neighborhood of any point in the open and dense principal orbit type. This result generalizes  the skew-product decompositions presented in~\cite{elworthy principal bundles} for regular free actions. Section~\ref{Projectable stochastic differential equations on associated bundles} studies stochastic differential equations on associated bundles; in this situation the local skew-product splitting induced by the Slice Theorem is  globally available.

Section~\ref{The Hamiltonian case} is dedicated to reduction and
reconstruction in the stochastic Hamiltonian category. Stochastic Hamiltonian
systems where introduced in~\cite{bismut 81} and generalized in~\cite{lo} to
accommodate non-Euclidean phase spaces and stochastic components modeled by
arbitrary semimartingales and not just Brownian motion. Given the generic
non-Euclidean character of reduced spaces, the generalization in~\cite{lo} is
in this context of much relevance. It is worth mentioning that, as it was
already the case for deterministic Hamiltonian systems, stochastic Hamiltonian
systems  are stable with respect to symplectic and Poisson
reduction; in short, the reduction of a stochastic Hamiltonian system is again a
stochastic Hamiltonian system. In Section~\ref{examples we see} we present
several (Hamiltonian) examples. The first one (Section~\ref{Stochastic collective
Hamiltonian motion}) has to do with deterministic systems in which a
stochastic perturbation  is added using the conserved quantities associated to
the symmetry (collective perturbation); such systems share the remarkable feature
that  symplectic reduction eliminates the stochastic part of the equation making the
reduced system deterministic. In Section~\ref{Stochastic mechanics on Lie groups} we
study the symmetries of stochastic mechanical systems on the cotangent bundles of Lie
groups. In this situation, the reduction and reconstruction equations can be written
down in a particularly explicit fashion that has to do with the Lie-Poisson
structure in the dual of the Lie algebra of the group in question. A particular case
of this is presented in Section~\ref{Stochastic perturbations of the free rigid
body} where we analyze two different stochastic  perturbations of the free rigid
body: one of them models the dynamics of a free rigid body subjected to small random
impacts and the other one an "unbolted"   rigid body that is not completely rigid.

\section{Symmetries and conservation laws of stochastic differential equations}
\label{Symmetries and conservation laws of stochastic differential equations}

Let $M $ and $N $ be two finite dimensional manifolds and let $(\Omega,
\mathcal{F},\{\mathcal{F}_t \,|\,t \geq 0\}, P)$ be a filtered probability space.
Let $X: \mathbb{R}_+ \times \Omega \rightarrow  N$ be a $N$-valued semimartingale.
Using the conventions in~\cite{emery}, 
a {\bfi  Stratonovich operator} from $N$ to $M$ is a family $\{
S(x,y)\}_{x \in N, y \in  M }$ such that $S(x,y) :T _xN \rightarrow T
_yM$ is a linear mapping that depends smoothly on its two entries. Let $S ^\ast
(x,y): T ^\ast _y M \rightarrow T _x ^\ast  N $ be the adjoint of $S(x,y)$. 

We recall that a
$M$-valued semimartingale $\Gamma$ is a solution of the the Stratonovich
stochastic differential equation
\begin{equation}
\label{equations differential form 1}
\delta \Gamma= S(X,\Gamma) \delta X
\end{equation}
associated to $X$ and $S$, if for any $\alpha \in  \Omega (M) $, the following
equality between Stratonovich integrals holds:
\[
\int   \langle \alpha, \delta \Gamma\rangle =\int  \langle  S ^\ast (X,\Gamma)
\alpha,
\delta X\rangle.
\]
We will refer to $X$ as the {\bfi  noise semimartingale} or the {\bfi  stochastic
component} of the stochastic differential equation~(\ref{equations differential form
1}). It can be shown~\cite[Theorem 7.21]{emery} that in this setup,
given a 
$\mathcal{F} _0$ measurable random variable $\Gamma _0$, there are a  stopping time $\zeta>0 $ and
a solution
$\Gamma $ of~(\ref{equations differential form 1}) with initial condition $\Gamma _0$ defined on the set
$\{(t, \omega)\in 
\mathbb{R}_+ \times  \Omega\mid t \in  [0, \zeta (\omega))\}$ that has the following
maximality and uniqueness property: if   $\zeta' $ is another
stopping time such that  $\zeta'< \zeta $ and $\Gamma' $ is another solution defined on  $\{(t,\omega)\in 
\mathbb{R}_+ \times  \Omega\mid t \in  [0, \zeta' (\omega))\}$, then $\Gamma'$ and $\Gamma$ coincide
in this set. If $\zeta $ is finite then $\Gamma$ explodes at time $\zeta$, that is, the path $\Gamma _t $
with $t \in  [0, \zeta) $ is not contained in any compact subset of $M$. If the manifold $M$ is compact then all the solutions of any stochastic differential equation defined on $M$ are defined for all time. Since this is a hypothesis that we are not willing to adopt, the reader should keep in mind that all the solutions that we will work with are defined only up to a maximal stopping time, even if this is not explicitly mentioned.

We also recall that stochastic differential equations can be formulated using
It\^o integration by associating a natural {\bfi  Schwartz operator} $ {\cal S}:
\tau_xN\rightarrow \tau _yM  $ on the second order tangent bundles, to the
Stratonovich operator
$S$; see~\cite{emery} and references therein for the details.

\begin{definition}
\label{definicion simetria dinamica}
Let $X: \mathbb{R}_+ \times \Omega \rightarrow  N$ be a $N$-valued semimartingale
and let $S:TN \times  M \rightarrow TM $ be a Stratonovich operator. Let
$\phi:M\rightarrow M$ be a diffeomorphism. We say that $\phi$ is a
{\bfseries\itshape symmetry} of the stochastic differential equation (\ref{equations
differential form 1}) if for any  $x \in N $ and $y\in M$
\begin{equation}
S\left( x,\phi\left( y\right) \right) =T_{y}\phi\circ S\left( x,y\right) . 
\label{definicion H invariante}
\end{equation}
\end{definition}

As it was already the case in standard deterministic context, the symmetries of a
stochastic differential equation imply degeneracies at the level of its solutions, as
we spell out in the following proposition.

\begin{proposition}
\label{symmetry implies degeneracy}
Let $X: \mathbb{R}_+ \times \Omega \rightarrow  N$ be a $N$-valued
semimartingale, $S:TN \times  M \rightarrow TM $ a Stratonovich operator, and let $\phi:M\rightarrow
M$ be a symmetry of the corresponding stochastic differential
equation~(\ref{equations differential form 1}).
If $\Gamma$  is
solution of~(\ref{equations
differential form 1}) then so is $\phi\left( \Gamma\right) $.
\end{proposition}

\noindent\textbf{Proof.\ \ }
Let $\Gamma$ be a solution of (\ref{equations differential form 1}). We need
to show that for any $\alpha\in\Omega\left( M\right) $,%
\begin{equation*}
\int\langle\alpha,\delta\phi\left( \Gamma\right) \rangle=\int\langle
S^{\ast}\left( X,\phi\left( \Gamma\right) \right) \alpha,\delta X\rangle. 
\end{equation*}
Since $\phi$ is a diffeomorphism, 
$\int\left\langle \alpha,\delta\phi\left( \Gamma\right) \right\rangle
=\int\left\langle \phi^{\ast}\alpha,\Gamma\right\rangle $ (see, for
instance,~\cite[\S 7.5]{emery}). Now, since
$\Gamma
$ is a solution of (\ref{equations differential form 1}), $\int\left\langle
\phi^{\ast}\alpha,\Gamma\right\rangle =\int\left\langle S^{\ast}\left(
X,\Gamma\right) (\phi^{\ast}\alpha),\delta X\right\rangle $. Since $\phi $ is a
symmetry, we have that $S^{\ast}\left(
x,\phi\left( y\right) \right) =S^{\ast}\left( x,y\right) \circ
T_{y}^{\ast}\phi$, for any $x \in N $, $ y \in  M $ and hence,
\begin{equation*}
\int\left\langle \phi^{\ast}\alpha,\Gamma\right\rangle =\int\left\langle
S^{\ast}\left( X,\Gamma\right) (\phi^{\ast}\alpha),\delta X\right\rangle
=\int\left\langle S^{\ast}\left( X,\phi\left( \Gamma\right) \right)
(\alpha),\delta X\right\rangle, 
\end{equation*}
which shows that $\phi\left( \Gamma\right) $ is a solution
of (\ref{equations differential form 1}).
\quad $\blacksquare$

\bigskip

The symmetries that we are mostly interested in are induced by the action of a Lie
group $G$ on the manifold $M$ via the map
$\Phi:G\times M\rightarrow M$. Given $\left( g,z\right) \in G\times M$, we will
usually write $g\cdot z$ to denote $\Phi\left( g,z\right) $. We also introduce the
maps%
\begin{equation*}
\begin{array}{rcc}
\Phi_{z}:G & \longrightarrow & M \\ 
g & \longmapsto & g\cdot z%
\end{array}
,\qquad
\begin{array}{rcc}
\Phi_{g}:M & \longrightarrow & M \\ 
z & \longmapsto & g\cdot z
\end{array}
. 
\end{equation*}
The Lie algebra of $G$ will be usually denoted by $\mathfrak{g}$ and we will write the tangent space to the orbit $G \cdot  m  $ that contains $m \in  M $  as $ \mathfrak{g}\cdot m:= T _m(G \cdot m) $.

\begin{definition}
We will say that the stochastic differential equation~(\ref{equations differential
form 1}) is {\bfseries\itshape$G$-invariant} if, for any $g\in G$, the
diffeomorphism $\Phi _{g}:M \rightarrow  M$ is a symmetry in the sense of Definition
\ref{definicion simetria dinamica}. In this situation we will also say that the
Stratonovich operator
$S$ is $G$-{\bfi  invariant}.
\end{definition}

\begin{remark}
\normalfont
Given a solution $\Gamma  $ of a $G$-invariant stochastic differential equation,
Proposition~\ref{symmetry implies degeneracy} provides an entire orbit of solutions
since for any $g \in  G $, the semimartingale $\Phi_g (\Gamma) $ is also a
solution. This degeneracy has also a reflection in the probability laws of the solutions in a form that we spell out in the following lines. Let $\Gamma:\left\{  0\leq t<\zeta\right\}
\rightarrow M$ be a solution of the $G$-invariant system $\left(
M,S,X,N\right)  $ defined up to the explosion time $\zeta$, which may be
finite if $M$ is not compact. In such case, $\Gamma$ can be
actually understood as a process that takes values in the Alexandroff one-point
compactification $\hat{M}:=M\cup\left\{  \infty\right\}  $ of $M$ and it is hence
defined in the whole space $\mathbb{R}_{+}\times\Omega$ (\cite[Chapter
V]{ikeda}). In this picture, the process $\Gamma$ is continuous and with the property that
$\Gamma_{t}\left(  \omega\right)  =\left\{  \infty\right\}  $, for any $\left(
t,\omega\right)  \in\mathbb{R}_{+}\times\Omega$ such that $t\geq\zeta\left(
\omega\right)  $.

Let now $\hat{W}(M)$ be the path space defined by
\begin{align*}
\hat{W}(M)  &  =\{w:\left[  0,\infty\right]  \rightarrow\hat{M}\text{
continuous such that }w\left(  0\right)  \in M\text{ and}\\
&  \hspace{0.55cm}\text{if }w(t)=\left\{  \infty\right\}  \text{ then
}w(t^{\prime})=\left\{  \infty\right\}  \text{ for any }t^{\prime}\geq t\}.
\end{align*}
Let $\left\{  P_{z}~|~z\in
M\right\}  $ be the family of probability measures on $\hat{W}(M)$  defined by  the solutions of $\left(
M,S,X,N\right)  $, that is, $P_{z}$ is the law of the random variable $\Gamma
^{z}:\Omega\rightarrow\hat{W}(M)$, where $\Gamma^{z}$ is the solution of
$\left(  M,S,X,N\right)  $ with initial condition $\Gamma_{t=0}^{z}=z$ a.s..
The action $\Phi:G\times M\rightarrow M$ may be extended to $\hat{M}$ just
putting $\Phi_{g}\left(  \left\{  \infty\right\}  \right)  =\left\{
\infty\right\}  $ for any $g\in G$. Since $\Phi_{g}\left(  \Gamma^{z}\right)
$ is the unique solution of the system $\left(  M,S,X,N\right)  $ with initial
condition $g\cdot z$ by Proposition~\ref{symmetry implies degeneracy} then  $P_{g\cdot z}=\Phi_{g}^{\ast}P_{z}$. More explicitly, for any measurable set $A \subset \hat{W}(M)$, $P_{g\cdot z}(A)=P_{z}\left( \Phi_{g} (A)\right)$.

The equivariance property of
the probabilities $\left\{  P_{z}~|~z\in M\right\}  $ can be found in
\cite{elworthy principal bundles} formulated in the context of equivariant diffusions on principal bundles.
In that setup, the authors replace the path space $\hat{W}(M)$ by $C\left(  l,r,M\right)  =\left\{
\sigma:\left[  l,r\right]  \rightarrow M~|~\sigma\text{ is continuous}
\right\}  $, $0\leq l<r<\infty$ and prove~\cite[Theorem
2.5]{elworthy principal bundles} that the probability laws $\{P_{z}%
^{l,r}~|~z\in M\}$ admit a factorization through probability kernels
$\{P_{z}^{H,l,r}~|~z\in M\}$ from $M$ to $C\left(  l,r,M\right)  $ and
$\{Q_{w}^{l,r}~|~w\in C\left(  l,r,M\right)  \}$ from $C\left(  l,r,M\right)
$ to $C_{e}\left(  l,r,G\right)  =\left\{  \sigma:\left[  l,r\right]
\rightarrow G~|~\sigma\text{ is continuous, }\sigma(l)=e\right\}  $ such that%
\[
P_{z}^{l,r}(U)=\int\int\mathbf{1}_{U}(g\cdot w)Q_{w}^{l,r}(dg)P_{z}%
^{H,l,r}(dw)
\]
for any Borel set $U\subseteq C\left(  l,r,M\right)  $. The prove of this fact uses a technique very close to the reduction-reconstruction scheme that we will introduce in the next section.

\end{remark}

Apart from degeneracies, the presence of symmetry in a stochastic
differential equation is also associated with the occurrence of conserved quantities
and, more generally, with the appearance of invariant submanifolds.  

\begin{definition}
Let $\Gamma
$ be a solution of the stochastic differential equation~(\ref{equations differential
form 1}) and let $L $ be an immersed submanifold of $M$. Let $\zeta $ be the maximal
stopping time of $\Gamma  $ and suppose that
$\Gamma_{0}(\omega)=Z_{0}$, where $Z_{0}$ is a random variable such that
$Z_{0}(\omega)\in L$, for all $
\omega\in\Omega$. We say that $L$ is an {\bfseries
\itshape invariant submanifold} of the stochastic differential equation if for
any stopping time $\tau<\zeta$ we have that $
\Gamma_{\tau}\in L$. 
\end{definition}

\begin{proposition}
\label{invariant manifolds proposition} 
Let $X: \mathbb{R}_+ \times \Omega \rightarrow  N$ be a $N$-valued semimartingale
and let $S:TN \times  M \rightarrow TM $ be a Stratonovich operator.
Let $L$ be an immersed submanifold of $M$ and
suppose that the Stratonovich operator $S$ is such that $\mbox{\rm Im}\left(
S(x,y)\right) \subset T_{y}L$, for any $y\in L$ and any $x\in N$. Then,  $L$ is an
invariant submanifold of the stochastic differential equation~(\ref{equations
differential form 1}) associated to $X$ and $S$.
\end{proposition}

\noindent\textbf{Proof.\ \ } By hypothesis, the Stratonovich operator $
S:TN \times  M\rightarrow TM$ induces another Stratonovich operator $%
S_{L}:TN \times L\rightarrow TL$, obtained from $S$ by restriction. It
is clear that if $i:L\hookrightarrow M$ is the inclusion then 
\begin{equation}
S_{L}^{\ast}(x,y)\circ T_{y}^{\ast}i=S^{\ast}(x,y), 
\label{one and other operator 1}
\end{equation}
for any $x\in N$ and $y\in L$. Let $\Gamma_{L}$ be the semimartingale in $L$
that is a solution of the Stratonovich stochastic differential equation 
\begin{equation}
\label{equation for l 1}
\delta\Gamma_{L}=S_{L}(X,\Gamma_{L})\delta X   
\end{equation}
with initial condition $\Gamma_{0}$ in $L$. We now show that $\overline{\Gamma }%
:=i\circ\Gamma_{L}$ is a solution of 
\begin{equation*}
\delta\overline{\Gamma}=S(X,\overline{\Gamma})\delta X. 
\end{equation*}
which proves the statement. Indeed, for any $\alpha
\in\Omega(M)$, 
\begin{equation*}
\int\langle\alpha,\delta\overline{\Gamma}\rangle=\int\langle\alpha
,\delta(i\circ\Gamma_{L})\rangle=\int\langle i^{\ast}\alpha,\delta\Gamma
_{L}\rangle. 
\end{equation*}
Since $\Gamma_{L}$ satisfies (\ref{equation for l 1}) and $i^{\ast}\alpha
\in\Omega(L)$, by (\ref{one and other operator 1}) this equals 
\begin{equation*}
\int\langle S_{L}^{\ast}(X,\Gamma_{L})(i^{\ast}\alpha),\delta X\rangle
=\int\langle S^{\ast}(X,i\circ\Gamma_{L})(\alpha),\delta
X\rangle=\int\langle S^{\ast}(X,\overline{\Gamma})(\alpha),\delta X\rangle, 
\end{equation*}
that is, $\delta\overline{\Gamma}=S(X,\overline{\Gamma})\delta X$, as
required. \quad$\blacksquare$

\medskip 

We now use Proposition~\ref{invariant manifolds proposition} to
show that the invariant manifolds that can be associated to deterministic
symmetric systems are also available in the stochastic context. Let $M$ be a
manifold acted properly upon by a Lie group $G$ via the map $\Phi:G\times
M\rightarrow M$. We recall that the action $\Phi $ is said to be proper when for any two convergent sequences $\{m _n \} $ and $\{g _n \cdot m _n:= \Phi(g _n, m _n) \}$ in $M$, there exists a convergent subsequence $\{g _{n_k}\} $ in $G$.
The properness hypothesis on the action implies implies that most of the useful features that compact group actions have, are still available. For example, proper group actions admit local slices, the isotropy subgroups are always compact, and  (the
connected components of) the {\bfseries\itshape isotropy type} submanifolds
defined by $M_{I}:=\{z\in M\mid G_{z}=I\}$, are embedded submanifolds of $M$
for any isotropy subgroup $I\subset G$ of the action.

\begin{proposition}[Law of conservation of the isotropy]
\label{law conservation isotropy}
Let $X: \mathbb{R}_+ \times \Omega \rightarrow  N$ be a $N$-valued semimartingale
and let $S:TN \times  M \rightarrow TM $ be a Stratonovich operator that is
invariant with respect to a proper action of the Lie group $G$ on the manifold $M$. 
Then, for any isotropy subgroup $I\subset G$, the isotropy type
submanifolds $M_{I}$ are invariant submanifolds of the stochastic differential
equation associated to $S$ and $X$.
\end{proposition}

\noindent\textbf{Proof.\ \ } The properness of the action guarantees that for any
isotropy subgroup $I\subset G$ and any $z \in M _I$,
\begin{equation}
\label{characterization of tangent isotropy}
T _z M _I=(T _zM)^I:=\{ v  \in  T _zM\mid T _z \Phi_g \cdot v=v,\text{ for any $g \in
I $}\}.
\end{equation}
Hence, for any $z \in M _I $ and $g \in I $, the $G$-invariance of the Stratonovich
operator
$S$ implies that
\begin{equation*}
T_{z}\Phi_{g}\circ S\left( x,z\right)= S\left( x,g\cdot z\right) =S\left( x,z\right),
\end{equation*}
which by~(\ref{characterization of tangent isotropy}) implies that ${\rm Im} \left(S\left(
x,z\right)\right) \subset  T _z M _I $. The invariance of the isotropy type
manifolds follows then from Proposition~\ref{invariant manifolds proposition}
. \quad$\blacksquare$

\begin{remark}
\normalfont
Some of the results that we just stated and others that will appear later on in the paper could be easily proved using their deterministic counterparts and the so called {\it Malliavin's Transfer Principle}~\cite{malliavin transfer} which says, roughly speaking, that results from the theory of ordinary differential equations are valid for stochastic differential equations in Stratonovich form. The unavailability of a metatheorem that explicitly proves and shows the range of applicability of this principle makes advisable its use with care.
\end{remark}

\section{Reduction and reconstruction}
\label{Reduction and reconstruction}

This section is the core of the paper. In the preceding paragraphs we explained how
the symmetries of a stochastic differential equation imply the existence of certain
conservation laws and degeneracies; reduction is a natural procedure to take
advantage of the former and eliminate the latter via a combination of restriction
and passage to the quotient operations. The end result of this strategy is the
formulation of a stochastic differential equation with the same noise semimartingale
but whose solutions take values in a manifold that is dimensionally smaller than the
original one, which justifies the term \emph{reduction} when we refer to this
process. Smaller dimension and the absence of symmetry induced degeneracies usually
make the reduced stochastic differential equation more tractable and easier to
solve. The gain is therefore clear if once we have found the solutions of the
reduced system, we know how to use them to find the solutions of the original system;
that task is feasible  and  is the \emph{reconstruction} process that will be
explained in the second part of this section.

\begin{theorem}[Reduction Theorem]
\label{reduction theorem statement general systems} 
Let $X: \mathbb{R}_+ \times \Omega \rightarrow  N$ be a $N$-valued semimartingale
and let $S:TN \times  M \rightarrow TM $ be a Stratonovich operator that is
invariant with respect to a proper action of the Lie group $G$ on the manifold $M$.
Let $I \subset G$ be an isotropy subgroup of the $G$-action on $M$, $M _I $ the
corresponding isotropy type submanifold, and $L _I:=N (I)/I $, with $N(I):=\{g \in
G\mid gI g ^{-1}=I\} $ the normalizer of $I$ in $G$. $L _I $  acts freely and
properly on
$M _I $  and hence the orbit space $M _I/L _I $ is a regular quotient manifold, that
is, the projection $\pi_I: M _I \rightarrow  M _I/ L _I  $ is a surjective
submersion. Moreover, there is a well defined 
Stratonovich operator $S_{M_I/L _I}:TN\times M_I/L _I\rightarrow T\left(M_I/L
_I\right)
$ given by
\begin{equation}
S_{M_I/L _I}(x,\pi_I(z))=T_{z}\pi_I\left( S(x,z)\right) ,\quad\text{for any $x\in N$
and $z\in M_I$}   \label{general reduced Stratonovich operator 1}
\end{equation}
such that if $\Gamma$ is a solution semimartingale of the stochastic differential
equation associated to $S$ and
$X$, with initial condition $\Gamma_{0}\subset M _I$, then so is $\Gamma_{M_I/L
_I}:=\pi_I\left( \Gamma\right) $ with respect to $S_{M_I/L _I}$ and  $X$, with
initial condition $\pi_I(\Gamma_{0})$. We will refer to  $S_{M_I/L _I} $ as the
{\bfi  reduced Stratonovich operator} and to $\Gamma_{M_I/L
_I} $ as the {\bfi  reduced solution}.
\end{theorem}

\noindent\textbf{Proof.\ \ } The statement about $M _I/L _I $ being a regular
quotient manifold is a standard fact about proper group actions on manifolds (see
for instance~\cite{dk}). Now, observe that $S_{M_I/L_I}:TN\times M_I/L_I\rightarrow
T\left( M_I/L_I\right) $ is well defined: if $z_{1}$, $z_{2}\in M_I$ are such that 
$\pi_I\left( z_{1}\right) =\pi_I\left( z_{2}\right) $, then there exists some $
g\in L _I$ satisfying $z_{2}=\Phi_{g}\left( z_{1}\right) $ (we use the same symbol
$\Phi $ to denote the $G$-action on $M$ and the induced $L _I$-action on  $M_I $).
Hence,
\begin{equation*}
S_{M_I/L_I}(x,\pi_I(z_{2}))=T_{z_{2}}\pi_I\circ S(x,z_{2})=T_{z_{2}}\pi_I\circ
T_{z_{1}}\Phi_{g}\circ S\left( x,z_{1}\right) =T_{z_{1}}\pi_I\circ S\left(
x,z_{1}\right) =S_{M_I/L_I}(x,\pi_I(z_{1})), 
\end{equation*}
where the $G$-invariance of $S$ has been used. Let now $\Gamma$ be a
solution semimartingale of the stochastic differential equation
associated to $S$ and $X $ with initial condition $\Gamma_{0}\subset M _I$. The
$G$-invariance of $S$ implies via Proposition~\ref{law conservation isotropy}
that $\Gamma \subset M _I $ and hence $\Gamma_{M_I/L
_I}:=\pi_I\left( \Gamma\right) $ is well defined.
In order
to prove the statement, we have to check that for any one-form
$\alpha\in\Omega(M_I/L_I)$ 
\begin{equation*}
\int\langle\alpha,\delta\Gamma_{M_I/L_I}\rangle=\int\langle S_{M_I/L_I}^{\ast
}(X,\Gamma_{M_I/L_I})\alpha,\delta X\rangle. 
\end{equation*}
This equality follows in a straightforward manner from (\ref{general reduced
Stratonovich operator 1}). Indeed, 
\begin{multline*}
\int\langle\alpha,\delta\Gamma_{M_I/L_I}\rangle=\int\langle\alpha,\delta\left(
\pi_I\circ\Gamma\right) \rangle=\int\langle\pi_I^{\ast}\alpha,\delta\Gamma
\rangle\\
=\int\langle S^{\ast}(X,\Gamma)\left( \pi_I^{\ast}\alpha\right) ,\delta
X\rangle=\int\langle S_{M_I/L_I}^{\ast}(X,\Gamma_{M_I/L_I})\alpha,\delta X\rangle, 
\end{multline*}
as required. \quad $\blacksquare$

\medskip

We are now going to carry out the reverse procedure, that is, given an isotropy
subgroup $I\subset G$ and a solution semimartingale $\Gamma_{M _I/L _I}$  of the
reduced stochastic differential equation with Stratonovich operator $S_{M_I/L _I}$
we will {\it reconstruct} a solution $\Gamma $ of the initial stochastic
differential equation  with Stratonovich operator $S $. In order to keep the
notation not too heavy we will assume in the rest of this section that the
$G$-action on $M$ is not only proper but also free, so that the only isotropy
subgroup is the identity element $e$ and hence there is only one isotropy type
submanifold, namely $M_e=M $. The general case can be obtained by replacing in the
following paragraphs $M $ by the isotropy type manifolds $M _I $, and $G $ by the
groups $L _I$. 

We now make our goal more precise. The freeness of the action $\Phi :G\times
M\rightarrow M$ guarantees that the canonical projection $\pi :M\rightarrow M/G$ 
is a principal bundle with structural group $G$.
We saw in the previous theorem that for any solution $\Gamma $ of a stochastic
differential equation associated to a $G$-invariant Stratonovich operator $S$  and a
$N$-valued noise semimartingale $X$, we can build a solution
$\Gamma _{M/G}=\pi
\left(
\Gamma
\right) $ of the reduced stochastic differential equation associated to the projected
Stratonovich operator
$S_{M/G}$ introduced in (\ref{general reduced Stratonovich operator 1}) and
to the stochastic component $X$. The main goal of the paragraphs that follow 
is to show how to reconstruct the dynamics of the initial system from solutions $\Gamma
_{M/G}$ of
the reduced system. As we will see in Theorem~\ref{teorema reconstruccion 1}, any
solution $\Gamma $ of the original stochastic differential equation may be written
as $\Gamma =\Phi _{g^{\Xi }}\left( d\right) $ where $d:\mathbb{R}_{+}\times
\Omega
\rightarrow M$ is a semimartingale such that $\pi \left( d\right) =\Gamma
_{M/G}$ and $g^{\Xi }:\mathbb{R}_{+}\times \Omega \rightarrow G$ is a $G$-valued
semimartingale which satisfies a suitable stochastic differential equation
on the group $G$.

We start by picking $A\in\Omega^{1}\left( M;\mathfrak{g}\right) $ ($\mathfrak{g} $
is the Lie algebra of $G$) an auxiliary principal connection on the left principal
$G$-bundle
$\pi:M\rightarrow M/G$ and let $TM={\rm Hor}\oplus{\rm Ver}$ be the decomposition of
the tangent bundle $TM $ into the Whitney sum of the horizontal and vertical bundles
associated to $A$. Analogously, the cotangent bundle $T^{\ast}M$ admits a
decomposition $T^{\ast}M={\rm Hor}^{\ast}\oplus{\rm Ver}^{\ast}$
where, by definition, ${\rm Hor}_{z}^{\ast}:=\left( {\rm Ver}
_{z}\right) ^{\circ}$ is the annihilator of the vertical subspace ${\rm Ver}%
_{z}$ at a point $z\in M$ and ${\rm Ver}^{\ast}_z:=\left( {\rm Hor}
_{z}\right) ^{\circ}$ is the annihilator of the horizontal subspace. Hence, any one
form $
\alpha\in\Omega\left( M\right) $ may be uniquely written as $%
\alpha=\alpha^{H}+\alpha^{V}$ with $\alpha^{H}\in{\rm Hor}^{\ast}$ and $%
\alpha^{V}\in {\rm Ver}^{\ast}$. A section of the bundle $%
\pi_{M}:T^{\ast }M\rightarrow M$ taking values in ${\rm Hor}^{\ast}$ is
called a horizontal one form. It is called vertical if $\alpha_{z}\in 
{\rm Ver}_{z}^{\ast}$ for any $z\in M$.

\smallskip Let $\Gamma_{M/G}\subset M_{M/G}$ be a solution of the reduced
stochastic differential equation associated to the Stratonovich operator $S_{M/G}$,
and with stochastic component
$X:\mathbb{R} _{+}\times\Omega\rightarrow V$ as in Theorem \ref{reduction theorem
statement general systems} . As we claimed, we are going to find a solution $\Gamma
$ to the original
$G$-invariant  stochastic differential equation associated to $S $, such that
$\pi\left(
\Gamma\right) =\Gamma_{M/G}$ with a given initial condition $\Gamma_{0}$. We
start by horizontally lifting $\Gamma_{M/G}$ to a $M$-valued semimartingale $%
d$. Indeed, by \cite[Theorem 2.1]{Shigekawa} (see also \cite{Catuogno}),
there exists a $M$-valued semimartingale $d: \mathbb{R}_+ \times \Omega \rightarrow 
M$ such that
$d_{0}=\Gamma_{0}$, $%
\pi\left( d\right) =\Gamma_{M/G}$ and that satisfies
\begin{equation}
\int\left\langle A,\delta d\right\rangle =0,   \label{ggdot 1}
\end{equation}
where (\ref{ggdot 1}) is a $\mathfrak{g}$-valued integral. More specifically, let
$\left\{
\xi_{1},...,\xi_{m}\right\} $ be a basis of the Lie algebra $\mathfrak{g}$ and let
$A\left( z\right) =\sum_{i=1}^{m}A^{i}\left( z\right) \xi_{i}$ the expression of $A$
in this basis. Then%
\begin{equation}
\int\left\langle A,\delta d\right\rangle :=\sum_{i=1}^m \int\left\langle
A^{i},\delta d\right\rangle \xi_i.   \label{ggdot 2}
\end{equation}
The condition (\ref{ggdot 1}) is equivalent to $\int\left\langle
\alpha,\delta d\right\rangle =0$ for any vertical one-form $\alpha\in\Omega\left(
M\right)
$ (see \cite[page 1641]{Catuogno}) which, in turn, implies
\begin{equation}
\int\left\langle \theta,\delta d\right\rangle =0   \label{ggdot 9}
\end{equation}
for any $T ^\ast M$-valued process $\theta:\mathbb{R}_{+}\times\Omega\rightarrow 
{\rm Ver}^{\ast}\subset T^{\ast}M$ over $d$. We want to find a $G$-valued
semimartingale $g^{\Xi}:\Omega\times\mathbb{R}_{+}\rightarrow G$ such
that $g_{0}^{\Xi}=e$ a.s. and $\Gamma=g^{\Xi}\cdot d$ is a solution of the
stochastic differential equation associated to the Stratonovich operator $S$ and the
$N$-valued noise semimartingale $X$.

\smallskip Let $g\in G,$ $z\in M$. It is easy to see that
\begin{equation}
\ker\left( T_{g}^{\ast}\Phi_{z}\right) =\left( T_{g\cdot z}\left( G\cdot
z\right) \right) ^{\circ}=\left( {\rm Ver}_{g\cdot z}\right) ^{\circ}=
{\rm Hor}_{g\cdot z}^{\ast}.   \label{ggdot 8}
\end{equation}
Where $G \cdot z  $ denotes the $G$-orbit that contains the point $z \in M $.
Therefore, the map
\begin{equation}
\label{tilde map important}
\widetilde{T_{g}^{\ast}\Phi_{z}}:=\left. T_{g}^{\ast}\Phi_{z}\right\vert _{%
{\rm Ver}_{g\cdot z}^{\ast}}:T_{g\cdot z}^{\ast}M\cap {\rm Ver}_{g\cdot
z}^{\ast}\longrightarrow T_{g}^{\ast}G 
\end{equation}
is an isomorphism. Let%
\begin{equation*}
\begin{array}{rcl}
\rho\left( g,z\right) :T_{g}^{\ast}G & \longrightarrow & T_{g\cdot z}^{\ast
}M\cap{\rm Ver}_{g\cdot z}^{\ast}\subset T_{g\cdot z}^{\ast}M
\\ 
\alpha_{g} & \longmapsto & \left( \widetilde{T_{g}^{\ast}\Phi_{z}}\right)
^{-1}\left( \alpha_{g}\right)%
\end{array}
\end{equation*}
and define $\psi^{\ast}\left( x,z,g\right) :T_{g}^{\ast}G\rightarrow
T_{x}^{\ast}N$ by
\begin{equation*}
\psi^{\ast}\left( x,z,g\right) =S^{\ast}\left(x,g\cdot z\right) \circ
\rho\left( g,z\right) . 
\end{equation*}
Finally, we define a dual Stratonovich operator between the manifolds $G$
and $M\times N$ as%
\begin{equation}
\label{Stratonovich operator for group}
\begin{array}{rcl}
K^{\ast}\left( \left( z,x\right) ,g\right) :T_{g}^{\ast}G & \longrightarrow
& T_{z}^{\ast}M\times T_{x}^{\ast}N \\ 
\alpha_{g} & \longmapsto & \left( 0,\psi^{\ast}\left(x,z,g\right) \left(
\alpha_{g}\right) \right) .%
\end{array}
\end{equation}

\begin{theorem}[Reconstruction Theorem]
\label{teorema reconstruccion 1}
Let $X: \mathbb{R}_+ \times \Omega \rightarrow  N$ be a $N$-valued semimartingale
and let $S:TN \times  M \rightarrow TM $ be a Stratonovich operator that is
invariant with respect to a free and proper action of the Lie group $G$ on the
manifold $M$. If we are given $\Gamma_{M/G}$ a solution semimartingale of the reduced
stochastic differential equation then $\Gamma=g^{\Xi}\cdot d$ is a solution of the
original stochastic differential equation such that  $\pi(\Gamma)= \Gamma_{M/G} $. 

In this statement, $d: \mathbb{R}_+ \times  \Omega \rightarrow M$ is the horizontal
lift of
$\Gamma_{M/G}$ using an auxiliary principal connection on $\pi:M \rightarrow  M/G $
such that  $\Gamma_0=d_0$, and 
$g^{\Xi}:\mathbb{R}_{+}\times\Omega
\rightarrow G$ is the semimartingale solution of the stochastic differential
equation%
\begin{equation}
\delta g^{\Xi}=K\left( \Xi,g\right) \delta\Xi   \label{ggdot 4}
\end{equation}
with initial condition $g_{0}^{\Xi}=e$, $K$ the Stratonovich operator introduced
in~(\ref{Stratonovich operator for group}), 
 and stochastic component
$\Xi=\left( X,d\right) $ We will refer to $d$ as the {\bfi  horizontal
lift} of
$\Gamma_{M/G}$ and to $\Gamma=g^{\Xi} $ as the {\bfi  stochastic phase} of the
reconstructed solution. 
\end{theorem}

\begin{remark}
\normalfont
As we already pointed out, Theorem~\ref{teorema reconstruccion 1} is also valid when
the group action is not free. In that situation, one is given a solution of the
reduced stochastic differential equation on the quotient $M _I/L _I $, with $I$ an
isotropy subgroup of the $G$-action on  $M$. The correct statement
(and the proof that follows) of the reconstruction theorem in this case can be
obtained from the one that we just gave by replacing $M $ by the isotropy type
manifold $M _I $ and
$G $ by the group $L _I$. 
\end{remark}

\medskip

\noindent\textbf{Proof of Theorem~\ref{teorema reconstruccion 1}.\ \ }
In order to check that $\Gamma=g^{\Xi}\cdot d$ is a solution of the
original stochastic differential equation we have to verify that for any $\alpha \in
\Omega(M) $,
\begin{equation}
\label{thing to verify phase}
\int\langle \alpha, \delta \Gamma\rangle=\int\langle S ^\ast (X, \Gamma)\alpha,
\delta X\rangle.
\end{equation}
Since $\Gamma=g^{\Xi}\cdot d= \Phi(g^{\Xi},d)$, the statement in~\cite[Lemma
3.4]{Shigekawa} allows us to write%
\begin{equation}
\int\left\langle \alpha,\delta\Gamma\right\rangle =\int\left\langle
\Phi_{g^{\Xi}}^{\ast}\alpha,\delta d\right\rangle +\int\left\langle \Phi
_{d}^{\ast}\alpha,\delta g^{\Xi}\right\rangle .   \label{ggdot 6}
\end{equation}
We split the verification of~(\ref{thing to verify phase}) into two
cases:
\begin{description}
\item[(i)] $\alpha\in\Omega\left( M\right) $ is horizontal or, equivalently, 
$\alpha=\pi^{\ast}\left( \eta\right) $ with $\eta\in\Omega\left( M/G\right) $%
. Since $\alpha$ is horizontal, then $\Phi_{d}^{\ast}\alpha=0$ by (\ref{ggdot
8}). Then, using (\ref{ggdot 6}),%
\begin{equation*}
\!\!\!\!\!\!\!\!\!\!\!\!\!\!\!\!\!\!\!\!\int\left\langle
\alpha,\delta\Gamma\right\rangle =\int\left\langle
\Phi_{g^{\Xi}}^{\ast}\alpha,\delta d\right\rangle 
=\int\left\langle
\Phi_{g^{\Xi}}^{\ast}\left(\pi^{\ast}\left( \eta\right)\right) ,\delta
d\right\rangle =
\int\left\langle (\pi\circ \Phi_{g^{\Xi}})^{\ast}\left( \eta\right) ,\delta
d\right\rangle =
\int\left\langle \pi^{\ast}\left( \eta\right) ,\delta
d\right\rangle =\int\left\langle \eta,\delta\Gamma_{M/G}\right\rangle . 
\end{equation*}
We recall that $\Gamma_{M/G}=\pi\left( d\right) $ is a solution of the
reduced system, that is,
\begin{equation*}
\int\left\langle \eta,\delta\Gamma_{M/G}\right\rangle =\int\left\langle
S_{M/G}^{\ast}\left( X,\Gamma_{M/G}\right) \left( \eta\right) ,\delta
X\right\rangle 
\end{equation*}
for any $\eta\in\Omega\left( M/G\right) $. This implies by (\ref{general
reduced Stratonovich operator 1}) that
\begin{equation*}
\int\left\langle \eta,\delta\Gamma_{M/G}\right\rangle =\int\left\langle
S_{M/G}^{\ast}\left( X,\Gamma_{M/G}\right) \left( \eta\right) ,\delta
X\right\rangle =\int\left\langle S^{\ast}\left( X,d\right) \left( \pi
^{\ast}(\eta)\right) ,\delta X\right\rangle . 
\end{equation*}
Now, due to the $G$-invariance of $S,$ we know that $S^{\ast}\left( x,g\cdot
z\right) =S^{\ast}\left( x,z\right) \circ T_{z}^{\ast}\Phi_{g}$, for any $
g\in G$, $x \in N $, $z \in  M $. Recall also that $T_{z}\Phi_{g}$ sends the
horizontal space
${\rm Hor}_{z}$ to ${\rm Hor}_{g\cdot z}$ and the vertical space ${\rm Ver}
_{z}$ to ${\rm Ver}_{g\cdot z}$. Moreover, since $\alpha$ is horizontal, 
$\Phi_{g}^{\ast}\alpha=\alpha$ for any $g\in G$. Therefore,%
\begin{align*}
\int\left\langle \eta,\delta\Gamma_{M/G}\right\rangle & =\int\left\langle
S^{\ast}\left( X,d\right) \left( \alpha\right) ,\delta X\right\rangle
=\int\left\langle S^{\ast}\left( X,d\right) \left( \Phi_{g^{\Xi}}^{\ast
}\alpha\right) ,\delta X\right\rangle \\
& =\int\left\langle S^{\ast}\left( X,g^{\Xi}\cdot d\right) \left(
\alpha\right) ,\delta X\right\rangle =\int\left\langle S^{\ast}\left(
X,\Gamma\right) \left( \alpha\right) ,\delta X\right\rangle
\end{align*}
and hence~(\ref{thing to verify phase}) holds.

\item[(ii)] $\alpha\in\Omega\left( M\right) $ is vertical.
Since $\alpha$ is vertical, so is $\Phi_{g^{\Xi}}^{\ast}\alpha$ as a $
T^{\ast }M$-valued process. Therefore, $\int\left\langle
\Phi_{g^{\Xi}}^{\ast}\alpha,\delta d\right\rangle =0$ by (\ref{ggdot 9}).
Thus, using (\ref{ggdot 6}),
\begin{equation*}
\int\left\langle \alpha,\delta\Gamma\right\rangle =\int\left\langle \Phi
_{d}^{\ast}\alpha,\delta g^{\Xi}\right\rangle . 
\end{equation*}
Now, as $g^{\Xi}$ is a solution of the stochastic differential equation (\ref%
{ggdot 4}),%
\begin{align}
\int\left\langle \Phi_{d}^{\ast}\alpha,\delta g^{\Xi}\right\rangle &
=\int\left\langle K^{\ast}\left( \Xi,g^{\Xi}\right) \left( \Phi_{d}^{\ast
}\alpha\right) ,\delta\Xi\right\rangle =\int\left\langle \left( 0,\psi
^{\ast}\left( g^{\Xi},X,d\right) \left( \Phi_{d}^{\ast}\alpha\right) \right)
,\delta\Xi\right\rangle  \notag \\
& =\int\left\langle \psi^{\ast}\left( g^{\Xi},X,d\right) \left( \Phi
_{d}^{\ast}\alpha\right) ,\delta X\right\rangle .   \label{ggdot 10}
\end{align}
Recall that $\psi^{\ast}\left( x,z,g\right) =S^{\ast}\left( x,g\cdot
z\right) \circ\rho\left( g,z\right) $. Moreover $\rho\left( g,z\right)
\left( \gamma_{g}\right) =\left( \widetilde{T_{g}^{\ast}\Phi_{z}}\right)
^{-1}\left( \gamma_{g}\right) $ for any $\gamma_{g}\in T_{g}^{\ast}G$. Hence,%
\begin{equation*}
\rho\left( g,z\right) \circ T_{g}^{\ast}\Phi_{z}\left( \alpha_{g\cdot
z}\right) =\left( \widetilde{T_{g}^{\ast}\Phi_{z}}\right) ^{-1}\left(
T_{g}^{\ast}\Phi_{z}\left( \alpha_{g\cdot z}\right) \right) =\alpha_{g\cdot
z}
\end{equation*}
for any $\alpha_{g\cdot z}\in T_{g\cdot z}^{\ast}M\cap{\rm Ver}_{g\cdot
z}^{\ast}$, since in that situation $T_{g}^{\ast}\Phi_{z}\left( \alpha_{g\cdot
z}\right) =\widetilde{T_{g}^{\ast}\Phi_{z}}\left( \alpha_{g\cdot z}\right) $%
. Therefore, expression~(\ref{ggdot 10}) equals
\begin{equation*}
\int\left\langle \psi^{\ast}\left( g^{\Xi},X,d\right) \left(
\Phi_{d}^{\ast}\alpha\right) ,\delta X\right\rangle =\int\left\langle
S^{\ast}\left( X,g^{\Xi}\cdot d\right) \left( \alpha\right) ,\delta
X\right\rangle =\int\left\langle S^{\ast}\left( X,\Gamma\right) \left(
\alpha\right) ,\delta X\right\rangle,
\end{equation*}
and hence~(\ref{thing to verify phase}) also holds whenever 
$\alpha\in\Omega\left( M\right) $ is vertical, as required. \quad $\blacksquare$
\end{description}

\medskip

The stochastic phase $g^{\Xi }$ introduced in the Reconstruction Theorem admits
another characterization that we present in the paragraphs that follow.
Let $\left\{
\xi_{1},...,\xi_{m}\right\} $ be a basis of $\mathfrak{g}$, the Lie algebra
of $G$ and write
$ A=\sum_{i=1}^{m}A^{i}\xi_{i}$, where $A^{i}\in\Omega\left( M\right) $ are the
components of the auxiliary connection $A\in\Omega^{1}\left( M;\mathfrak{g}\right)
$ in this basis. Consider the $\mathfrak{g}$-valued semimartingale
\begin{equation}
Y=\sum_{i=1}^m\int\left\langle S^{\ast}\left( X,d\right) \left( A^{i}\right)
,\delta X\right\rangle \xi_i  \label{ggdot 15}.
\end{equation}

\begin{proposition}
\label{prop ggdot 1}Let $Y:\mathbb{R}_{+}\times \Omega \rightarrow \mathfrak{g}$ be
the $\mathfrak{g}$-valued semimartingale defined in~(\ref{ggdot 15}). Then, the
stochastic phase
$g^{\Xi }:\mathbb{R} _{+}\times
\Omega
\rightarrow G$ introduced in~(\ref{ggdot 4})  is the unique
solution of the stochastic differential equation
\begin{equation}
\delta g =L\left( Y,g \right) \delta Y  \label{ggdot 11}
\end{equation}
associated to the Stratonovich operator $L$ given by
\begin{equation*}
\begin{array}{rcl}
L\left( \xi,g\right) :T_{\xi}\mathfrak{g}& \longrightarrow  & T_{g}G \\ 
\eta & \longmapsto  & 
T_{e}L_{g}\left( \eta\right),
\end{array}
\end{equation*}
with initial condition $g_{0} =e$. The symbol $L _g:G \rightarrow G $ denotes the
left translation map by $g \in G $.
\end{proposition}

In the proof of this proposition, we will denote by $\xi_M (z):=
\left.\frac{d}{dt}\right|_{t=0}\exp t \xi \cdot z $ the infinitesimal vector field
associated to $\xi \in  \mathfrak{g}$  by  the $G$-action on $M$ evaluated at $z\in
M$. Analogously, we will write $\xi_G $ for the infinitesimal generators of the
$G$-action on itself by left translations. We recall (see~\cite{hsr} for a proof)
that for any $g\in G$, $\xi\in\mathfrak{g}$, and $z\in M$,
\begin{equation}
T_{z}\Phi_{g}\left(  \xi  _{M}\left( z\right) \right) =\left( 
{\rm Ad}_{g}\xi\right) _{M}\left( g\cdot z\right) . 
\label{ggdot 13}
\end{equation}
Moreover, $T_{g}\Phi _{z}\left(   \xi   _{G}\left( g\right) \right)
=T_{z}\Phi _{g}\left( \xi _{M}\left( z\right) \right) $ or,
in other words,%
\begin{equation}
 \xi  _{G}\left( g\right) =\widetilde{T_{g}\Phi _{z}}^{-1}\circ
T_{z}\Phi _{g}\left(  \xi _{M}\left( z\right) \right) ,
\label{ggdot 12}
\end{equation}%
where $\widetilde{T_{g}\Phi _{z}}^{-1}:T_{g\cdot z}M\cap {\rm Ver}
_{g\cdot z}\rightarrow T_{g}G$ is the isomorphism introduced in~(\ref{tilde map
important}). 

\medskip

\noindent \textbf{Proof of Proposition \ref{prop ggdot 1}}.
A result in~\cite{Shigekawa} shows that in order to prove the statement it suffices
to check that $\int \left\langle
\theta ,\delta g^{\Xi }\right\rangle =Y$, where $
\theta $ is the canonical $\mathfrak{g}$-valued one form on $G$ defined by $
\theta _{g}\left(  \xi _{G}\left( g\right) \right) =\xi
$, for any $g\in G$ and $\xi\in \mathfrak{g}$. Indeed, Lemmas 3.2 and  
3.3 in~\cite{Shigekawa} show that a $G$-valued semimartingale $g^G$  is such that
$\int \left\langle \theta ,\delta g^{G}\right\rangle =Y
$ if and only if $g^{G}$ is a solution of~(\ref{ggdot 11}).
Now, suppose that $g^{\Xi }$ is a solution of~(\ref{ggdot 4}), 
\begin{equation*}
\int \left\langle \theta ,\delta g^{\Xi }\right\rangle  =\int \left\langle
\psi ^{\ast }\left( g^{\Xi },X,d\right) \left( \theta \right) ,\delta
X\right\rangle =\int \left\langle S^{\ast }\left( X,g^{\Xi }\cdot d\right) \circ \rho
\left( g^{\Xi },d\right) \left( \theta \right) ,\delta X\right\rangle .
\end{equation*}
We are now going to verify that for any $g \in  G $ and $z \in M $,
\begin{equation}
\rho \left( g,z\right) \left( \theta \right) =\left( \Phi _{g^{-1}}^{\ast
}A\right) \left( g\cdot z\right) .  \label{ggdot 14}
\end{equation}
First of all notice that
as $\rho \left( g,z\right) \left( \gamma _{g}\right) =\left( \widetilde{
T_{g}^{\ast }\Phi _{z}}\right) ^{-1}\left( \gamma _{g}\right) \in T_{g\cdot
z}^{\ast }M\cap {\rm Ver}_{g\cdot z}^{\ast }$, for any $\gamma
_{g}\in T_{g}^{\ast }G$ and since $A$ vanishes when acting
on horizontal vector fields,
it suffices to verify (\ref{ggdot 14})
when acting on vector fields of the form $ \xi  _{M}$, for some $%
\xi \in \mathfrak{g}$. Using~(\ref{ggdot 13}), the right hand side of (\ref{ggdot
14}) then reads%
\begin{equation*}
\left( \Phi _{g^{-1}}^{\ast }A\right) \left( g\cdot z\right) \left( 
\xi  _{M}\left( g\cdot z\right) \right)  =A\left( z\right) \left(
T_{z}\Phi _{g^{-1}}\left(  \xi _{M}\left( g\cdot z\right)
\right) \right) 
 =A\left( z\right) \left( \left( {\rm Ad}
_{g^{-1}}\xi \right) _{M}\left( z\right) \right) 
={\rm Ad}_{g^{-1}}\xi .
\end{equation*}
As to the left hand side, we can write using~(\ref{ggdot 13})  and~(\ref{ggdot 12}), 
\begin{align*}
\rho \left( g,z\right) \left( \theta(g) \right) \left(   \xi 
_{M}\left( g\cdot z\right) \right) & =\left[ \left( \widetilde{T_{g}^{\ast
}\Phi _{z}}\right) ^{-1} \theta (g)  \right] \left( \xi
  _{M}\left( g\cdot z\right) \right)=\theta \left( g\right) \left[
\widetilde{T_{g}\Phi _{z}}^{-1}\left(
 \xi _{M}\left( g\cdot z\right) \right) \right]  \\
& =\theta \left( g\right) \left[ \widetilde{T_{g}\Phi _{z}}^{-1}\circ
T_{z}\Phi _{g}\circ T_{g\cdot z}\Phi _{g^{-1}}\left(  \xi 
_{M}\left( g\cdot z\right) \right) \right]  \\
 & =\theta \left( g\right) \left[ \widetilde{%
T_{g}\Phi _{z}}^{-1}\circ T_{z}\Phi _{g}\circ \left( {\rm Ad}%
_{g^{-1}}\xi \right) _{M}\left( z\right) \right]  \\
 & =\theta \left( g\right) \left[ \left( {\rm 
Ad}_{g^{-1}}\xi \right) _{G}\left( g\right) \right] ={\rm Ad}_{g^{-1}}\xi .
\end{align*}
Thus,
\begin{equation*}
\int \left\langle S^{\ast }\left( X,g^{\Xi }\cdot d\right) \circ \rho \left(
g^{\Xi },d\right) \left( \theta \right) ,\delta X\right\rangle =\int
\left\langle S^{\ast }\left( X,g^{\Xi }\cdot d\right) \left( \Phi _{\left(
g^{\Xi }\right) ^{-1}}^{\ast }A\right) ,\delta X\right\rangle .
\end{equation*}%
Now, since the Stratonovich operator $S$ is $G$-invariant, we have that
$S^{\ast }\left( x,g\cdot z\right) =S^{\ast }\left( x,z\right) \circ T_{z}^{\ast
}\Phi _{g}$, for any $x \in N $, $z \in M $, and $g \in  G $, and hence
\begin{equation*}
S^{\ast }\left(x,g\cdot z\right) \left( \left( \Phi _{g^{-1}}^{\ast
}A\right) \left( g\cdot z\right) \right) =S^{\ast }\left(x,z\right) \circ
T_{z}^{\ast }\Phi _{g}\circ T_{g\cdot z}^{\ast }\Phi _{g^{-1}}\left( A\left(
z\right) \right)  
=S^{\ast }\left(x,z\right) \left( A\left( z\right) \right) .
\end{equation*}
Therefore,
\begin{equation*}
\int \left\langle \theta ,\delta g^{\Xi }\right\rangle =\int \left\langle
S^{\ast }\left( X,g^{\Xi }\cdot d\right) \left( \Phi _{\left( g^{\Xi
}\right) ^{-1}}^{\ast }A\right) ,\delta X\right\rangle =\int \left\langle
S^{\ast }\left( X,d\right) \left( A\right) ,\delta X\right\rangle =Y,
\end{equation*}
and consequently $g^{\Xi }$ solves (\ref{ggdot 11}). The argument that we just
gave can be easily reversed to prove that if $g^{\Xi }$ is a solution
of (\ref{ggdot 11}) then it is also a solution of (\ref{ggdot 4}). \quad
$\blacksquare $

\bigskip

\noindent The combination of the reduction and the reconstruction of the
solution semimartingales of a symmetric stochastic differential equation
can be seen as a method to split the problem of finding its solutions into
three simpler tasks which we summarize as follows:

\begin{description}
\item[Step 1:] Find a solution $\Gamma_{M/G}$ for the reduced stochastic differential
equation associated to the reduced Stratonovich operator $S_{M/G}$ on the
dimensionally smaller space
$M/G $.
\item[Step 2:] Take an auxiliary principal connection $A\in \Omega
^{1}\left( M;\mathfrak{g}\right) $ for the principal bundle $\pi:M \rightarrow 
M/G $ and a horizontally lifted semimartingale
$d:\mathbb{R} _{+}\times \Omega \rightarrow M$, that is $\int \left\langle A,\delta
d\right\rangle =0$,  such that $d_{0}=\Gamma
_{0}$ and 
$\pi
\left( d\right) =\Gamma _{M/G}$.
\item[Step 3:] Let $g ^\Xi : \mathbb{R}_+ \times \Omega \rightarrow G$ be the
solution semimartingale of the stochastic differential equation (\ref{ggdot 11}) on
$G$
\begin{equation*}
\delta g =L\left( Y,g \right) \delta Y
\end{equation*}
with initial condition $g_{0}=e$ a.s. and with noise semimartingale $Y=\int
\left\langle S^{\ast}\left( X,d\right) \left( A\right) ,\delta X\right\rangle $. The
solution of the original stochastic differential equation associated to the
Stratonovich operator $S$ with initial condition $\Gamma_{0}$ is then
$\Gamma=\Phi_{g^\Xi}\left( d\right) $.
\end{description}

\begin{remark}
\label{remark completeness reduction}
\normalfont
Theorem~\ref{teorema reconstruccion 1} has as a consequence that the maximal existence times $\zeta $ and $\zeta_{M/G  }$ of $\pi $-related solutions $\Gamma  $ and $\Gamma_{M/G}$ of the original symmetric and reduced systems, coincide. Indeed,  if we write $\Gamma_{t}=g_{t}\cdot d_{t}$, with $d_{t}$ a horizontal lift of $\Gamma_{M/G}$, then first, $d_{t}$ is defined up to the same (maybe
finite) explosion time $\zeta_{M/G}$ of $\Gamma_{M/G}$. Second, as the semimartingale $g_{t}$ is the solution of the left-invariant stochastic differential equation~(\ref{ggdot 11}) then it  is in principle stochastically complete (\cite[Chapter VII \S 6, Example (i) page 131]{elworthy}) if its stochastic forcing is. Since in our case, the stochastic component $Y$~(\ref{ggdot 15})  depends on $d_{t}$, we can conclude that $g_t $ is defined again on the stochastic interval  $[0,\zeta_{M/G} )$. We consequently conclude that the maximal existence
time of the solutions of the initial symmetric system $\left(  M,S,X,N\right)  $ coincides with that of the corresponding solutions of the reduced system $(M/G,S_{M/G},X,N)$. Notice that this in particular implies that if the reduced manifold $M/G $ is compact then all the solutions of the original symmetric system are defined for all time, even if $M$ is not compact.
\end{remark}

\section{Symmetries and skew-product decompositions}
\label{Symmetries and skew-product decompositions}

The skew-product decomposition of second order differential operators is a factorization technique that has been used in the stochastic processes literature in order to split the semielliptic and, in particular, the diffusion operators, associated to certain stochastic differential equations (see, for instance,~\cite{Pauwels-Rogers, liao skew-product, taylor skew-product}, and references therein). This splitting has important consequences as to the  properties of the solutions of these equations, like certain factorization properties of their probability laws and of the associated stochastic flows.

Symmetries are a natural way to obtain this kind of decompositions as it has already been exploited in~\cite{elworthy principal bundles}. Our goal in the following pages consists of generalizing the existing results in two ways: first, we will generalize the notion of skew-product to arbitrary stochastic differential equations by working with the notion of skew-product decomposition of the Stratonovich operator; we will indicate below how our approach coincides with the traditional one in the case of diffusions. Second, we will show that the skew-product decompositions presented in~\cite{elworthy principal bundles} for regular free action are also available (at least locally)  for singular proper group actions.

\begin{definition}
\label{def skew-product}
Let $N$, $M_{1}$, and $M_{2}$ be three smooth manifolds
and $S\left(  x,m\right)  :T_{x}N\rightarrow T_{m}\left(  M_{1}\times
M_{2}\right)  $, $x\in N$, $m=\left(  m_{1},m_{2}\right)  \in M_{1}\times
M_{2}$, a Stratonovich operator from $N$ to the product manifold $M_{1}\times
M_{2}$. We will say that $S$ admits a {\bfi  skew-product decomposition} if there exists a Stratonovich
operator $S_{2}\left(  x,m_{2}\right)  :T_{x}N\longrightarrow T_{m_{2}}M_{2}$
from $N$ to $M_{2}$ and a $M_{2}$-dependent Stratonovich operator
$S_{1}\left(  x,m_{1},m_{2}\right)  :T_{x}N\rightarrow T_{m_{1}}M_{1}$ such
that%
\[
S\left(  x,m\right)  =\left(  S_{1}\left(  x,m_{1},m_{2}\right)  ,S_{2}\left(
x,m_{2}\right)  \right)  \in\mathcal{L}\left(  T_{x}N,T_{m_{1}}M_{1}\times
T_{m_{2}}M_{2}\right)
\]
for any $m=\left(  m_{1},m_{2}\right)  \in M_{1}\times M_{2}$. The operators
$S_{1}$ and $S_{2}$ will be called the \textbf{factors} of $S$.
\end{definition}

In order to show the relation between this definition  and the classical one used in the papers that we just quoted, we first have to briefly recall the relation between the global Stratonovich and It\^o formulations for the stochastic differential equations (see
\cite{emery} for a detailed presentation of this subject). Given $M$ and $N$ two manifolds, a Schwartz operator is
a family of Schwartz maps (see \cite[Definition 6.22]{emery})
$\mathcal{ S}\left( x,z\right) :\tau _{x}N\rightarrow \tau _{z}M$ between the
tangent bundles of second order $\tau N$ and $\tau M$. In this context, the It\^{o}
stochastic differential equation defined by the Schwartz
operator $\mathcal{S}$ with stochastic component a continuous
semimartingale $ X:\mathbb{R}_{+}\times \Omega \rightarrow N$  is
\begin{equation}
d\Gamma =\mathcal{S}\left( X,\Gamma \right) dX.
\label{Ito differential form}
\end{equation}
Given a Stratonovich operator $S$, there is a unique Schwartz
operator $\mathcal{S}:\tau N\times M\rightarrow \tau M$ that is an extension
of $S$ to the tangent bundles of second order and
which
makes the It\^o and Stratonovich stochastic differential equations associated to $S
$ and
$\mathcal{S}$  equivalent, in the
sense that they have the same semimartingale solutions. $\mathcal{S} $ is
constructed as follows. Let $\gamma (t)=(x(t),z(t))\in N\times M$ be a smooth curve
that verifies
$S(x(t),z(t))(\dot{x}(t))=\dot{z}(t)$, for all $t$. We define
$\mathcal{S}(x(t),z(t))\left( L_{\ddot{x}(t)}\right) :=\left( L_{%
\ddot{z}(t)}\right) $, where the second order differential operators $\left(
L_{\ddot{x}(t)}\right) \in \tau _{x\left( t\right) }N$ and $\left( L_{\ddot{z%
}(t)}\right) \in \tau _{z\left( t\right) }M$ are defined as $\left( L_{\ddot{%
x}(t)}\right) \left[ h\right] :=\frac{d^{2}}{d^{2}t}h\left( x\left( t\right)
\right) $ and $\left( L_{\ddot{z}(t)}\right) \left[ g\right] :=\frac{d^{2}}{
d^{2}t}g\left( z\left( t\right) \right) $, for any $h\in C^{\infty }(N)$ and 
$g\in C^{\infty }\left( M\right) $. This relation completely determines $%
\mathcal{S}$ since the vectors of the form $L_{\ddot{x}(t)}$ span $\tau
_{x\left( t\right) }M$. 

It is easy to show that if $S:TN\times\left(
M_{1}\times M_{2}\right)  \rightarrow T\left(  M_{1}\times M_{2}\right)  $ is a
Stratonovich operator that admits a skew-product decomposition with factors $S_{1}$ and $S_{2}$ then the equivalent Schwartz operator
$\mathcal{S}:\tau N\times\left(  M_{1}\times M_{2}\right)  \rightarrow
\tau\left(  M_{1}\times M_{2}\right)  $  can be written as
\begin{equation}
\mathcal{S}\left(  x,\left(  m_{1},m_{2}\right)  \right)  = 
\mathcal{S}_{1}\left(  x,m_{1},m_{2}\right)+\mathcal{S}
_{2}\left(  x,m_{2}\right)  ,\label{eq skew 1}%
\end{equation}
for any $x\in N$ and any $m=\left(  m_{1},m_{2}\right)  \in M_{1}\times M_{2}%
$. In this expression, $\mathcal{S}_{1}\left(  x,m_{1},m_{2}\right) :\tau_{x}N\rightarrow\tau_{m}(M_1\times M _2) $ and 
$\mathcal{S}_{2}\left(  x,m_{2}\right)  :\tau_{x}N\rightarrow\tau_{m}(M_1\times M _2)$ are
the equivalent Schwartz operators of the Stratonovich operators $\widetilde{S}_{1}, \widetilde{S}_{2}:TN\times\left(M_{1}\times M_{2}\right)  \rightarrow T\left(  M_{1}\times M_{2}\right)  $  defined by  $\widetilde{S} _1(x,m):=T _{m _1}i_{m _2}\left( S _1(x,m _1, m _2)\right) $ and $\widetilde{S} _2(x,m):=T _{m _2}i_{m _1}\left( S _2(x,m _2)\right) $. The maps $i_{m _1}:M _2 \rightarrow M _1 \times M _2 $ and $i_{m _2}:M _1 \rightarrow M _1 \times M _2 $ are the natural inclusions obtained by fixing $m _1 $ and $m _2 $, respectively.

Now, the notion of  skew-product decomposition of a second
order differential operator $L\in\mathfrak{X}_{2}\left(  M_{1}\times
M_{2}\right)  $ on $M_{1}\times M_{2}$ that one finds in the literature (see for instance~\cite{taylor skew-product}) consists on the existence of two
smooth maps $L_{1}:M_{2}\rightarrow\mathfrak{X}_{2}\left(  M_{1}\right)  $ and
$L_{2}\in\mathfrak{X}_{2}\left(  M_{2}\right)  $ such that for any $f\in C^{\infty}\left(  M_{1}\times M_{2}\right)  $
\begin{equation}
\label{skew-product operators}
L\left[  f\right]  \left(  m_{1},m_{2}\right)  =\left(  L_{1}\left(
m_{2}\right)  \left[  f\left(  \cdot,m_{2}\right)  \right]  \right)  \left(
m_{1}\right)  +\left(  L_{2}\left[  f\left(  m_{1},\cdot\right)  \right]
\right)  \left(  m_{2}\right).
\end{equation}
The relation between this notion and the one introduced in Definition~\ref{def skew-product} is very easy to establish for semielliptic diffusions. Indeed, suppose that the Stratonovich operator associated to a semielliptic diffusion admits a skew-product decomposition;  we just saw that this implies in general the existence of a skew-product decomposition~(\ref{eq skew 1}) of the corresponding Schwartz operator, which in turn implies  the availability of a skew-product decomposition of the infinitesimal generator associated to~(\ref{Ito differential form}) in the sense of~(\ref{skew-product operators}). See~\cite[page 15]{taylor skew-product} for a sketch of the proof of this fact.

In conclusion, since in the cases that have already been studied, the skew-product decompositions of Stratonovich operators carry in their wake the skew-product decompositions as differential operators  of the associated infinitesimal generators, we can focus in what follows on the more general situation that consists of  adopting Definition~\ref{def skew-product}.

\subsection{Skew-products on principal fiber bundles. Free actions.}

Let $M$, $N$ be two manifolds, $G$ a Lie group, and $\Phi:G\times M\rightarrow
M$ a proper and free action. We already know that $M/G$ is a smooth manifold
under these hypotheses and that $\pi_{M/G}:M\rightarrow M/G$ is a principal fiber
bundle with structural group $G$. The goal of the following paragraphs is to
show that any $G$-invariant Stratonovich operator $S:TN\times M\rightarrow TM$
on $M$ admits a local skew-product decomposition. This result is also true even
if the action $\Phi$ is not free, as we will see in the next section. However, what makes this local decomposition possible in this simpler case   is not the fact that the $G$-action is free and proper but that
$\pi_{M/G}:M\rightarrow M/G$ is a principal fiber bundle.
Consequently, in order to keep our exposition as general
as possible, we will adopt as the setup for the rest of this subsection a
$G$-invariant Stratonovich operator $S:TN\times P\rightarrow TP$ on an
arbitrary (left) $G$-principal fiber bundle $\pi:P\rightarrow Q$. 
This setup has been studied in detail in~\cite{elworthy principal bundles} for invariant diffusions. In the following proposition we generalize the vertical-horizontal splitting in that paper to arbitrary Stratonovich operators  and we formulate it in terms of skew-products.

\begin{proposition}
\label{skew-product caso libre}Let $N$ be a manifold, $\pi:P\rightarrow Q$ a
(left) principal bundle with structure group $G$, $S:TN\times P\rightarrow TP$
a $G$-invariant Stratonovich operator, $X:\mathbb{R}_{+}\times\Omega
\rightarrow N$ a $N$-valued semimartingale, and $\sigma:U\rightarrow\pi
^{-1}\left(  U\right)  \subseteq P$ a local section of  $ \pi  $ defined on an open neighborhood
$U\subseteq Q$. Then, $S$ admits a skew-product decomposition on $\pi^{-1}\left(  U\right)  $.
More explicitly, there exists a diffeomorphism $F:G\times U\rightarrow\pi^{-1}\left(
U\right)  $ and a skew-product split Stratonovich operator $S_{G \times U}
:TN\times\left(  G\times U\right)  \rightarrow T\left(  G\times U\right)  $
such that $F$ establishes a bijection between semimartingales $\Gamma$
starting on $\pi^{-1}\left(  U\right)  $ which are solutions of the stochastic
system $\left(  P,S,X,N\right)  $ up to time $\tau=\inf\left\{  t>0~|~\Gamma
_{t}\notin \pi^{-1}(U)\right\}  $ and the $\left(  G\times U\right)  $-valued
semimartingales $(\widetilde{g}_{t},\Gamma_{t}^{Q})$ that solve $\left(  G\times U,S_{G \times U},X,N\right)  $,
\begin{equation}
\delta(\widetilde{g}_{t},\Gamma_{t}^{Q})=S_{G \times U}\left(  X,(\widetilde
{g}_{t},\Gamma_{t}^{Q})\right)  \delta X_{t}. \label{eq tan-nor 11}%
\end{equation}

\end{proposition}

\noindent\textbf{Proof.\ \ } Let $U\subseteq Q$ be an open neighborhood and $\sigma:U\rightarrow\pi^{-1}\left(  U\right)  \subseteq P$ a local section of $\pi:P\rightarrow Q$. Given that $G$ acts freely on $P$, the map
\[
\begin{array}
[c]{rrr}%
F:G\times U & \longrightarrow & \pi^{-1}\left(  U\right) \\
\left(  g,q\right)  & \longmapsto & g\cdot\sigma\left(  q\right)
\end{array}
\]
is a $G$-equivariant diffeomorphism, where $g\cdot\sigma\left(  q\right)  =\Phi_{g}\left(
\sigma\left(  q\right)  \right)  $ denotes the (left) action of $g\in
G$ on $\sigma\left(  q\right)  \in P$ via $\Phi:G\times P\rightarrow P$ and the product manifold $G \times  U $ is considered as a left $G$-space with the action defined by $g \cdot (h, q):=(g \cdot h, q) $. Thus,
we can use $F$ to identify $\pi^{-1}\left(  U\right)  \subseteq P$ with the
product manifold $G \times U$. 

Now, given $p=g\cdot\sigma\left(  q\right)  \in
\pi^{-1}\left(  U\right)  $, define $\operatorname*{Hor}_{p}\subseteq T_{p}P$ as  $\operatorname*{Hor}_{p}:=T_{\sigma\left(  q\right)  }\Phi_{g}\circ
T_{q}\sigma\left(  T_{q}Q\right)  $. It is straightforward to see that the family of
horizontal spaces $\{\operatorname*{Hor}_{p}~|~p\in\pi^{-1}\left(  U\right)
\}$ is invariant by the $G$-action and hence defines a principal
connection $A_{\sigma}\in\Omega^{1}\left(  \pi^{-1}\left(  U\right)
;\mathfrak{g}\right)  $ on the open neighborhood $\pi^{-1}\left(  U\right)  $.
Moreover, if $\Gamma^{Q}:\mathbb{R}_{+}\times\Omega\rightarrow Q$ is a
$Q$-valued semimartingale starting at $q$, then $\sigma\left(  \Gamma
^{Q}\right)  $ is the unique horizontal lift on $P$ of $\Gamma^{Q}$ associated
to the connection $A_{\sigma}$ starting at $\sigma\left(  q\right)  \in
\pi^{-1}\left(  q\right)  $ and defined up to time $\tau_{U}=\inf
\{t>0~|~\Gamma_{t}^{Q}\notin U\}$.

Consider now the skew-product split Stratonovich operator $S_{G \times U}\left(
x,\left(  g,q\right)  \right)  :TN\times\left(  G\times U\right)  \rightarrow
T\left(  G\times U\right)  $ such that, for any $x\in N$, $g\in G$, $q\in U$%
\[
S_{G \times U}\left(  x,\left(  g,q\right)  \right)  =\left(  K\left(
(\sigma(q),x\right)  ,g),S_{P/G}\left(  x,q\right)  \right)  \in
\mathcal{L}\left(  T_{x}N,T_{g}G\times T_{q}U\right)  ,
\]
where $K$ is the Stratonovich operator introduced in~(\ref{Stratonovich operator for group}) and $S_{P/G}$
the reduced Stratonovich operator constructed out of $S$ as in~(\ref{general reduced Stratonovich operator 1}). Let $(\widetilde{g}_{t}
,\Gamma_{t}^{Q})$ a $\left(  G\times U\right)  $-valued semimartingale
solution of the stochastic system (\ref{eq tan-nor 11}), i.e.
\[
\delta(\widetilde{g}_{t},\Gamma_{t}^{Q})=S_{G \times U}\left(  X,(\widetilde
{g}_{t},\Gamma_{t}^{Q})\right)  \delta X,
\]
with initial condition $\left(  g,q\right)  \in G\times U$. We claim that
$\Gamma_{t}=F(\widetilde{g}_{t},\Gamma_{t}^{Q})=\widetilde{g}_{t}\cdot
\sigma(\Gamma_{t}^{Q})$ is a solution of the stochastic system $\left(
P,S,X,N\right)  $ with initial condition $g\cdot\sigma(q)$ up to the first
exit time $\tau_{U}=\inf\{t>0~|~\Gamma_{t}^{Q}\notin U\}$. This is a
consequence of the Reconstruction Theorem~\ref{teorema reconstruccion 1} and the fact that $\sigma
(\Gamma_{t}^{Q})$ is the horizontal lift of a solution of the reduced system
$\left(  Q,S_{P/G},X,N\right)  $. Conversely, let $\Gamma$ be a solution of the stochastic system $\left(  P,S,X,N\right)  $ with initial condition $p=g\cdot\sigma\left(  q\right)
\in\pi^{-1}\left(  U\right)  $. By the Reconstruction Theorem~\ref{teorema reconstruccion 1}, $\Gamma$ can be written as $\Gamma_{t}=\widetilde{g}_{t}\cdot d_{t}$. We recall that $d_{t}$ the horizontal lift with respect to an
arbitrary connection $A\in\Omega^{1}\left(  Q;\mathfrak{g}\right)  $ of the
solution $\Gamma_{t}^{Q}= \pi(\Gamma _t)$ of the reduced system $\left(  Q,S_{P/G},X,N\right)
$ (see Theorem~\ref{reduction theorem statement general systems}) with initial condition $\sigma(q) $. On the other hand, $\widetilde{g}_{t}$ is the solution of the stochastic
system~(\ref{ggdot 4}) with initial condition $g\in G$. If we take in this procedure $A_{\sigma}\in\Omega^{1}\left(  \pi^{-1}(U);\mathfrak{g}\right)  $ as the auxiliary connection, that is, the one given by the local section $\sigma:U\rightarrow\pi^{-1}\left(  U\right)  $, then   $d_{t}=\sigma(\Gamma_{t}^{Q})$ and  it is straightforward to
check that $(\widetilde{g}_{t},\Gamma_{t}^{Q})$ is a solution of
(\ref{eq tan-nor 11}) with initial condition $\left(  g,q\right)  \in G\times
U$. $\blacksquare$

\begin{example}
\normalfont
Let $G$ be a Lie group, $H\subseteq G$ a closed subgroup, and $R$ a
smooth manifold. In \cite{Pauwels-Rogers}, Pauwels and Rogers show several
examples of skew-product decompositions of Brownian motions on manifolds of the
type $R\times G/H$ which share a common feature, namely, they are obtained from
skew-product split Brownian motions on $R\times G$ via the reduction $\pi:R\times
G\rightarrow R\times G/H$. The $H$-action on $R\times G$ is 
$h\cdot\left(  r,g\right)  :=\left(  r,gh\right)  $, for any $h\in H$, $r\in
R$, and $g\in G$. An important result in this paper is Theorem 2 which reads as
follows: suppose that $R\times G/H$ is a Riemannian manifold with Riemannian
metric $\eta$ and that the tensor $\pi^{\ast}\eta$ is $G$-invariant. Furthermore, suppose that
the decomposition $T_{\left(  r,g\right)  }\left(  R\times G\right)
=T_{r}R\oplus T_{g}G$ is orthogonal with respect to $\pi^{\ast}\eta$, for any
$r\in R$, $g\in G$, and that the Lie algebra $\mathfrak{g}$ of $G$ admits an
$\operatorname*{Ad}\nolimits_{H}$-invariant inner product. Under these
hypotheses,  $R\times G$ admits a $G$-invariant
Riemannian metric $\hat{\eta}$ such that  if $\Gamma$ is a Brownian motion on
$R\times G$ with respect to $\hat{\eta}$ then $\Gamma$ has a skew-product decomposition and moreover, $\pi\left(
\Gamma\right)  $ is a Brownian motion on $\left(  R\times G/H,\eta\right)  $.
This result is repeatedly used in \cite{Pauwels-Rogers} to obtain skew-product
decompositions of Brownian motions on various manifolds of matrices.
\end{example}

\begin{example}[Brownian motion on symmetric spaces]
\normalfont
Let $\left(
M,\eta\right)  $ be a  Riemannian symmetric space with Riemannian metric $\eta$.
We want to define Brownian motions on   $\left(
M,\eta\right)  $ by reducing   a suitable process defined on the connected component  containing the identity of its group of isometries. The notation and most of the results in
this example, in addition to a comprehensive exposition on symmetric spaces,
can be found in \cite{helgason} and \cite{kobayashi-nomizu}. The reader is
encouraged to check with \cite{elworthy diffusion book} to learn more about stochastics in the context of homogeneous spaces. 

We start by recalling that a $M$-valued process $\Gamma$ is a Brownian motion whenever
\[
f(\Gamma)-f\left(  \Gamma_{0}\right)  -\frac{1}{2}\int\Delta(f)\left(
\Gamma_{s}\right)  ds
\]
is a real valued  local semimartingale for any $f\in C^{\infty}(M)$, where $\Delta$
denotes the Laplacian. The Laplacian is defined as the
trace of the Hessian associated to the Levi-Civita connection $\nabla$ of
$\eta$, that is,
\[
\Delta\left(  f\right)  \left(  m\right)=\sum_{i=1}^{r}\left(
\mathcal{L}_{Y_{i}}\circ\mathcal{L}_{Y_{i}}-\nabla_{Y_{i}}Y_{i}\right)
(f)(m)
\]
where $\left\{  Y_{1},...,Y_{r}\right\}  \subset\mathfrak{X}\left(  M\right)
$ is family or vector fields such that $\left\{  Y_{1}(m),...,Y_{r}%
(m)\right\}  $ is an orthonormal basis of $T_{m}M$, $m\in M$.

Let $G$ be the connected component containing the identity of the
isometries group $I(M)\subseteq\operatorname{Diff}(M)$ of $M$. Take $o\in M$ a
fixed point and let $s$ be a geodesic symmetry at $o$. The
Lie group $G$ acts on $M$ transitively and, if $K$ denotes the isotropy group
of $o$, $M$ is diffeomorphic to $G/K$ (\cite[Chapter IV, Theorem 3.3]%
{helgason}). Denote by $\pi:G\rightarrow G/K$ the canonical projection and
suppose that $\dim\left(  G\right)  <\infty$. Let $\sigma:G\rightarrow G$ be
the involutive automorphism of $G$ defined by $\sigma\left(  g\right)
=s\circ\Phi_{g}\circ s$ for any $g\in G$, where $\Phi:G\times M\rightarrow G$
denotes as usual the left action of $G$ on $M$. $T_{e}\sigma:\mathfrak{g}%
\rightarrow\mathfrak{g}$ induces an involutive automorphism of $\mathfrak{g}$.
That is, $T_{e}\sigma\circ T_{e}\sigma=\operatorname*{Id}$ but $T_{e}%
\sigma\neq\operatorname*{Id}$. Let $\mathfrak{k}$ and $\mathfrak{m}$ be the
the eigenspaces in $\mathfrak{g}$ associated to the eigenvalues $1$ and $-1$
of $T_{e}\sigma$, respectively, such that $\mathfrak{g}=\mathfrak{k}%
\oplus\mathfrak{m}$. It can be checked that $\mathfrak{k}$ is a Lie subalgebra
of $\mathfrak{g}$ and that (see~\cite[Chapter XI Proposition 2.1]{kobayashi-nomizu}).
\[
\left[  \mathfrak{k},\mathfrak{k}\right]  \subseteq\mathfrak{k}\text{,
\ \ }\left[  \mathfrak{k},\mathfrak{m}\right]  \subseteq\mathfrak{m}\text{,
\ and \ }\left[  \mathfrak{m},\mathfrak{m}\right]  \subseteq\mathfrak{k}.
\]
Since the
infinitesimal generators $\xi^{M}\in\mathfrak{X}\left(  M\right)  $ of the
$G$-action $\Phi$ on $M$, with $\xi\in\mathfrak{m}$, span the tangent space at any point
$gK\in G/K$, any affine connection  is fully characterized by its value  on the left-invariant
vector fields $\xi^{M}$ with $\xi\in\mathfrak{m}$. In the particular case of the Levi-Civita connection $\nabla$ associated to the metric $\eta$, its $G$-invariance  implies via~\cite[Chapter XI, Theorem 3.3]{kobayashi-nomizu} that
\begin{equation}
\nabla_{\xi^{M}}\zeta^{M}\left(  gK\right)  =0\label{eq tan-nor 26}%
\end{equation}
for any pair of left-invariant vector fields $\xi^{M}$ and $\zeta^{M}$. A
consequence of (\ref{eq tan-nor 26}) is that the Laplacian $\Delta$
takes the expression   $\Delta\left(  f\right)  \left(
gK\right)  =\sum_{i=1}^{r}\mathcal{L}_{\xi_{i}^{M}}\circ\mathcal{L}_{\xi
_{i}^{M}}(f)(gK)$,  $gK\in G/K$, where $\left\{  \xi_{1}^{M}(gK),...,\xi_{r}^{M}
(gK)\right\}  $ is an orthonormal basis of $T_{gK}(G/K)$.

Let $\left\{  \xi_{1},...,\xi_{r}\right\}  $ be a basis of $\mathfrak{m}$ such
that $\left\{  T_{e}\pi\left(  \xi_{1}\right)  ...,T_{e}\pi\left(  \xi
_{r}\right)  \right\}  $ is an orthonormal basis of $T_{K}(G/K) \simeq T_o M$ with respect
to $\eta_{o}$ and let $\left\{  \xi_{1}^{G},...,\xi_{r}^{G}\right\}
\subset\mathfrak{X}\left(  G\right)  $ the corresponding family of
\textit{right-invariant} vector fields built from $\left\{  \xi_{1}%
,...,\xi_{r}\right\}  $. Observe that $\left\{  \xi_{1}^{M},...,\xi_{r}%
^{M}\right\}  $ is an orthonormal basis of the tangent space at any point
$gK\in G/K$ due to the transitivity of the $G$-action on $M$ and to the $G$-invariance of
the metric $\eta$. Consider now the Stratonovich stochastic differential
equation on $G$
\begin{equation}
\delta g_{t}=\sum_{i=1}^{r}\xi_{i}^{G}(g_{t})\delta B_{t}^{i}%
,\label{eq tan-nor 24}%
\end{equation}
where $\left(  B_{t}^{1},...,B_{t}^{r}\right)  $ is a $\mathbb{R}^{r}$-valued
Brownian motion. The equation (\ref{eq tan-nor 24}) is by construction
$K$-invariant with respect to the \textit{right action} $R:K\times
G\rightarrow G$, $R_{k}\left(  g\right)  =gk$.
In addition, it is straightforward to check that the projection $\pi
:G\rightarrow G/K$ send any right-invariant vector field $\xi^{G}%
\in\mathfrak{X}\left(  G\right)  $, $\xi\in\mathfrak{g}$, to the infinitesimal
generator $\xi^{M}\in\mathfrak{X}\left(  M\right)  $ of the $G$-action
$\Phi:G\times M\rightarrow M$. Indeed,
for any $\xi\in\mathfrak{g}$, $g\in G$, and $k\in K$\[
T_{g}\pi\left(  \xi^{G}(g)\right)  =T_{g}\pi\circ T_{e}R_{g}\left(
\xi\right)  =\left.  \frac{d}{dt}\right\vert _{t=0}\pi\left(  \exp\left(
t\xi\right)  g\right)  =\left.  \frac{d}{dt}\right\vert _{t=0}\Phi\left(
\exp\left(  t\xi\right)  ,\pi\left(  g\right)  \right)  =\xi^{M}(gK),
\]
and hence (\ref{eq tan-nor 24}) projects to the stochastic differential equation
\begin{equation}
\delta\Gamma_{t}=\sum_{i=1}^{r}\xi_{i}^{M}(\Gamma_{t})\delta B_{t}%
^{i}\label{eq tan-nor 25}%
\end{equation}
on $M$ by the Reduction Theorem~\ref{reduction theorem statement general systems}. A straightforward computation shows that that the solution semimartingales of
(\ref{eq tan-nor 25}) have as infinitesimal generator the Laplacian
$\Delta=\sum_{i=1}^{r}\mathcal{L}_{\xi_{i}^{M}}\circ\mathcal{L}_{\xi
_{i}^{M}}$ and hence by the It\^o formula
\[
f(\Gamma)-f\left(  \Gamma_{0}\right)  -\frac{1}{2}\int\Delta(f)\left(
\Gamma_{s}\right)  ds=\sum_{i=1}^r\int  \xi_i^M [f](\Gamma) dB ^i
\]
which allows us to conclude that they are Brownian motions. It
is worth noticing that since right-invariant systems such that
(\ref{eq tan-nor 24}) are stochastically complete (see \cite[Chapter VII
\S 6]{elworthy}) and by the Reduction and Reconstruction Theorems~\ref{reduction theorem statement general systems} and~\ref{teorema reconstruccion 1} any solution of (\ref{eq tan-nor 25}) may be written
as $\Gamma_{t}=\pi\left(  g_{t}\right)  $ for a suitable solution $g_{t}$ of
(\ref{eq tan-nor 24}), the Brownian motion on a symmetric space is
stochastically complete.
\end{example}

\subsection{Skew-products induced by non-free actions. The tangent-normal decomposition}

In this section we will show how the results that we just presented for free actions can be generalized to the non-free case by using the notion of {\bfi  slice}~\cite{koszul53, pa} and a generalization to the context of Stratonovich operators of the so-called {\bfi  tangent-normal decomposition} of $G$-equivariant vector fields with respect to proper group actions~\cite{krupa, field91}. 

Let $\Phi:G\times M\rightarrow M$ be a proper action of the Lie group $G$ on the manifold $M$ and let
$M/G$ be the associated orbit space, $M/G$. Observe that as the group action is not necessarily free, the orbit space $M/G$ needs not be a smooth
manifold. 

In order to introduce the notion of slice we start by considering a subgroup
$H\subset G$  of $G$.
Suppose that $H$ acts on the left on a certain manifold $A$. The {\bfi
twisted action\/}
of $H$ on the product $G\times A$ is defined by
\[h\cdot (g,\,a)=(gh,\,h^{-1}\cdot a), \quad h \in H,\; g \in G,\text{ and } a \in A.\]
Note  that this action is free and proper by the freeness and properness of the 
action on the $G$-factor. The 
{\bfi twisted product\/}
$G\times_{H} A$ is defined as the orbit space $(G \times A)/H $ corresponding to the
twisted action. The elements of $G\times_{H} A$ will be denoted by $[g,a]$, $g  \in  G $,
$a\in A $. The twisted product $G\times_H A$ is a
$G$-space relative to the left action defined by 
$g'\cdot [g,a] := [g'g, a]$. Also, it can be shown that the action of $H$ on $A$ is proper if and only if the $G$-action on $G\times_H A$ just defined is proper (see~\cite[Proposition 2.3.17]{hsr}).

Let now $m\in M$ and denote $H:=G_m$. A {\bfi tube\/}
around the orbit
$G\cdot m$ is a $G$--equivariant diffeomorphism
\[\varphi:G\times_{H} A\longrightarrow U,\]
where $U$ is a $G$--invariant neighborhood of the orbit $G\cdot m$ and $A$ 
is some manifold on which $H$ acts.  Note that the $G$--action on the twisted
product 
$G\times_{H} A$ is proper since by the properness of the $G$-action on $M$,
the isotropy subgroup $H$ is
compact and, consequently, its action on $A$ is proper.

\begin{definition}
\label{slice}
Let $M$ be a manifold and $G$  a Lie group acting properly on 
$M$. Let $m\in M$ and denote $G_m:=H$. Let $W$ be a 
submanifold of $M$ such that $m\in W$ and $H\cdot W =W$. We say
that $W$ is a {\bfi slice\/}
\index{slice}
at $m$ if the $G$--equivariant map
\[\begin{array}{cccc}
\varphi: &G\times_{H} W&\longrightarrow&U\\
 &\left[ g,\,s\right]& \longmapsto &g\cdot s
\end{array}\]
is a tube about $G\cdot m$ for some $G$--invariant open 
neighborhood $U$ of $G\cdot m$. Notice that if $W$ is a slice at $m$ then $\Phi_g(W) $ is a
slice at the point $\Phi_g(m)$.
\end{definition}

The Slice Theorem of Palais~\cite{pa} proves that there exists a slice at any point of a proper $G$-manifold. The following theorem, whose proof can be found in~\cite{hsr} provides several equivalent characterizations of the concept of
slice that are available in the literature.
\begin{theorem}
\label{slice theorem}
Let $M$ be a manifold and $G$  a Lie group acting properly 
on $M$. Let $m\in M$, denote $H:=G _m$, $\mathfrak{h}$ the Lie algebra of $H$, and let
$W$ be a submanifold of $M$ containing $m$. Then the following
statements are equivalent:
\begin{description}
\item[(i)] There is a tube $\varphi:G\times_{H} A\longrightarrow 
U$ about $G\cdot m$ such that $\varphi [e,\,A]=W$. 
\item[(ii)] $W$ is a slice at $m$.
\item [(iii)] The submanifold $W$ satisfies the following properties:
\begin{description}
\item [{\bf (a)}] The set $G \cdot W $ is an open neighborhood of the orbit $G \cdot m $
and $ W $ is closed in $G \cdot W $.
\item [{\bf (b)}] For any $z \in W $ we have that $T _z M = \mathfrak{g}\cdot z + T _zW $.
Moreover, $ \mathfrak{g}\cdot z \cap T _zW=\mathfrak{h}\cdot z $. In particular, for
$z=m$ the sum $\mathfrak{g}\cdot z + T _zW $ is direct.
\item [{\bf (c)}] $W$ is $H$--invariant. Moreover, if $z \in W $ and $g \in G $ are such
that $g \cdot z \in W$, then $g \in H $.
\item [{\bf (d)}] Let $\sigma: U \subset G/H \rightarrow G $ be a local section of the
submersion $G \rightarrow G/H $. Then, the map $F:U \times W \rightarrow M $ given by
$F(u,z):= \sigma(u) \cdot z $ is a diffeomorphism onto an open subset of $M$.
\end{description}
\item[(iv)] $G\cdot W$ is an open neighborhood of $G\cdot m$ 
and there is an equivariant smooth retraction
\[r:G\cdot W \longrightarrow G\cdot m\]
of the injection $G \cdot m \hookrightarrow G \cdot W $ such that $r^{-1}(m)=W$.
\end{description}
\end{theorem}

Let now $S:TN\times M\rightarrow TM$ be a $G$-invariant Stratonovich operator.
The existence of slices for the $G$-action allow us to carry out two decompositions of $S$. The first one, that we will call {\bfi  tangent-normal decomposition} is semi-global in the sense that it shares the properties that the Slice Theorem has in this respect, which is global in the  orbit directions and local in the directions transversal to the orbits; this decomposition consists of writing $S$
as the sum of two Stratonovich operators such that, roughly speaking,
one is tangent to the orbits of the $G$-action  and the other one is transversal
to them. The second one is purely local and yields a {\bfi skew-product decomposition} of $S$ in the sense of Definition~\ref{def skew-product}, provided that an additional hypothesis on the isotropies in the slice is present. This hypothesis, whose impact will be explained in detail later on, is generically satisfied and hence the following theorem shows that $S$ admits a skew product decomposition in a neighborhood of most points in $M$ (those points form an open and dense subset of $M$). 
These statements are rigorously proved in the following theorem.

\begin{theorem}
\label{teorema tan-nor decomposition}Let $X:\mathbb{R}_{+}\times
\Omega\rightarrow N$ be a $N$-valued semimartingale, $\Phi:G\times
M\rightarrow M$ a proper Lie group action, and $S:TN\times M\rightarrow TM$ a
$G$-invariant Stratonovich operator. Let $m\in  M$  and $W$ a slice at $m$. Then, there exist
two Stratonovich operators $S_{N}:TN\times
W\rightarrow TW$ and $S_{T}:TN\times G\cdot W\rightarrow T(G\cdot W)$ such
that the following statements hold:

\begin{description}
\item[(i)] Let $\operatorname*{Lie}\left(  N(G_{z})\right)  $ denote the Lie
algebra of the normalizer $N(G _z )$ in $G$ of the isotropy group $G_{z}$, $z\in G\cdot W$. The
Stratonovich operator $S_{T}$ is $G$-invariant and $S_{T}(x,z)\in
\mathcal{L}\left(  T_{x}N,\operatorname*{Lie}\left(  N(G_{z})\right)  \cdot
z\right)  $ for any $x\in N$ and any $z\in G\cdot W$. Moreover, there exists an adjoint
$G$-equivariant map $\xi:TN\times G\cdot W\rightarrow\mathfrak{g}$, (that is, $\xi\left(
x,g\cdot z\right)  =\operatorname*{Ad}_{g}\circ \xi\left(  x,z\right)  $, for
any $g\in G$) such that $S_{T}(x,z)=T_{e}\Phi_{z}\circ \xi\left(  x,z\right)$.

\item[(ii)] The Stratonovich operator $S_{N}:TN\times W\rightarrow TW$ is
$G_{m}$-invariant.
\item [(iii)]  If $z=g\cdot w\in G\cdot W$, with $g\in G$ and $w\in W$,
then
\begin{equation}
\label{decomposition for Stratonovich}
S\left(  x,z\right)  =S_{T}\left(  x,z\right)  +T_{w}\Phi_{g}\circ
S_{N}\left(  x,w\right)  =T_{w}\Phi_{g}\circ\left(  S_{T}\left(  x,w\right)
+S_{N}\left(  x,w\right)  \right)  .
\end{equation}
This sum of Stratonovich operators will be referred to as the
\textbf{tangent-normal decomposition} of $S$.

\item[(iv)] Let $\varphi$ be the flow of the stochastic system $\left(
W,S_{N},X,N\right)  $ so that $\varphi\left(  w\right)  $ denotes the solution
of
\begin{equation}
\delta\Gamma=S_{N}(X,\Gamma)\delta X\label{eq tan-nor 6}%
\end{equation}
with initial condition $\Gamma_{t=0}=w$ a.s.. Let
$S_{\mathfrak{g}\times W}:TN\times\left(  \mathfrak{g}\times W\right)
\rightarrow T\left(  \mathfrak{g}\times W\right)  $ be the Stratonovich
operator defined as $S_{\mathfrak{g}\times W}\left(  x,\left(  \eta,w\right)
\right)  =\xi\left(  x,w\right)  \times S_{N}\left(  x,w\right)
\in\mathcal{L}\left(  T_{x}N,\mathfrak{g}\times T_{w}W\right)  $ and let
$\left(  \eta^{w},\Gamma^{w}\right)  $ be the solution semimartingale of the
stochastic system $\left(  \mathfrak{g}\times W,S_{\mathfrak{g}\times
W},X,N\right)  $ with initial condition $\left(  0,w\right)  \in
\mathfrak{g}\times W$. Finally, let $\widetilde{g}:\left\{  0\leq
t<\tau_{\varphi}\right\}  \rightarrow G$ be the solution of the stochastic
system $\left(  G,L,\eta^{w},\mathfrak{g}\right)  $ with initial condition $g\in
G$ and where $L:T\mathfrak{g}\times
G\rightarrow TG$ is such that $L\left(  \eta,g\right)  (\nu)=T_{e}L_{g}(\nu)$.
Then, the semimartingale
\[
\Gamma_{t}=\widetilde{g}_{t}\cdot\varphi_{t}\left(  w\right)
\]
is a solution up to time $\tau_{\varphi}$ of the stochastic system $(M,S,X,N) $ with initial condition $z=g\cdot w\in G\cdot W$.

\item[(v)] Suppose now that  $G _w=G _z $, for any $w \in W $. Then $S$ admits a local skew-product decomposition. More specifically, for any point $m \in  M $,
there exists an open neighborhood $V\subseteq G/G_{m}$ of $G_{m}$, a
diffeomorphism $F:V\times W\rightarrow U\subseteq M$, and a skew-product
split Stratonovich operator $S_{V\times W}:TN\times\left(  V\times W\right)
\rightarrow T\left(  V\times W\right)  $ such that $F$ establishes a bijection
between semimartingales $\Gamma$ starting on $U$ which are solution of the
stochastic system $\left(  U,S,X,N\right)  $ and semimartingales on $V\times
W$ solution of the stochastic system $\left(  V\times W,S_{V\times
W},X,N\right)  $. Moreover,
\[
S_{V\times W}\left(  x,\left(  gG_{m},w\right)  \right)  =T_{g}\pi_{G_{m}
}\circ T_{e}L_{g}\left(\xi\left(  x,w\right) \right) \times S_{N}\left(  x,w\right)
\]
for any $x\in N$, $gG_{m}\in V \subset G/G_{m}$, and any $g\in G$ such that $\pi_{G_{m}%
}\left(  g\right)  =gG_{m}$.
\end{description}
\end{theorem}

\begin{remark}
\normalfont
The last point in this theorem shows that  proper symmetries of Stratonovich operators imply the availability of skew-products decompositions around most points in the manifold where the solutions take place. Indeed, the Principal Orbit Type Theorem (see for instance~\cite{dk}) shows that there exists an isotropy subgroup $H$ whose associated isotropy type manifold $M_{(H)}:=\{z \in  M\mid G _z=kHk ^{-1}, k \in  G\} $  is open and dense in $M$. Hence, for any point $m \in M_{(H)}$ there exist slice coordinates around the orbit $G \cdot m  $ in which the manifold $M$ looks locally like $G \times _H W= G \times _H W _H\simeq G/H \times W _H $. This local trivialization of the manifold $M$ into two factors and  the results in part {\bf (v)} of the theorem  can be used  to split the Stratonovich operator $S$, in order to obtain a locally defined skew-product around all the points in the open dense subset $M_{(H)}$ of $M$.
\end{remark}

\noindent\textbf{Proof.\ \ }As we already said, this construction is much inspired by a similar one available in the context of equivariant vector fields~\cite{krupa, field91}. In this proof we will mimic the strategy for that result followed in~\cite[Theorem 3.3.5]{hsr}. 

We start by noting that the properness of the action guarantees that  the isotropy subgroup $G _m$ is compact and hence there exists an open $G _m $-invariant neighborhood $V\subseteq
G/G_{m}$ of $G _m  $ and a local section $\sigma:V\subseteq G/G_{m}\rightarrow G$ with the following equivariance
property~\cite{field91}: $\sigma(h\cdot gG_{m})=h\sigma(gG_{m})h^{-1}$, for any $h \in  G _m  $ and $gG _m \in V $. If we now construct with this section the map
$F:V\times W\rightarrow U\subseteq M$ introduced in Theorem~\ref{slice theorem}, that is
\begin{equation}
F(gG_{m},w):=\sigma(gG_{m})\cdot w\label{eq tan-nor 10},
\end{equation}
we obtain a $ G _m $-equivariant map by considering  the  diagonal $G _m $-action in $V\times W $. Since for any $w\in W$ we have that $F^{-1}(w)=(G_{m},\sigma(G_{m}
)^{-1}\cdot w)$,
\begin{equation}
T_{w}F^{-1}\circ S(x,w)=:S_{V}(x,w)\times S_{W}(x,w)\in\mathcal{L}\left(
T_{x}N,T_{G_{m}}V\times T_{\sigma(G_{m})^{-1}\cdot w}W\right)
.\label{eq tan-nor 9}%
\end{equation}
Define%
\begin{subequations}
\begin{align}
S_{N}(x,w) &  :=T_{\sigma(G_{m})^{-1}\cdot w}\Phi_{\sigma(G_{m})}\circ
S_{W}(x,w)\in T_{w}W\label{eq tan-nor 2}\\
S_{T}(x,g\cdot w) &  :=T_{w}\Phi_{g}\circ T_{e}\Phi_{w}\circ T_{\sigma(G_{m}%
)}R_{\sigma(G_{m})^{-1}}\circ T_{G_{m}}\sigma\circ S_{V}(x,w)\nonumber\\
&  =T_{e}\Phi_{g\cdot w}\circ\operatorname*{Ad}\nolimits_{g}\circ
T_{\sigma(G_{m})}R_{\sigma(G_{m})^{-1}}\circ T_{G_{m}}\sigma\circ
S_{V}(x,w).\label{eq tan-nor 3}%
\end{align}
\end{subequations}

\medskip

\noindent {\bf (i)} Let $z=g\cdot w\in G\cdot W$, $g\in G$, $w\in
W$, $x \in N $, and define $\xi:TN\times G\cdot W\rightarrow\mathfrak{g}$ by%
\begin{equation}
\xi(x,z)=\operatorname*{Ad}\nolimits_{g}\circ T_{\sigma(G_{m})}R_{\sigma
(G_{m})^{-1}}\circ T_{G_{m}}\sigma\circ S_{V}(x,w).\label{eq tan-nor 13}%
\end{equation}
It can be seen that $\xi\left(  x,z\right)  $ is well defined by reproducing the
steps taken in \cite[Theorem 3.3.5 (i)]{hsr}. More specifically, it can be shown that if $z$ is
written as $z=g^{\prime}\cdot w^{\prime}$ for some other $g^{\prime}\in G$ and
$w^{\prime}\in W$ then%
\[
\operatorname*{Ad}\nolimits_{g}\circ T_{\sigma(G_{m})}R_{\sigma(G_{m})^{-1}%
}\circ T_{G_{m}}\sigma\circ S_{V}(x,w)=\operatorname*{Ad}\nolimits_{g^{\prime
}}\circ T_{\sigma(G_{m})}R_{\sigma(G_{m})^{-1}}\circ T_{G_{m}}\sigma\circ
S_{V}(x,w^{\prime}).
\]
Using~(\ref{eq tan-nor 3}) and~(\ref{eq tan-nor 13}) we have that
\[
S_{T}(x,g\cdot w)=T_{w}\Phi_{g \cdot w}\circ \xi(x, g \cdot  w)\]
It is an exercise to check that $\xi\left(
x,g\cdot w\right)  =\operatorname*{Ad}_{g}\circ \xi\left(  x,w\right)  $, for
any $g\in G$, and hence
the Stratonovich operator $S_{T}$ is $G$-invariant. This $G$-invariance implies by Proposition~\ref{law conservation isotropy} that the image
of $S_{T}(x,z)$ is such that $\operatorname*{Im}(S_{T}(x,z))\subseteq
T_{z}M_{G_{z}}$. On the other hand, $\operatorname*{Im}(S_{T}%
(x,z))=\operatorname*{Im}(T_{e}\Phi_{z}\circ\xi(x,z))\subseteq\mathfrak{g}%
\cdot z$, therefore
\[
\operatorname*{Im}\left(  S_{T}(x,z)\right)  \subseteq T_{z}M_{G_{z}}%
\cap\mathfrak{g}\cdot z=T_{z}\left(  N(G_{z})\cdot z\right)
\]
by \cite[Proposition 2.4.5]{hsr} and hence ${\rm Im}(\xi(x,z))\subset 
\operatorname*{Lie}\left(  N(G_{z})\right)  $.

\medskip

\noindent {\bf (ii)} and {\bf (iii)} It is immediate to see that the Stratonovich operator $S_{N}%
:TN\times W\rightarrow TW$ defined in (\ref{eq tan-nor 2}) is $G_{m}%
$-invariant. Let $w\in W$; using (\ref{eq tan-nor 9}) and 
(\ref{eq tan-nor 10})
\begin{align*}
S\left(  x,w\right)   &  =T_{\left(  G_{m},\sigma\left(  G_{m}\right)
^{-1}\cdot w\right)  }F\circ\left(  S_{V}(x,w)\times S_{W}(x,w)\right)  \\
&  =T_{e}\Phi_{w}\circ T_{\sigma(G_{m})}R_{\sigma(G_{m})^{-1}}\circ T_{G_{m}%
}\sigma\circ S_{V}(x,w)+T_{\sigma(G_{m})^{-1}\cdot w}\Phi_{\sigma(G_{m})}\circ
S_{W}(x,w)\\
&  =S_{T}\left(  x,w\right)  +S_{N}\left(  x,w\right)  ,
\end{align*}
where (\ref{eq tan-nor 2}) and (\ref{eq tan-nor 3}) have been used. The
equality~(\ref{decomposition for Stratonovich}) then follows from the $G$-invariance of $S$ and
$S_{T}$.

\medskip

\noindent {\bf (iv)} First of all observe that if $\left(  \eta^{w},\Gamma^{w}\right)  $
is the $\mathfrak{g}\times W$-valued semimartingale solution of the stochastic
system $\left(  \mathfrak{g}\times W,S_{\mathfrak{g}\times W},X,N\right)  $
with constant initial condition $\left(  0,w\right)  \in\mathfrak{g}\times W$,
then
\[
\left\langle \mu,\eta^{w}\right\rangle =\int\left\langle \xi\left(
X,\varphi_{t}(w)\right)  ^{\ast}(\mu),\delta X\right\rangle
\]
for any $\mu\in\mathfrak{g}^{\ast}$. In other words, $\eta^{w}$ may be
regarded as the solution of the stochastic differential equation
\begin{equation}
\delta\eta^{w}=\xi\left(  X,\varphi_{t}(w)\right)  \delta X
\label{eq tan-nor 7}%
\end{equation}
with initial condition $\eta_{t=0}^{w}=0$ a.s.. Notice that $\eta^{w}$ is defined up to time $\tau_{\varphi(w)}$, that is, the time of existence of the solution $\varphi (\omega) $.
Let now $\Gamma_{t}=\widetilde{g}_{t}\cdot\varphi_{t}\left(  w\right)  $ be the
$M$-valued semimartingale in the statement. Applying the rules of Stratonovich
differential calculus and the Leibniz rule we obtain
\begin{equation}
\label{in two parts stochastic differential equation}
\delta\Gamma_{t}=T_{\widetilde{g}_{t}}\Phi_{\varphi_{t}(w)}(\delta
\widetilde{g}_{t})+T_{\varphi_{t}(w)}\Phi_{\widetilde{g}_{t}}(\delta
\varphi_{t}(w))
\end{equation}
We rewrite the first summand in this expression as
\begin{align*}
T_{\widetilde{g}_{t}}\Phi_{\varphi_{t}(w)}(\delta\widetilde{g}_{t})  &
=T_{\varphi_{t}(w)}\Phi_{\widetilde{g}_{t}}\circ T_{e}\Phi_{\varphi_{t}%
(w)}\circ T_{\widetilde{g}_{t}}L_{\widetilde{g}_{t}^{-1}}(\delta\widetilde
{g}_{t})\\
&  =T_{\varphi_{t}(w)}\Phi_{\widetilde{g}_{t}}\circ T_{e}\Phi_{\varphi_{t}%
(w)}(\delta\eta_{t}^{w})\\
&  =T_{\varphi_{t}(w)}\Phi_{\widetilde{g}_{t}}\circ T_{e}\Phi_{\varphi_{t}%
(w)}\circ\xi(X_{t},\varphi_{t}(w))\delta X_{t}\\
&  =T_{\varphi_{t}(w)}\Phi_{\widetilde{g}_{t}}\circ S_{T}(X,\varphi
_{t}(w))\delta X_{t},
\end{align*}
where in the second and third line we have used that $\widetilde{g}_{t}$ is a
solution of $\left(  G,L,\eta^{w},\mathfrak{g}\right)  $ and equation
(\ref{eq tan-nor 7}), respectively. The second summand of ~(\ref{in two parts stochastic differential equation}) can be written as
\[
T_{\varphi_{t}(w)}\Phi_{\widetilde{g}_{t}}(\delta\varphi_{t}(w))=T_{\varphi
_{t}(w)}\Phi_{\widetilde{g}_{t}}\circ S_{N}(X,\varphi_{t}(w))\delta X_{t}%
\]
because $\varphi_{t}(w)$ is a solution of (\ref{eq tan-nor 6}). Therefore,
using~(\ref{decomposition for Stratonovich}) we can conclude that
\begin{align*}
\delta\Gamma_{t}  &  =T_{\varphi_{t}(w)}\Phi_{\widetilde{g}_{t}}\circ\left(
S_{N}(X,\varphi_{t}(w))+S_{T}(X,\varphi_{t}(w))\right)  \delta X_{t}\\
&  =S\left(  X,\widetilde{g}_{t}\cdot\varphi_{t}(w)\right)  \delta
X_{t}=S\left(  X,\Gamma_{t}\right)  \delta X_{t}%
\end{align*}
which shows that $\Gamma_{t}$ is a solution up to time $\tau_{\varphi}$ of the stochastic system $(M,S,X,N) $ with initial condition $z=g\cdot w\in G\cdot W$

\medskip

\noindent {\bf (v)} Let $w\in W$ and $h\in G_{m}=G_{w}$. Let $\Psi$ be the twisted action of $G _m  $ on $W$, that is,  $\Psi:G_{m}\times\left(  G\times W\right)
\rightarrow\left(  G\times W\right)  $ defined as $\Psi_{h}\left(  g,w\right)
:=\left(  gh,h^{-1}\cdot w\right)  $, and whose orbit space is the twisted product $G\times_{G_{m}}%
W$. The hypothesis $G_{m}=G_{w}$, for any $w \in W $, implies that  $G\times_{G_{m}}W$ can be easily identified with
$G/G_{m}\times W$ using the diffeomorphism
\[%
\begin{array}
[c]{rrl}%
G\times_{G_{m}}W & \longrightarrow & G/G_{m}\times W\\
\left[  g,w\right]   & \longmapsto & \left(  gG_{m},w\right)  .
\end{array}
\]
Consider now the Stratonovich operator defined by
\[
S_{G\times W}\left(  x,\left(  g,w\right)  \right)  =T_{e}L_{g}\circ\xi\left(
x,w\right)  \times S_{N}\left(  x,w\right).
\]
We are going to show that $S_{G\times W} $ is $G_{m}$-invariant under the action defined by $\Psi$. Indeed, given that $G_{w}=G_{m}$, $\Psi_{h}\left(
g,w\right)  =\left(  gh,w\right)  $ for any $h\in G_{m}$, $g\in G$, and $w\in
W$,  we have%
\begin{align}
S_{G\times W}\left(  x,\Psi_{h}\left(  g,w\right)  \right)   &  =S_{G\times
W}\left(  x,\left(  gh,w\right)  \right)  =T_{e}L_{gh}\circ\xi\left(
x,w\right)  \times S_{N}\left(  x,w\right)  \nonumber\\
&  =T_{h}L_{g}\circ T_{e}L_{h}\circ\xi\left(  x,w\right)  \times S_{N}\left(
x,w\right)  \nonumber\\
&  =T_{h}L_{g}\circ T_{e}R_{h}\circ\operatorname*{Ad}\nolimits_{h}\circ
\xi\left(  x,w\right)  \times S_{N}\left(  x,w\right)  \nonumber\\
&  =T_{g}R_{h}\circ T_{e}L_{g}\circ\operatorname*{Ad}\nolimits_{h}\circ
\xi\left(  x,w\right)  \times S_{N}\left(  x,w\right)  .\label{eq tan-nor 5}%
\end{align}
But due to the $G$-equivariance of $\xi$ we have  $\xi\left(  x,w\right)  =\xi\left(  x,h\cdot w\right)
=\operatorname*{Ad}\nolimits_{h}\circ\xi\left(  x,w\right)  $, for any $h\in G_{m}$. In addition, $T_{\left(
g,w\right)  }\Psi_{h}=T_{g}R_{h}\times\operatorname*{Id}$, so
(\ref{eq tan-nor 5}) equals%
\[
T_{\left(  g,w\right)  }\Psi_{h}\circ\left(  T_{e}L_{g}\circ\xi\left(
x,w\right)  \times S_{N}\left(  x,w\right)  \right)  =T_{\left(  g,w\right)
}\Psi_{h}\circ S_{G\times W}\left(  x,(g,w)\right),
\]
which shows that  $S_{G\times W}$ is $G_{m}$-invariant.

We can therefore apply the Reduction Theorem~\ref{reduction theorem statement general systems} to conclude that $S_{G\times W}$ projects onto a 
stochastic system $\left(  G/G_{m}\times W,S_{G/G_{m}\times W},X,N\right)  $
on $G\times_{G_{m}}W\simeq G/G_{m}\times W$ with Stratonovich operator%
\begin{equation}
S_{G/G_{m}\times W}\left(  x,\left(  gG_{m},w\right)  \right)  :=T_{g}%
\pi_{G_{m}}\circ S_{G\times W}\left(  x,\left(  g,w\right)  \right)  =T_{g}%
\pi_{G_{m}}\circ T_{e}L_{h}\circ\xi\left(  x,w\right)  \times S_{N}\left(
x,w\right)  ,\label{eq tan-nor 14}%
\end{equation}
where $x\in N$, $w\in W$, and $g\in G$ is any element such that $\pi_{G_{m}%
}\left(  g\right)  =gG_{m}$. 
Notice that by~(\ref{decomposition for Stratonovich}), expression~(\ref{eq tan-nor 14}) proves that the Stratonovich operator
$S_{G/G_{m}\times W}$ is a local skew-product decomposition of $S$ on $G/G_{m}\times W$.

Concerning the solutions, by {\bf (iv)} any solution of the stochastic system $(M,S,X,N)$ starting at
some point $z=g\cdot w\in U\subseteq G\cdot W$ can be written as the image by
the action $\Phi$ of the solution $\left( g_{t},\varphi_{t}\left(  w\right)
\right)  $ of the stochastic system $\left(  G\times W,S_{G\times
W},X,N\right)  $ starting at $\left(  g,w\right)  \in G\times W$ and defined
up to time $\tau_{\varphi(w)}$. Then, the Reduction Theorem~\ref{reduction theorem statement general systems} guarantees that this solution 
can be  projected to a solution of $\left(  G/G_{m}\times W,S_{G/G_{m}\times
W},X,N\right)  $ starting at $\left(  gG_{m},w\right)  \in G/G_{m}\times W$, also defined up to time $\tau_{\varphi(w)}$. 
Conversely, in order to recover a solution of the original system from a
solution $\left(  \left(  gG_{m}\right)  _{t},w_{t}\right)  $ of 
$\left(  G/G_{m}\times W,S_{G/G_{m}\times W},X,N\right)  $ we need to
invoke the Reconstruction Theorem~\ref{teorema reconstruccion 1} by choosing an auxiliary connection $A\in\Omega^{1}\left(  G;\mathfrak{g}_{m}\right)  $. 
This will yield a solution $\left(
g_{t},w_{t}\right)  $ of $\left(  G\times W,S_{G\times W},X,N\right)  $ with 
$g_{t}$ a $G$-valued semimartingale that can be written as%
\[
g_{t}=d_{t}h_{t},
\]
where $d_{t}:\mathbb{R}_{+}\times\Omega\rightarrow G$ is the horizontal lift
of $\left(  gG_{m}\right)  _{t}$ with respect to $A$ and $h_{t}:\mathbb{R}
_{+}\times\Omega\rightarrow G_{m}$ is a suitable semimartingale on $G_{m}$.
The key point is that the image by the action $\Phi  $ of the solution $\left(
g_{t},w_{t}\right)  $ of $\left(  G\times W,S_{G\times W},X,N\right)  $, that is,
\[
\Phi\left(  g_{t},w_{t}\right)  =g_{t}\cdot w_{t}=d_{t}h_{t}\cdot w_{t}%
=d_{t}\cdot w_{t}%
\]
yields a solution of $(M,S,X,N)$. Notice that the semimartingale $h_{t}$
plays no role. Indeed, let
$\sigma:V\subseteq G/G_{m}\rightarrow G$ be the local $G _m$-equivariant section introduced
in the beginning of the proof. We already saw in Proposition
\ref{skew-product caso libre} that if $\left(  gG_{m}\right)  _{t}%
:\mathbb{R}_{+}\times\Omega\rightarrow G/G_{m}$ is a $G/G_{m}$-valued
semimartingale then $\sigma\left(  \left(  gG_{m}\right)  _{t}\right)  $ is
the horizontal lift with respect to the connection $A_{\sigma}\in\Omega
^{1}\left(  \pi_{G_{m}}^{-1}\left(  V\right)  ;\mathfrak{g}_{m}\right)  $
induced by the local section $\sigma$. Consequently, any solution $\Gamma_{t}$
of the initial stochastic system $\left(  M,S,X,N\right)  $ with initial
condition $\Gamma_{t=0}=g\cdot w\in U\subset G\cdot W$ can be locally
expressed as $\sigma\left(  \left(  gG_{m}\right)  _{t}\right)  \cdot w_{t}$
where $\left(  \left(  gG_{m}\right)  _{t},w_{t}\right)  $ is a solution of
the stochastic system $\left(  G/G_{m}\times W,S_{G/G_{m}\times W},X,N\right)
$ with initial condition $\left(  \pi_{G_{m}}(g),w\right)  \in G/G_{m}\times
W$.
\quad $\blacksquare$

\begin{example}[Liao decomposition of Markov processes]
\normalfont
The possibility of decomposing stochastic processes using a group invariance property has been used beyond the context of stochastic differential equations. For example, Liao~\cite{liao
decomposition} has used what he calls the {\bfi  transversal submanifolds} of a compact group action to carry out an {\bfi  angular-radial decomposition} of  the Markov processes that are equivariant with respect to those actions. To be more specific, let $M$ be a manifold acted upon by a Lie group $G$ and let $\Gamma:\mathbb{R}_{+}\times\Omega\rightarrow M$ be a $M$-valued Markov process with transition semigroup $P _t $; that is, $\Gamma$ is a process with c\`adl\`ag paths that satisfies the simple Markov property
\[
E\left[  f\left(  \Gamma_{t+s}\right)  |\mathcal{F}_{t}\right]  =P_{s}f\left(
\Gamma_{t}\right)
\]
a.s. for $s<t$ and $f\in C_{b}^{\infty}\left(  M\right)  $, where
$C_{b}^{\infty}\left(  M\right)  $ is the space of bounded smooth functions on
$M$, and $\left\{  \mathcal{F}_{t}\right\}  _{t\in\mathbb{R}_{+}}$ is the
natural filtration induced by $\Gamma$. Furthermore, suppose that the
Markov process $\Gamma$ or, equivalently, its transition semigroup $P_{t}$ is
$G$-equivariant in the sense that
\[
P_{t}\left(  f\circ\Phi_{g}\right)  =\left(P_{t}f \right)\circ\Phi_{g}%
\]
for any $g\in G$. Additionally, in~\cite{liao
decomposition} it is assumed the existence of a submanifold $W\subseteq
M$ which is \textbf{globally transversal} to the $G$-action. This means that $W$
intersects each $G$-orbit at exactly one point, that is, for any $w \in W $,  $G\cdot w\cap W=\left\{  w\right\}
$ and $M=%
{\textstyle\bigcup_{w\in W}}
G \cdot  w  $. The existence of such global transversal section  is a strong hypothesis that only a limited number of actions satisfy. A larger range of applicability of the results in~\cite{liao decomposition} can be obtained if one is willing to work locally using the slices introduced in this section. Indeed, suppose now that the group $G$ is not compact but just that the group action is proper; let $m\in M$ and $\varphi:G\times_{G_{m}}W\rightarrow U\subseteq
M$ a tube around the orbit $G \cdot m $ where, additionally, we assume that $G_{w}=G_{m}$ for any $w\in W$. With this hypothesis which, incidentally is the same one that in part {\bf (v)} of Theorem~\ref{teorema tan-nor decomposition} allowed us to obtain a skew-product decomposition of the invariant Stratonovich operator, the slice $W$ \emph{is a local transversal manifold} in the sense of~\cite{liao decomposition}. 

Let now $J:U\subseteq M\rightarrow W$ be the projection that associates to each point, the unique element in its orbit that intersects $W $. Liao proves~\cite[Theorem 1]{liao decomposition} that  the {\bfi  radial part} $y:=J\left(  \Gamma\right)
$ of the Markov process $\Gamma$ is also a Markov process with transition semigroup
$Q_{t}:=J^{\ast}P_{t}$. Moreover, if the group $G$ is compact and $ \Gamma $ is Feller then so is $y$ and its generator is fully determined by that of $\Gamma $.

Let now $\pi_{G_{m}}:G\rightarrow G/G_{m}$ be the canonical projection and let
$\phi:V\times W\rightarrow U$ be the diffeomorphism associated to the
local section $\sigma:V\rightarrow\pi_{G_{m}}^{-1}\left(  V\right)  \subseteq
G$ such that $\phi\left(  gG_{m},w\right)  =\sigma\left(  gG_{m}\right)
\cdot w$. Let $\Gamma$ be $U$-valued Markov process starting at $m$ and
$y=J\left(  \Gamma\right)  $ its radial part. Let $\bar{\Gamma}:\left\{  0\leq
t<\tau_{U}\right\}  \rightarrow V\subseteq G/G_{m}$ be the process such that
$\Gamma_{t}=\sigma\left(  \bar{\Gamma}_{t}\right)  \cdot y_{t}$, where
$\tau_{U}=\inf\left\{  t>0~|~\Gamma_{t}\notin U\right\}  $. $\bar{\Gamma}$ is
called the {\bfi  angular part} of $\Gamma$. Liao shows (see~\cite[Theorem
3]{liao decomposition}) that the angular process $\bar{\Gamma}_{t}$ is a
{\bfi  nonhomogeneous L\'{e}vy process} under the conditional probability
built by conditioning with respect to the $\sigma$-algebra generated by the
radial process. The reader is encouraged to check with \cite{liao
decomposition} for precise definitions and statements (see also \cite{liao book}). 
\end{example}

\section{Projectable stochastic differential equations on associated bundles}
\label{Projectable stochastic differential equations on associated bundles}

In the previous section we saw how the availability of the slices associated to a proper group action allows the local splitting of the invariant Stratonovich operators using what we called the tangent-normal decomposition. Additionally, this decomposition yields generically a local skew-product splitting of the invariant Stratonovich operator in question. The key idea behind these splittings was the possibility of locally modeling the manifold where the solutions of the stochastic differential equation take place as a twisted product. A natural setup that we could consider are the manifolds $M$ where this product structure is \emph{global}, that is $M=P\times_{G}W$, with $P $ and $W$ two $G$-manifolds. The most standard situation where such manifolds are encountered is when $M$ is the {\bfi  associated bundle} to the $G$-principal bundle $ \pi: P \rightarrow  Q $: let $W$ be an effective left $G$-space and $\bar{\pi}:P\times_{G}W\rightarrow Q$, $\bar{\pi}([p,w])= \pi (p ) $. A classical theorem in bundle theory shows that such construction is a principal $G$-bundle with fiber $W$ and it is usually referred to as the bundle associated to $ \pi: P \rightarrow  Q $ with fiber $W$. To be more specific, consider the commutative diagram that defines $\bar \pi$:
\begin{equation}%
\begin{array}
[c]{rrl}%
P\times W & \overset{\kappa}{\longrightarrow} & P\times_{G}W\\
_{pr_{1}}\downarrow &  & \downarrow_{\bar{\pi}}\\
P & \overset{\pi}{\longrightarrow} & Q.
\end{array}
\label{eq tan-nor 17}%
\end{equation}
In this diagram, $\kappa_{p}:\{p\}\times W\rightarrow\bar{\pi}^{-1}%
(\pi(p))=:\left(  P\times_{H}W\right)  _{\pi(p)}$ is a diffeomorphism
(see for instance~\cite[10.7]{michor}). Hence, the correspondence $p\rightarrow\kappa_{p}$,
$p\in P$, allows us to consider the elements of $P$ as diffeomorphisms from the
typical fiber $W$ of $P\times_{G}W$ to $\bar{\pi}^{-1}\left(  q\right)  $,
with $q=\pi(p)$. 

Stochastic processes and diffusions on associated bundles have
deserved certain attention in the literature (see \cite{liao skew-product} for
example) because, as we will see in the following paragraphs, the available geometric structure makes possible a Reduction-Reconstruction procedure that in some cases implies the existence of a global skew-product decomposition. In this context, the notion of invariance is replaced by what we will call $\bar{\pi}$-projectability: if $N$ is a manifold and
$S:TN\times M\rightarrow TM$ a Stratonovich operator from $N$ to $M$, we
say that $S$ is $\bar{\pi}$\textbf{-projectable} if the Stratonovich operator $S _Q$
from $N$ to $Q$%
\[
S_Q\left(  x,q\right)  :=T_{\left[  p,w\right]  }\bar{\pi}\circ S\left(
x,\left[  p,w\right]  \right)  \in\mathcal{L}\left(  T_{x}N,T_{\left[
p,w\right]  }M\right)
\]
is well defined, where $\left[  p,w\right]  \in M$ is any point such that
$\bar{\pi}\left(  \left[  p,w\right]  \right)  =q\in Q$.

\begin{theorem}
\label{teorema tan-norm associated bundles}Let $\bar \pi:M=P\times
_{G}W\rightarrow Q$ be the associated bundle introduced in the previous discussion. Let $N$ be a manifold, $S:TN\times
M\rightarrow TM$ a $\bar{\pi}$-projectable Stratonovich operator onto $Q$, and
$X:\mathbb{R}_{+}\times\Omega\rightarrow N$ a $N$-valued semimartingale. Then
there exist a Stratonovich operator $S_{P\times W}:TN\times\left(  P\times
W\right)  \rightarrow TP\times TW$ with the property that if $\left(
p_{t},w_{t}\right)  $ is any solution of the stochastic system $\left(
P\times W,S_{P\times W},X,N\right)  $ with initial condition $\left(
p,w\right)  \in P\times W$, then $\Gamma _t:=\kappa\left(  p_{t},w_{t}\right)  $ is the
solution of $\left(  M,S,X,N\right)  $ starting at $\left[  p,m\right]  $.
Furthermore, $p_{t}$ can be written as the horizontal lift of $\bar{\pi}\left(  \Gamma
_{t}\right)  $ with respect to an auxiliary connection $A\in\Omega^{1}\left(
P;\mathfrak{g}\right)  $. Conversely, if $\Gamma_{t}$ is a solution of $\left(  M,S,X,N\right)  $ and $p_{t}$ the
horizontal lift of $\bar{\pi}\left(  \Gamma_{t}\right)  $ with respect to $A$,
then $\left(  p_{t},\kappa_{p_{t}}^{-1}\left(  \Gamma_{t}\right)  \right)  $
is a solution of $\left(  P\times W,S_{P\times W},X,N\right)  $.
\end{theorem}

\noindent\textbf{Proof.\ \ }Let $A\in\Omega^{1}\left(  P;\mathfrak{g}\right)  $ be an auxiliary principal
connection for $\pi:P\rightarrow Q$ and let $\widehat{A}_{p}:T_{\pi
(p)}Q\rightarrow\operatorname*{Hor}_{p}P\subseteq T_{p}P$ be the inclusion of
the tangent space $T_{q}Q$ at $q=\pi(p)$ into the horizontal space
$\operatorname*{Hor}_{p}P$ at $p\in P$ defined by $A$. Consider the family of
linear maps $\widehat{\mathbf{A}}_{\left[  p,w\right]  }:T_{\bar{\pi}(\left[
p,w\right]  )}Q\rightarrow T_{\left[  p,w\right]  }M$ for any $\left[
p,w\right]  \in P\times_{G}W$ as%
\begin{equation}
\widehat{\mathbf{A}}_{\left[  p,w\right]  }=T_{p}\kappa_{w}\circ\widehat
{A}_{p},\label{eq tan-nor 15}%
\end{equation}
where $\kappa_{w}(p):=\kappa\left(  p,w\right)  $ for any $w\in W$. The family of maps $\{\widehat{\mathbf{A}}_{\left[  p,w\right]  }~|~\left[
p,w\right]  \in M\}$ define what is called the \textbf{induced connection}
$\mathbf{A}$ (\cite[11.8]{michor}) on $P\times_{G}W$ by $A\in\Omega
^{1}\left(  P;\mathfrak{g}\right)  $. It can be easily checked that
$\mathbf{A}$ is well-defined, that is, the expression (\ref{eq tan-nor 15})
does not depend on the particular choice of $p\in P$ and $w\in W$ in the
class $\left[  p,w\right]  \in P\times_{G}W$ used to define it. Indeed, if $\left[
p,w\right]  =\left[  p^{\prime},w^{\prime}\right]  $ then there exists some
$g\in G$ such that $p^{\prime}=g\cdot p$ and $w^{\prime}=g^{-1}\cdot w$. Since
the connection $A$ is principal, $\widehat{A}_{p^{\prime}}=T_{p}R_{g}%
\circ\widehat{A}_{p}$, where $R:G\times P\rightarrow P$ denotes the $G$-right
action on $P$. On the other hand, since $\kappa\left(  p^{\prime}=p\cdot
g,w^{\prime}\right)  =\kappa\left(  p,g\cdot w^{\prime}\right)  $, we have%
\begin{equation}
T_{p^{\prime}}\kappa_{w^{\prime}}\circ T_{p}R_{g}=T_{p}\kappa_{g\cdot
w^{\prime}}\text{ \ or, equivalently, \ }T_{p^{\prime}}\kappa_{w^{\prime}%
}=T_{p}\kappa_{g\cdot w^{\prime}}\circ T_{p^{\prime}}R_{g^{-1}}%
.\label{eq tan-nor 21}%
\end{equation}
Therefore, $\widehat{\mathbf{A}}_{\left[  p^{\prime},w^{\prime}\right]
}=T_{p^{\prime}}\kappa_{w^{\prime}}\circ\widehat{A}_{p^{\prime}}=T_{p}%
\kappa_{g\cdot w^{\prime}}\circ T_{p^{\prime}}R_{g^{-1}}\circ T_{p}R_{g}%
\circ\widehat{A}_{p}=T_{p}\kappa_{w}\circ\widehat{A}_{p}=\widehat{\mathbf{A}%
}_{\left[  p,w\right]  }$.

Let $S_{Q}:TN\times Q\rightarrow TQ$ be the Stratonovich operator defined as%
\begin{equation}
S_{Q}\left(  x,q\right)  :=T_{\left[  p,w\right]  }\bar{\pi}\circ S\left(
x,\left[  p,w\right]  \right)  , \label{eq tan-nor 18}%
\end{equation}
where $\left[  p,w\right]  \in P\times_{G}W$ is any point such that $\bar{\pi
}\left(  \left[  p,w\right]  \right)  =q$, $x\in N$, and $w\in W$. This
Stratonovich operator is well-defined because $S$ is by hypothesis $\bar{\pi}%
$-projectable. Let $\widehat{H}_{\left[  p,w\right]  }:T_{\left[  p,w\right]
}M\rightarrow\operatorname*{Hor}_{\left[  p,w\right]  }M\subseteq T_{\left[
p,w\right]  }M$ and $\widehat{V}_{\left[  p,w\right]  }:T_{\left[  p,w\right]
}M\rightarrow\operatorname*{Ver}_{\left[  p,w\right]  }M\subseteq T_{\left[
p,w\right]  }M$ be the projections onto the horizontal and vertical spaces
associated to $\mathbf{A}$, respectively, at $\left[  p,w\right]  \in
P\times_{G}W$. Define the Stratonovich operator $S_{P\times W}:TN\times
(P\times W)\rightarrow TP\times TW$ as%
\begin{equation}
S_{P\times W}\left(  x,\left(  p,w\right)  \right)  =\widehat{A}_{p}\circ
S_{Q}\left(  x,\pi\left(  p\right)  \right)  \times\left(  T_{w}\kappa
_{p}\right)  ^{-1}\circ\widehat{V}_{\left[  p,w\right]  }\circ S\left(
x,\left[  p,w\right]  \right)  \in\mathcal{L}\left(  T_{x}N,T_{(p,w)}(P\times
W)\right)  \label{eq tan-nor 16}%
\end{equation}
for any $x\in N$, $w\in W$, and $p\in P$. Recall from (\ref{eq tan-nor 17})
that $\kappa_{p}:W\rightarrow M_{\pi(p)}$ is a diffeomorphism for any $p\in P$
and hence $\left(  T_{w}\kappa_{p}\right)  ^{-1}$ exists as a map. Now, we claim that if
$\left(  p_{t},w_{t}\right)  $ is a $\left(  P\times W\right)  $-valued
semimartingale solution of the stochastic system $\left(  P\times W,S_{P\times
W},X,N\right)  $ then $\Gamma_{t}:=\kappa_{p_{t}}(w_{t})$ is a solution of
$\left(  M,S,X,N\right)  $. Indeed, applying the Stratonovich rules for
differential calculus,%
\begin{align*}
\delta\Gamma_{t}  &  =T_{w_{t}}\kappa_{p_{t}}\left(  \delta w_{t}\right)
+T_{p_{t}}\kappa_{w_{t}}\left(  \delta p_{t}\right) \\
&  =T_{w_{t}}\kappa_{p_{t}}\circ\left(  T_{w_{t}}\kappa_{p_{t}}\right)
^{-1}\circ\widehat{V}_{\left[  p_{t},w_{t}\right]  }\circ S\left(
X_{t},\left[  p_{t},w_{t}\right]  \right)  \delta X_{t}+T_{p_{t}}\kappa_{w_{t}}\circ\widehat{A}_{p_{t}}\circ S_{Q}\left(
X_{t},\pi\left(  p_{t}\right)  \right)  \delta X_{t}\\
&  =\widehat{V}_{\left[  p_{t},w_{t}\right]  }\circ S\left(  X_{t},\left[
p_{t},w_{t}\right]  \right)  \delta X_{t}+\widehat{\mathbf{A}}_{\left[
p_{t},w_{t}\right]  }\circ S_{Q}\left(  X_{t},\pi\left(  p_{t}\right)
\right)  \delta X_{t}\\
&  =\widehat{V}_{\left[  p_{t},w_{t}\right]  }\circ S\left(  X_{t},\left[
p_{t},w_{t}\right]  \right)  \delta X_{t}+\widehat{\mathbf{A}}_{\left[
p_{t},w_{t}\right]  }\circ T_{\left[  p_{t},w_{t}\right]  }\bar{\pi}\circ
S\left(  X_{t},\left[  p_{t},w_{t}\right]  \right)  \delta X_{t}\\
&  =\widehat{V}_{\left[  p_{t},w_{t}\right]  }\circ S\left(  X_{t},\left[
p_{t},w_{t}\right]  \right)  \delta X_{t}+\widehat{H}_{\left[  p_{t}%
,w_{t}\right]  }\circ S\left(  X_{t},\left[  p_{t},w_{t}\right]  \right)
\delta X_{t} =S\left(  X_{t},\left[  p_{t},w_{t}\right]  \right)  \delta X_{t}=S\left(
X_{t},\Gamma_{t}\right)  \delta X_{t},
\end{align*}
and hence $\Gamma_{t}$ is a solution of $\left(  M,S,X,N\right)  $.

Conversely, let $\Gamma_{t}$ be a solution of $\left(  M,S,X,N\right)  $ such
that $\Gamma_{t=0}=\left[  p,m\right]  $ a.s. and let $p_{t}$ be the horizontal
lift of $\bar{\pi}\left(  \Gamma_{t}\right)  $ with respect to the auxiliary
connection $A\in\Omega^{1}\left(  P;\mathfrak{g}\right)  $ starting at some
$p_{0}\in\pi^{-1}\left(  \bar{\pi}\left(  \left[  p,w\right]  \right)
\right)  $. Define $\widetilde{w}_{t}:=\kappa_{p_{t}}^{-1}\left(  \Gamma
_{t}\right)  $. Observe that $\widetilde{w}_{t=0}=\widetilde{w}_{0}$ is such
that $\left[  p_{0},\widetilde{w}_{0}\right]  =\left[  p,m\right]  $. Since
$\kappa_{p}:W\rightarrow M_{\pi\left(  p\right)  }$ is a diffeomorphism,
$\widetilde{w}_{t}$ is uniquely determined a.s. by $\Gamma_{t}$ once $p_{t}$
is fixed. Indeed, $\widetilde{w}_{t}$ is the unique semimartingale such that
$\kappa_{p_{t}}\left(  \widetilde{w}_{t}\right)  =\Gamma_{t}$. But we have already seen that the
solution of $\left(  P\times W,S_{P\times W},X,N\right)  $ starting at
$\left(  p_{0},\widetilde{w}_{0}\right)  \in P\times W$ may be expressed as $\left(  p_{t},w_{t}\right)  $, with $p_{t}$ the fixed
horizontal lift of $\bar{\pi}\left(  \Gamma_{t}\right)  $ that we have been using all along. Therefore $w_{t}=\widetilde{w}_{t}$ a.s. necessarily and $w_{t}=\kappa_{p_{t}}%
^{-1}\circ\kappa_{p_{t}}\left(  w_{t}\right)  =\kappa_{p_{t}}^{-1}\left(
\Gamma_{t}\right)  $.
\quad $\blacksquare$
\begin{corollary}
\label{corollario associated bundles}
Using the same notation as in the proof of Theorem~\ref{teorema tan-norm associated bundles}, suppose that $\left(  T_{w}%
\kappa_{p}\right)  ^{-1}\circ\widehat{V}_{\left[  p,w\right]  }\circ S\left(
x,\left[  p,w\right]  \right)  $ in (\ref{eq tan-nor 16}) does not depend on
$p\in P$. In such case there exists a unique $G$-invariant Stratonovich operator $S_{W}:TN\times
W\rightarrow TW$ from $N$ to $W$ determined by the relation
\begin{equation}
T_{w}\kappa_{p}\circ S_{W}\left(  x,w\right)  =\widehat{V}_{\left[
p,w\right]  }\circ S\left(  x,\left[  p,w\right]  \right)
\label{eq tan-nor 20}%
\end{equation}
for any $x\in N$, $w\in W$, and $p\in P$. Moreover, $S_{P\times W}$ in
(\ref{eq tan-nor 16}) admits the  skew-product decomposition
\[
S_{P\times W}\left(  x,\left(  p,w\right)  \right)  =\widehat{A}_{p}\circ
S_{Q}\left(  x,\pi\left(  p\right)  \right)  \times S_{W}\left(  x,w\right).
\]
\end{corollary}

\noindent\textbf{Proof.\ \ }First of all notice that as $\left(  T_{w}%
\kappa_{p}\right)  ^{-1}\circ\widehat{V}_{\left[  p,w\right]  }\circ S\left(
x,\left[  p,w\right]  \right)  $  does not depend on
$p\in P$, the expression~(\ref{eq tan-nor 20}) is a good definition that uniquely determines $S _W $. The only non-trivial point in the statement that needs proof is the $G$-invariance of $S _W $: 
let $g \in G $ and $\left(  p^{\prime},w^{\prime}\right)  $, $\left(  p,w\right)  \in
P\times W$  such that $p^{\prime
}=p\cdot g$ and $w^{\prime}=g^{-1}\cdot w$. Since $\widehat{V}_{\left[
p^{\prime},w^{\prime}\right]  }\circ S\left(  x,\left[  p^{\prime},w^{\prime
}\right]  \right)  =\widehat{V}_{\left[  p,w\right]  }\circ S\left(  x,\left[
p,w\right]  \right)  $, we necessarily have
\[
T_{w}\kappa_{p}\circ S_{W}\left(  x,w\right)  =T_{w^{\prime}}\kappa
_{p^{\prime}}\circ S_{W}\left(  x,w^{\prime}\right)  .
\]
As $\kappa\left(  p\cdot g,w\right)  =\kappa\left(
p,g\cdot w\right)  $, we have that $T_{g\cdot w}\kappa_{p}\circ T_{w}l_{g}=T_{w}%
\kappa_{p\cdot g}$, where $l:G\times W\rightarrow W$ is the $G$-action on $W$.
Thus,%
\[
T_{w^{\prime}}\kappa_{p^{\prime}}\circ S_{W}\left(  x,w^{\prime}\right)
=T_{w}\kappa_{p}\circ T_{g^{-1}\cdot w}l_{g}\circ S_{W}\left(  x,g^{-1}\cdot
w\right)  .
\]
Since $T_{w}\kappa_{p}:T_{w}W\rightarrow T_{\left[  p,w\right]  }\left(
P\times_{G}W\right)  $ is an isomorphism, we conclude comparing the two
previous relations that%
\[
S_{W}\left(  x,w\right)  =T_{g^{-1}\cdot w}l_{g}\circ S_{W}\left(
x,g^{-1}\cdot w\right)  ,
\]
necessarily, which amounts to $S_{W}$ being $G$-invariant.
\quad $\blacksquare$
\begin{remark}
\normalfont
It is worth noticing that, under the hypotheses of Corollary
\ref{corollario associated bundles} and unlike Theorem
\ref{teorema tan-nor decomposition}, the skew-product decomposition of
$S_{P\times W}:TN\times\left(  P\times W\right)  \rightarrow T\left(  P\times
W\right)  $ is now global.
\end{remark}

\begin{remark}
\normalfont
If the hypotheses of Corollary \ref{corollario associated bundles} hold, we
can solve a stochastic system $\left(  M,S,X,N\right)  $ on the associated
bundle $\bar{\pi}:M=P\times_{G}W\rightarrow Q$ with $\bar{\pi}$-projectable
Stratonovich operator $S$ using the following reduction-reconstruction scheme.
On one hand, we find the solution starting at $\bar{\pi}\left(  \left[
p,w\right]  \right)  $ on the base space system $\left(  Q,S_{Q},X,N\right)  $,
where $S_{Q}$ was given in (\ref{eq tan-nor 18}). We lift then this solution
to the principal bundle $P$ using an auxiliary connection
$A\in\Omega^{1}\left(  P;\mathfrak{g}\right)  $. We choose the lift $p_{t}$
starting at some $p_{0}\in\pi^{-1}\left(  \bar{\pi}\left(  \left[  p,w\right]
\right)  \right)  $. On the other hand, we find the solution $w_{t}$ of the
independent stochastic system $\left(  W,S_{W},X,N\right)  $ with initial condition
$w_{0}$ such that $\kappa\left(  p_{0},w_{0}\right)  =\left[  p,w\right]  $.
Then $\kappa_{p_{t}}\left(  w_{t}\right)  $ is the solution of $\left(
M,S,X,N\right)  $ starting at $\left[  p,w\right]  $. 
\end{remark}

\begin{example}
\normalfont
{\bf Projectable SDEs and the horizontal-vertical factorization of diffusion operators}. In this example we show how some of the results in~\cite{liao skew-product} on the factorization of certain semielliptic differential operators on associated
bundles can be rethought in the light of the results in
Theorem \ref{teorema tan-norm associated bundles} and Corollary
\ref{corollario associated bundles}. We recall that a second order differential operator $L_{Q}%
\in\mathfrak{X}_{2}\left(  Q\right)  $ on a manifold $Q$ is called
semielliptic if any point $q \in Q$ has an open neighborhood $U$ where $L _Q $ can be locally written as
\begin{equation}
\label{infinitesimal local}
L_{Q}|_{U}=\sum_{i=1}^{s}\mathcal{L}_{Y_{i}}\mathcal{L}_{Y_{i}}+\mathcal{L}_{Y_{0}}%
\end{equation}
for some $Y_{0}$, $Y_{i}\in\mathfrak{X}\left(  U\right)  $, $i=1,...,s$. Such
a semielliptic operator can be seen as the infinitesimal generator for the laws of the  solution semimartingales
of the following stochastic system $\left(  Q,S_{Q},X,\mathbb{R}\times\mathbb{R}^{s}\right)  $ (see for instance~\cite[Theorem 1.2, page 238]{ikeda}): let $X:\mathbb{R}
_{+}\times\Omega\rightarrow\mathbb{R}\times\mathbb{R}^{s}$ be the
semimartingale%
\[
X_{t}\left(  \omega\right)  =\left(  t,B_{t}^{1}\left(  \omega\right)
,...,B_{t}^{s}\left(  \omega\right)  \right)  ,
\]
where $\left(  B^{1},...,B^{s}\right)  $ is a $s$-dimensional Brownian motion
and consider the Stratonovich operator%
\[%
\begin{array}
[c]{rrl}%
S_{Q}\left(  x,q\right)  :T_{x}\left(  \mathbb{R}\times\mathbb{R}^{s}\right)
& \longrightarrow & T_{q}U\subseteq T_{q}Q\\
\left(  u,v^{1},...,v^{s}\right)   & \longmapsto & uY_{0}+\sum_{i=1}^{s}%
v^{i}Y_{i}.
\end{array}
\]

Let now $G$ be a Lie group, $\pi:P\rightarrow Q$ a principal $G$-bundle, and
consider a manifold $W$ acted upon by $G$ via the map $l:G\times W\rightarrow W$. Let
$L_{W}\in\mathfrak{X}_{2}\left(  W\right)  $ be the semielliptic
differential operator on $W$ given by
\[
L_{W}=\sum_{i=1}^{n}\mathcal{L}_{Z_{i}}\mathcal{L}_{Z_{i}}+\mathcal{L}_{Z_{0}}%
\]
where $Z_{0},Z_{1},...,Z_{n}\in\mathfrak{X}\left(  V\right)  $ on some
$V\subseteq W$. As we just did, we will consider $L _W$ as the generator for the laws  of the solutions
of the stochastic system $\left(  W,S_{W},X^{\prime},\mathbb{R}%
^{n+1}\right)  $, where $X^{\prime}:\mathbb{R}_{+}\times\Omega\rightarrow
\mathbb{R}^{n+1}$ is a noise semimartingale constructed using the time process $t$ and
$n$ independent Brownian motions, and $S_{W}$ is the Stratonovich operator given by
\[%
\begin{array}
[c]{rrl}%
S_{W}\left(  x,w\right)  :T_{x}\left(  \mathbb{R}\times\mathbb{R}^{n}\right)
& \longrightarrow & T_{w}V\subseteq T_{w}W\\
\left(  u,v^{1},...,v^{n}\right)   & \longmapsto & uZ_{0}+\sum_{i=1}^{n}%
v^{i}Z_{i}.
\end{array}
\]
In addition, we will assume that both $L_{W}$ and  $S_{W}$ are $G$-invariant.
Let  $\widehat{\mathbf{A}}$ be a connection on the associated bundle
$M=P\times_{G}Q$ and define the Stratonovich operator $S:T\mathbb{R}%
^{n+s+1}\times M\rightarrow TM$ as%
\[
S\left(  x,\left[  p,w\right]  \right)  =T_{w}\kappa_{p}\circ S_{W}\left(
x,w\right)  +\widehat{\mathbf{A}}_{[p,w]}\circ S_{Q}\left(  x,\pi\left(
p\right)  \right)
\]
consistently with the notation introduced so far. Taking $\left(  B_{t}%
^{1},...,B_{t}^{n+s}\right)  $ a $\left(  n+s\right)  $-dimensional Brownian
motion, the stochastic system $(M,S,\widetilde{X},\mathbb{R}^{n+s+1})$ with
stochastic component $\widetilde{X}:\mathbb{R}_{+}\times\Omega\rightarrow
\mathbb{R}^{n+s+1}$ given by $\widetilde{X}_{t}\left(  \omega\right)  =\left(
t,B_{t}^{1}\left(  \omega\right)  ,...,B_{t}^{n+s}\left(  \omega\right)
\right)  $ satisfies  by construction the hypotheses of Theorem
\ref{teorema tan-norm associated bundles} and Corollary
\ref{corollario associated bundles}. The projected stochastic
system of $(M,S,\widetilde{X},\mathbb{R}^{n+s+1})$ onto $Q$ is obviously
$\left(  Q,S_{Q},X,\mathbb{R}^{s+1}\right)  $ and the one induced in the typical
fiber $W$ is $\left(  W,S_{W},X^{\prime},\mathbb{R}^{n+1}\right)  $. It is
straightforward to check that the probability laws of the solutions of $(M,S,\widetilde{X},\mathbb{R}
^{n+s+1})$ have as infinitesimal generator
\begin{equation}
L_{M}=\widetilde{L}_{Q}+L_{W}^{\ast},\label{eq tan-nor 22}%
\end{equation}
where $\widetilde{L}_{Q}$ is what Liao~\cite{liao skew-product} calls the
\textbf{horizontal lift} of $L_{Q}$ and $L_{W}^{\ast}$ the \textbf{vertical
operator} induced by $L_{W}$. 

Many of the results presented in \cite{liao
skew-product} about the factorization (\ref{eq tan-nor 22}) of semielliptic
operators on associated bundles and their related diffusions can be understood
from the perspective of stochastic systems and stochastic differential
equations that we have adopted here using
Theorem \ref{teorema tan-norm associated bundles} and Corollary
\ref{corollario associated bundles}. 
In order to illustrate this point consider the following
result in Liao's article about Riemannian submersions (see also
\cite{elworthy kendall}): let $\left(  M,\eta\right)  $ be a complete
Riemannian space with Riemann metric tensor $\eta$ and let $\bar{\pi}:M\rightarrow Q$ be
a Riemannian submersion with totally geodesic fibers. In this setup, $\bar{\pi
}:M\rightarrow Q$ is an associated bundle whose structure group $G$ is the
group of isometries of the standard fiber $W:=\bar{\pi}^{-1}\left(
q_{0}\right)  $ for some $q_{0}\in Q$~\cite{hermann}. Indeed, it can be checked that all the fibers of
$\bar{\pi}:M\rightarrow Q$ are isometric, so we can take any of them as a
standard fiber, and that $G$ has finite dimension~\cite[Remark 1.10, page 185]{Illi. Jour.
Math.}. Let
$\pi:P\rightarrow Q$ be the corresponding principal bundle. Additionally,
since $\kappa_{p}:W\rightarrow\bar{\pi}^{-1}\left(  q\right)  $ is an isometry
for any $p\in P$, the restriction $\eta_{\bar{\pi}^{-1}(q)}$ of the metric
$\eta$ to $\bar{\pi}^{-1}(q)$ may be considered as induced from the metric
$\eta_{\bar{\pi}^{-1}(q_{0})}$ of $W$ by $\kappa_{p}$ which, in addition, is
invariant by $G$. Then,
\begin{equation}
\Delta_{M}=\widetilde{\Delta}_{Q}+\Delta_{W}^{\ast}\label{eq tan-nor 23}
\end{equation}
where $\Delta_{Q}$ is the Laplacian on $Q$ and $\Delta_{W}$ the Laplacian on
$W$ (\cite[Proposition 3]{liao skew-product}). As a consequence of
(\ref{eq tan-nor 23}), if $\Gamma_{t}$ is a $M$-valued Brownian motion associated to
the Laplacian $\Delta_{M}$ on $M$ then $\bar{\pi}\left(  \Gamma_{t}\right)  $
is a Brownian motion on $Q$ with generator $\Delta_{Q}$ (see also
\cite[Theorem 10E]{elworthy}).
Let now $A\in\Omega^{1}\left(  P,\mathfrak{g}\right)  $ be the principal
connection on $P$ whose associated connection $\mathbf{A}$ on $\bar{\pi
}:M\rightarrow Q$ is such that $\operatorname*{Hor}_{m}^{\perp}
=\operatorname*{Ver}_{m}$ for any $m\in M$, that is, the horizontal subspace $\operatorname*{Hor}%
_{m}\subset T_{m}M$ of $\mathbf{A}$ is the orthogonal complement of
$\operatorname*{Ver}_{m}$, $m\in M$. Then, if $p_{t}$ denotes the horizontal
lift of $\bar{\pi}\left(  \Gamma_{t}\right)  $ to $P$ with respect to $A$ then
$\kappa_{p_{t}}^{-1}\left(  \Gamma_{t}\right)  $ is a Brownian motion on $W$
with generator $\Delta_{W}^\ast $ \cite[Propositions 6]{liao skew-product}.
\end{example}

\section{The Hamiltonian case}
\label{The Hamiltonian case}

Hamiltonian dynamical systems are a class of differential equations in the
non-stochastic deterministic context in which reduction techniques have been much
developed. This is mainly due to their central role in mechanics and applications to
physics and also to the added value that symmetries usually have in this category. As
we saw in Proposition~\ref{law conservation isotropy} the
symmetries of a stochastic differential equation bring in their wake certain
invariance properties of its flow that have to do with the preservation of the
isotropy type submanifolds. Symmetric Hamiltonian deterministic systems also
preserve isotropy type submanifolds but they usually exhibit additional invariance
features caused by the presence of symmetry induced {\bfi  first integrals} or
{\bfi  constants of motion}, usually encoded as components of a {\bfi  momentum
map}. 

The goal in this section is to show that the reduction and reconstruction techniques
that have been developed for deterministic Hamiltonian dynamical systems can be
extended to the stochastic Hamiltonian systems that have been introduced
in~\cite{lo} as a generalization of those in~\cite{bismut 81} and that we now
briefly review. The reader is encouraged to look at the original
references~\cite{bismut 81,lo}  for more details. In the following paragraphs we
will assume certain familiarity with standard deterministic Hamiltonian systems and
reduction theory (see for instance~\cite{fom, hsr} and references therein).

Let $(M,\{\cdot ,\cdot \})$ be a finite dimensional Poisson manifold,
$X:\mathbb{R}_{+}\times \Omega
\rightarrow V$ a continuous semimartingale that takes values on the vector space $V$
with $ X_{0}=0$, and let $h:M\rightarrow V^{\ast }$ be a smooth function with values
in $V ^\ast $, the dual of $V$. Let
$\{\epsilon ^{1},\ldots ,\epsilon ^{r}\}$ be a basis of $ V^{\ast }$ and let
$h _1, \ldots, h _r \in C^\infty(M) $ be such that $h=\sum_{i=1}^{r}h_{i}%
\epsilon ^{i}$. The
{\bfi  stochastic Hamiltonian system} associated to $h$ with {\bfi  stochastic
component}
$X 
$ is the stochastic differential equation
\begin{equation}
\delta \Gamma ^{h}=H(X,\Gamma )\delta X  \label{equations differential form}
\end{equation}
defined by the Stratonovich operator
$H(v,z):T_{v}V\rightarrow T_{z}M$  defined by%
\begin{equation}
H(v,z)\left( u\right) :=\sum_{i=1}^{r}\left\langle \epsilon
^{i},u\right\rangle X_{h_{i}}\left( z\right) ,
\label{Hamiltonian equations differential form}
\end{equation}%
where $X_{h_{i}}$ is the Hamiltonian vector field associated to $h_{i}\in
C^{\infty }\left( M\right) $. In this case, the dual Stratonovich operator $%
H^{\ast }(v,z):T_{z}^{\ast }M\rightarrow T_{v}^{\ast }V$ of $H(v,z)$ is
given by $H^{\ast }(v,z)(\alpha _{z})=-\mathbf{d}h(z)\cdot B^{\sharp
}(z)(\alpha _{z})$, where $B^{\sharp }:T^{\ast }M\rightarrow TM$ is the
vector bundle map naturally associated to the Poisson tensor $B\in \Lambda
^{2}(M)$ of $\{\cdot ,\cdot \}$ and $\mathbf{d}h=\sum_{i=1}^{r}\mathbf{d}%
h_{i}\otimes \epsilon ^{i}$. We will usually summarize this construction 
by saying that  $\left( M,\{ \cdot , \cdot \},h,X\right) $
is a {\bfseries\itshape stochastic Hamiltonian system}. A case to which we will
dedicate particular attention is the one in which the
Poisson manifold
$(M,\{\cdot ,\cdot
\})$ is actually symplectic with symplectic form $\omega$ and the bracket $\{\cdot ,\cdot
\} $ is
obtained from $\omega$ via the expression $\{f,h\}= \omega(X _f, X _h)$, $f,h \in
C^\infty(M) $.

\subsection{Invariant manifolds and conserved quantities of a stochastic Hamiltonian
system}

As we already said, the presence of symmetries in a Hamiltonian system forces the
appearance of invariance properties that did not use to occur for arbitrary
symmetric dynamical systems. Before we proceed with the study of those conservation
laws in the stochastic Hamiltonian case, we extract some conclusions on invariant
manifolds that can be obtained from Proposition~\ref{invariant manifolds
proposition} in that situation.

\begin{proposition}
Let $\left( M,\{ \cdot , \cdot \},h:M \rightarrow V ^\ast ,X\right) $ be a stochastic
Hamiltonian system. Let$
\{\epsilon^{1},\ldots,\epsilon^{r}\}$ be a basis of $V^{\ast}$ and write $
h=\sum_{i=1}^{r}h_{i}\epsilon^{i}$. Consider the following situations:

\begin{description}
\item[(i)] Suppose that $M$ is symplectic (respectively, Poisson) and let $z
\in M $ be such that $\mathbf{d}h(z)= 0$ (respectively, $X _{h _{i}}(z)=0 $,
for all $i\in\{1, \ldots,r\}$). Then, the Hamiltonian semimartingale $%
\Gamma^{h} $ with constant initial condition $\Gamma_{0}(\omega)=z $%
, for all $\omega\in\Omega$, is an equilibrium, that is $\Gamma^{h}=
\Gamma_{0}$.

\item[(ii)] Let $S _{1}, \ldots, S _{r}$ be submanifolds of $M$ with
transverse intersection $S:=S _{1}\cap\ldots\cap S _{r} $, such that $X _{h
_{i}}(z _{i})\in T _{z _{i}}S _{i} $, for all $z _{i} \in S _{i} $ and $i
\in\{1, \ldots,r\} $. Then $S $ is an invariant submanifold of the
stochastic Hamiltonian system $\left( M,\{ \cdot , \cdot \},h:M \rightarrow V ^\ast
,X\right) $.

\item[(iii)] The symplectic leaves of $(M, \{ \cdot, \cdot\} )$ are
invariant submanifolds of the stochastic Hamiltonian system $\left( M,\{ \cdot ,
\cdot \},h:M \rightarrow V ^\ast ,X\right) $.
\end{description}
\end{proposition}

\noindent\textbf{Proof.\ \ }It is a direct consequence of Proposition~\ref
{invariant manifolds proposition} and of the fact that the Stratonovich
operator is given by $H(v, z)(u):=\sum_{i=1}^{r}\langle
\epsilon^{j},u\rangle X_{h _{j}}(z) $. In \textbf{(i)} the hypothesis $%
\mathbf{d}h (z)=0 $ implies in the symplectic case that $X _{h _{i}}(z)=0 $,
for all $i\in\{1, \ldots,r\}$. Hence, both in the symplectic and in the
Poisson cases $H(v, z)=0 $ and hence by Proposition~\ref{invariant manifolds
proposition}, the point $z $ is an invariant submanifold and consequently an
equilibrium. For \textbf{(ii)} it suffices to recall that the transversality
hypothesis implies that $T _{z}S=T _{z} S _{1}\cap\ldots\cap T _{z} S _{r} $%
, for any $z \in S $. \textbf{(iii)} follows from the fact that the tangent
space to the symplectic leaves is spanned by the Hamiltonian vector fields
and hence $\mbox{\rm Im} \left( H(v,z)\right) \subset T _{z} \mathcal{L}_{z}$%
, for any $z\in M $ and any $v \in V $, with $\mathcal{L} _{z} $ the
symplectic leaf that contains the point $z $. \quad$\blacksquare$

\medskip

In the Hamiltonian case, most of the invariant manifolds of a system come as the
level sets of a conserved quantity (also called first integral) of the motion. In the
next definition we come back for a second to the case of general stochastic
differential equations.

\begin{definition}
Let $M$ and $N$ be two manifolds, let $X: \mathbb{R}_+ \times \Omega \rightarrow 
N$ be a $N$-valued semimartingale, and let $S:TN \times  M \rightarrow TM $ be a
Stratonovich operator. A function $f \in  C^\infty(M) $  is said to be a {\bfi 
conserved quantity} (respectively {\bfi   strongly conserved quantity}) of the
stochastic differential equation associated to $X$ and $S$ when for any solution
semimartingale $\Gamma $ we have that $f (\Gamma)= f (\Gamma _0)$ (respectively,
when $S^{\ast}\left(x,z\right) \left( \mathbf{d}f(z)\right) =0$, for any $x \in N $,
$z
\in M $).
\end{definition}

It is immediate to check that any strongly conserved quantity is a conserved
quantity. The concept of strongly conserved quantity can be equally defined for
Schwartz operators. Indeed, it can be shown that if
$S(x,z):T_{x}N%
\rightarrow T_{z}M$ is a Stratonovich operator and $\mathcal{S}\left(
x,z\right) :\tau_{x}N\rightarrow\tau_{z}M$ is the Schwartz operator induced by $S$,
then
$f\in C^{\infty}\left( M\right) $ is a strongly conserved conserved quantity for $S$
if and only if
\begin{equation}
\label{strongly conserved itos}
\mathcal{S}^{\ast}\left( x,z\right) \left( d_{2}f(z)\right) =0.
\end{equation}
for any $z\in M$ and any $x\in N$. We recall that the second order one-form $d_{2}f\in \Omega_2 (M)$
is defined as
$d_{2}f\left( L\right) \left( z\right) =L\left[ f\right] \left( z\right) $, for any
$ L\in\tau_{z}M$.

We now go back to the Hamiltonian category. Hamiltonian conserved quantities have an
interesting partial characterization in terms of Poisson commutation relations with
the components of the Hamiltonian function that the reader can find as Proposition
2.11 of~\cite{lo}. In the case of strongly conserved quantities the situation is much
simpler, as the next proposition shows.

\begin{proposition}
Let $\left( M,\{ \cdot , \cdot \},h:M \rightarrow V ^\ast ,X\right) $ be a stochastic
Hamiltonian system. Let $
\{\epsilon^{1},\ldots,\epsilon^{r}\}$ be a basis of $V^{\ast}$ and write $
h=\sum_{i=1}^{r}h_{i}\epsilon^{i}$.
Consider the Stratonovich operator $H$ given by (\ref{Hamiltonian
equations differential form}). Then, $f\in C^{\infty}\left( M\right) $ is a
strongly conserved quantity of $H$ if and only if $\left\{ f,h_{i}\right\} =0$ for
all $i=1,...,r$.
\normalfont%
(\cite{lo}).
\end{proposition}

\noindent\textbf{Proof.\ \ }
Let $f\in C^{\infty}\left( M\right) $, $v\in V$, $z\in M,$ and $u\in T_{v}V$. By
(\ref{Hamiltonian equations differential form}),%
\begin{equation*}
\left\langle H^{\ast}\left( v,z\right) \left( \mathbf{d}f(z)\right)
,u\right\rangle =\left\langle \mathbf{d}f(z),H\left( v,z\right) \left( u\right)
\right\rangle =\sum_{i=1}^{r}\left\{ f,h_{i}\right\} \left( z\right)
\left\langle u,\epsilon^{i}\right\rangle . 
\end{equation*}
Since $u\in T_{v}V$ is arbitrary, $H^{\ast}\left( v,z\right) \left( \mathbf{d%
}f(z)\right) =0$ if and only if $\left\{ f,h_{i}\right\} \left( z\right) =0$.
\quad $\blacksquare$

\medskip

We now concentrate on the conserved quantities that one can associate to the
invariance of a Hamiltonian system with respect to a group action. We recall that
given a Lie group $G$ acting on the Poisson manifold $(M, \{ \cdot , \cdot \}) $
(respectively, symplectic $(M , \omega)$) via the map $ \Phi:G \times M \rightarrow M
$, we will say that the action is {\bfi  canonical} when for any $g \in G $ and $f,h
\in C^\infty(M) $, $\{ f,h\} \circ \Phi_g=\{ \Phi ^\ast _g f ,\Phi ^\ast _g h\} $
(respectively, $\Phi^\ast _g \omega= \omega$). In this context, we will say that the
Hamiltonian system $\left( M,\{ \cdot , \cdot \},h:M \rightarrow V ^\ast ,X\right) $
is $G$-{\bfi  invariant} whenever the $G$-action on $M$ is canonical and the
Hamiltonian function $h:M \rightarrow V ^\ast $ is $G$-invariant. Notice that the
invariance of $h$  and the canonical character of the action imply that the
associated Stratonovich operator $H$ is also $G$-invariant. Indeed, Let $
\{\epsilon^{1},\ldots,\epsilon^{r}\}$ be a basis of $V^{\ast}$ and write $
h=\sum_{i=1}^{r}h_{i}\epsilon^{i}$; if $h$ is $G$-invariant, then so are the
components $h_i$, $i \in \{1, \ldots, r \} $, that is $h _i \in C^\infty(M)^G $,
and hence, for any $g \in G$ we have that $T \Phi _g \circ X _{h _i}= X _{h _i}
\circ \Phi_g$, which implies that $H(v,z)\left( u\right) :=\sum_{i=1}^{r}\left\langle \epsilon
^{i},u\right\rangle X_{h_{i}}\left( z\right) $is $G$-invariant.

Now suppose that $M$ is a Poisson manifold $(M,\{\cdot,\cdot\})$ acted
properly and canonically upon by a Lie group $G$. We also recall that the {%
\bfseries\itshape optimal momentum map}~\cite{optimal} $\mathcal{J}%
:M\rightarrow M/D_{G}$ of the $G$-action on $(M,\{\cdot,\cdot\})$ is the
projection onto the leaf space of the integrable distribution $D_{G}\subset
TM$ (in the generalized sense of Stefan-Sussmann) given by $%
D_{G}:=\{X_{f}\mid f\in C^{\infty}(M)^{G}\}$.

\begin{proposition}
\label{several noether} Let $(M,h,X,V)$ be a standard Hamiltonian system
acted properly and canonically upon by a Lie group $G$ via the map $%
\Phi:G\times M\rightarrow M$. Suppose that $h:M\rightarrow V^{\ast}$ is a $G$%
-invariant function.

\begin{description}
\item[(i)] \textbf{Law of conservation of the isotropy}: The isotropy type
submanifolds $M_{I}$ are invariant submanifolds of the stochastic
Hamiltonian system associated to $h$ and $X$, for any isotropy subgroup $%
I\subset G$.

\item[(ii)] \textbf{Noether's Theorem}: If the $G$-action on $(M,\{\cdot
,\cdot\})$ has a momentum map associated $\mathbf{J}:M\rightarrow \mathfrak{g%
}^{\ast}$ then its level sets are left invariant by the stochastic
Hamiltonian system associated to $h$ and $X$. Moreover, its components are
conserved quantities.

\item[(iii)] \textbf{Optimal Noether's Theorem}: The level sets of the
optimal momentum map $\mathcal{J}:M\rightarrow M/D_{G}$ are left invariant
by the stochastic Hamiltonian system associated to $h$ and $X$.
\end{description}
\end{proposition}

\noindent\textbf{Proof.\ \ } \textbf{(i)} As we already saw, the $G$-invariance of
$h$ implies that $H(v,z)(u):=\sum_{i=1}^{r}\langle\epsilon^{j},u\rangle
X_{h_{j}}(z)$ is 
$G$-invariant. The statement follows from
Proposition~\ref{law conservation isotropy}. \textbf{(ii)} Let $%
\xi\in\mathfrak{g}$ be arbitrary and let $\mathbf{J}^{\xi}:=\langle \mathbf{J%
},\xi\rangle\in C^{\infty}(M)$ be the corresponding component. The $G$%
-invariance of the components $h_{i}$ of the Hamiltonian implies that $\{%
\mathbf{J}^{\xi},h_{i}\}=-\mathbf{d}h_{i}\cdot\xi_{M}=0$, where $\xi_{M}\in%
\mathfrak{X}(M)$ is the infinitesimal generator associated to the element $%
\xi$. By formula (2.8) in \cite{lo} we have that 
\begin{equation*}
\mathbf{J}^{\xi}(\Gamma^{h})-\mathbf{J}^{\xi}(\Gamma_{0})=\sum_{j=1}^{r}\int%
\{\mathbf{J}^{\xi},h_{j}\}\,\delta X^{j}=0, 
\end{equation*}
where $X_{j}$, $j\in\{1,\ldots,r\}$, are the components of $X$ in the basis $%
\{e_{1},\ldots,e_{r}\}$ of $V$ dual to the basis
$\{\epsilon^{1},\ldots,\epsilon^{r}\}$ of $V^{\ast}$. Since this equality
holds for any $\xi \in%
\mathfrak{g}$, we have that $\mathbf{J}(\Gamma^{h})=\mathbf{J}(\Gamma_{0})$
and the result follows. \textbf{(iii)} It is a straightforward consequence
of the construction of the optimal momentum map and Proposition \ref%
{invariant manifolds proposition}. \quad$\blacksquare$

\begin{remark}
\normalfont When the manifold $M$ is symplectic and the group action has a
standard momentum map $\mathbf{J}:M\rightarrow\mathfrak{g}^{\ast}$
associated, part \textbf{(iii)} in the previous proposition implies the
first two since it can be shown that in that situation (see \cite{optimal})
the level sets of the optimal momentum map coincide with the connected
components of the intersections $\mathbf{J}^{-1}(\mu)\cap M_{I}$, with $%
\mu\in\mathfrak{g}^{\ast}$ and $I$ an isotropy subgroup of the $G$-action on $M$.
\end{remark}

\begin{remark}
\label{Stratonovich tangente a momento}\normalfont The level sets of the
momentum map $\mathbf{J}$ may not be submanifolds of $M$ unless the $G$-action
is, in addition to proper and canonical, also free (\cite[Corollary 4.6.2]%
{hsr}). If this is the case, the relation $\{\mathbf{J}^{\xi
},h_{i}\}=-\mathbf{d}h_{i}\cdot\xi_{M}=0$, which stems from the $%
G$-invariance of $h$, implies then that $\mbox{\rm Im}\left( H(v,z)\right)
\subset T_{z}\mathbf{J}^{-1}(\mu)$, for any $z\in\mathbf{J}^{-1}(\mu)$ and
any $v\in V$. Then Proposition \ref{invariant manifolds proposition} may be
invoked to prove the invariance of the fibers $\mathbf{J}^{-1}(\mu)$ under the
stochastic Hamiltonian system associated to $H$.
\end{remark}

\subsection{Stochastic Hamiltonian  reduction and reconstruction}

The goal of this section is showing that stochastic Hamiltonian systems share with
their deterministic counterpart a good behavior with respect to symmetry
reduction. The main idea that our following theorem tries to convey to the reader is
that \emph{the symmetry reduction of a stochastic Hamiltonian system yields a
stochastic Hamiltonian system}, that is, the stochastic Hamiltonian category is
stable under reduction.

The following theorem spells out, in the simplest possible case, how to
reduce symmetric Hamiltonian stochastic systems. In a remark below we give the
necessary prescriptions to carry  this procedure out in more general situations.
The main simplifying hypothesis is the freeness of the action. We recall that 
in this
situation, the orbit space $M/G$ inherits from $M$ a Poisson structure $%
\{\cdot,\cdot\}_{M/G}$ naturally obtained by projection of that in $M$, that
is, $\{f,g\}_{M/G}\circ\pi:=\{f\circ\pi,g\circ\pi\}$, for any $f,g\in
C^{\infty}(M/G)$, with $\pi:M\rightarrow M/G$ the orbit projection.
Moreover, if $M$ is actually symplectic with symplectic form $\omega$, and the
action has a coadjoint equivariant momentum map $\mathbf{J}: M \rightarrow 
\mathfrak{g}^\ast$, then the symplectic leaves of this Poisson structure are
naturally symplectomorphic to the (connected components) of the {\bfi 
Marsden-Weinstein}
\cite{mwr} {\bfi  symplectic quotients}
$(M_{\mu}:=\mathbf{J}^{-1}(\mu)/G_{\mu},\omega_{
\mu})$, with $\mu \in\mathfrak{g}^{\ast}$ and $G_{\mu}$ the coadjoint
isotropy of $\mu$. The symplectic structure $\omega_{\mu}$ on $M_{\mu}$ is
uniquely determined by the expression $\pi_{\mu}^{\ast}\omega_{\mu}=i_{%
\mu}^{\ast}\omega$, with $i_{\mu }:\mathbf{J}^{-1}(\mu)\hookrightarrow M$
the injection and $\pi_{\mu }:\mathbf{J}^{-1}(\mu)\rightarrow\mathbf{J}%
^{-1}(\mu)/G_{\mu}$ the projection. See \cite{fom, hsr} and references
therein for a general presentation of reduction theory.

\begin{theorem}
\label{reduction theorem statement} 
Let $\left( M,\{ \cdot , \cdot \},h:M \rightarrow V ^\ast ,X\right) $ be a stochastic
Hamiltonian system that is invariant with respect to the canonical, free , and
proper action
$\Phi:G\times M\rightarrow M$ of the Lie group $G$ on $M$. 
\begin{description}
\item[(i)] \textbf{Poisson reduction}: The projection $h_{M/G}$ of the
Hamiltonian function $h$ onto $M/G $, uniquely determined by $h_{M/G}\circ\pi=h$,
with
$%
\pi:M\rightarrow M/G$ the orbit projection, induces a stochastic Hamiltonian
system on the Poisson manifold $(M/G,\{\cdot,\cdot\}_{M/G})$ with stochastic
component $X$ and whose Stratonovich operator $H_{M/G}:TV\times
M/G\rightarrow T\left( M/G\right) $ is given by 
\begin{equation}
H_{M/G}(v,\pi(z))(u)=T_{z}\pi\left( H(v,z)(u)\right)=\sum_{i=1}^{r}\langle
\epsilon^i,u\rangle X_{h _i^{M/G}}(\pi (z))
 ,\quad\text{ $u,v\in V$
and $z\in M$.}   \label{Poisson reduced Stratonovich operator}
\end{equation}
In the previous expression $
\{\epsilon^{1},\ldots,\epsilon^{r}\}$ is a basis of $V^{\ast}$,
$h_{M/G}=\sum_{i=1}^r h _i^{M/G} \epsilon^i$, and $
h=\sum_{i=1}^{r}h_{i}\epsilon^{i}$; notice that the functions $h _i^{M/G} \in
C^\infty(M/G)$ are the projections of the components $h _i\in C^\infty(M)^{G} $,
that is $h _i^{M/G}\circ \pi=h _i
$. Moreover, if
$\Gamma$ is a solution semimartingale of the Hamiltonian system associated to $H$
with initial condition $\Gamma_{0}$, then so is
$%
\Gamma_{M/G}:=\pi\left( \Gamma\right) $ with respect to $H_{M/G}$, with
initial condition $\pi(\Gamma_{0})$.
\item[(ii)] \textbf{Symplectic reduction}: Suppose that $M$ is now symplectic
and that the group action
has a coadjoint equivariant momentum map $\mathbf{J}:M\rightarrow\mathfrak{g}%
^{\ast}$ associated. Then for any $\mu\in\mathfrak{g}^{\ast}$, the function $%
h_{\mu}:M_{\mu}:=\mathbf{J}^{-1}(\mu)/G_{\mu}\rightarrow V^{\ast}$ uniquely
determined by the equality $h_{\mu}\circ\pi_{\mu}=h\circ i_{\mu}$ , induces
a stochastic Hamiltonian system on the symplectic reduced space $(M_{\mu }:=%
\mathbf{J}^{-1}(\mu)/G_{\mu},\omega_{\mu})$ with stochastic component $X$
and whose Stratonovich operator $H_{\mu}:TV\times M_{\mu}\rightarrow TM_{\mu}
$ is given by 
\begin{equation}
H_{\mu}(v,\pi_{\mu}(z))(u)=T_{z}\pi_{\mu}\left( H(v,i_{\mu}(z))(u)\right)
=\sum_{i=1}^{r}\langle
\epsilon^i,u\rangle X_{h _i^{\mu}}(\pi_\mu (z)) ,\quad
\text{$u,v\in V$ and $z\in\mathbf{J}^{-1}(\mu)$,} 
\label{symplectic reduced Stratonovich operator}
\end{equation}
where Remark \ref{Stratonovich tangente a momento} has been implicitly used.
In the previous expression, the functions $h _i^\mu \in
C^\infty(\mathbf{J}^{-1}(\mu)/G _\mu) $ are the coefficient functions in the linear
combination $h_{\mu}=\sum_{i=1}^r h _i^{\mu} \epsilon^i$ and are related to the
components $h _i \in C^\infty(M)^{G}$ of $h $ via the relation $h _i^\mu \circ
\pi_\mu=h _i \circ i _\mu $.
Moreover, if $\Gamma$ is a solution semimartingale of the Hamiltonian system
associated to $H$ with initial condition $\Gamma_{0}\subset\mathbf{J}%
^{-1}(\mu)$, then so is $\Gamma_{\mu}:=\pi_{\mu}\left( \Gamma\right) $ with
respect to $H_{\mu}$, with initial condition $\pi_\mu(\Gamma_{0})$.
\end{description}
\end{theorem}

\begin{remark}
\normalfont In the absence of freeness of the action the orbit spaces $M/G$
and $\mathbf{J}^{-1}(\mu)/G_{\mu}$ cease to be regular quotient manifolds.
Moreover, it could be that (even for free actions) there is no standard
momentum map available (this is generically the case for Poisson manifolds).
This situation can be handled by using the so called optimal momentum map 
\cite{optimal} and its associated reduction procedure \cite{symplectic
reduced}. Given that by part \textbf{(iii)} of Proposition \ref{several
noether} the fibers of the optimal momentum map are preserved by the
Hamiltonian semimartingales associated to invariant Hamiltonians one can
formulate, for any proper group action on a Poisson manifold, a theorem
identical to part \textbf{(ii)} of Theorem \ref{reduction theorem statement}
with the standard momentum map replaced by the optimal momentum map. In the
particular case of a (non-necessarily free) symplectic proper action that
has a standard momentum map associated, such result guarantees the good
behavior of the symmetric stochastic Hamiltonian systems with respect to the
singular reduced spaces in \cite{sl}; see also~\cite{cylinder reduction, cylinder
reduction singular} for the symplectic case without a standard momentum map. 
\end{remark}

\noindent\textbf{Proof of Theorem \ref{reduction theorem statement}.\ \ }
\textbf{(i)} can be proved by mimicking the proof of Theorem~\ref{reduction theorem
statement general systems} by simply taking into account the fact that the
$G$-invariance of $h$  implies that of $H$ and that for any $i \in \{1, \ldots, r \}
$, one has that $T \pi \circ X_{h _i}=X_{h _i^{M/G}} \circ \pi $. 

\noindent \textbf{(ii)} 
Expression~(\ref{symplectic reduced Stratonovich operator}) is guaranteed by the
fact that 
$X_{h_{i}^{\mu}}\circ\pi_{
\mu}=T\pi_{\mu}\circ X_{h_{i}}\circ i_{\mu}$, for any $i\in\{1,\ldots,r\}$
(see for instance~\cite[Theorem 6.1.1]{hsr}).
Let now $\Gamma$ be a solution semimartingale of the Hamiltonian system
associated to $H$ with initial condition $\Gamma_{0}\subset\mathbf{J}%
^{-1}(\mu)$. Notice first that by part \textbf{(ii)} in Proposition \ref%
{several noether}, $\Gamma\subset\mathbf{J}^{-1}(\mu)$ and hence the
expression $\Gamma_{\mu }:=\pi_{\mu}\left( \Gamma\right) $ is well defined.
In order to prove the statement, we have to check that for any one-form $%
\alpha_{\mu}\in \Omega(M_{\mu})$ 
\begin{equation*}
\int\langle\alpha_{\mu},\delta\Gamma_{\mu}\rangle=\int\langle H_{\mu}^{\ast
}(X,\Gamma_{\mu})\alpha_{\mu},\delta X\rangle. 
\end{equation*}
This equality follows in a straightforward manner from~(\ref{symplectic
reduced Stratonovich operator}). Indeed, 
\begin{equation*}
\int\langle\alpha_{\mu},\delta\Gamma_{\mu}\rangle=\int\langle\alpha_{\mu
},\delta\left( \pi_{\mu}\circ\Gamma\right) \rangle=\int\langle\pi_{\mu
}^{\ast}\alpha_{\mu},\delta\Gamma\rangle=\int\langle
H^{\ast}(X,\Gamma)\left( \pi_{\mu}^{\ast}\alpha_{\mu}\right) ,\delta
X\rangle=\int\langle H_{\mu }^{\ast}(X,\Gamma_{\mu})\alpha_{\mu},\delta
X\rangle, 
\end{equation*}
as required. \quad$\blacksquare$

\medskip

As to the reconstruction problem of solutions of  a symmetric stochastic
differential equation starting from a solution of the Poisson or symplectic reduced
stochastic differential equation, Theorem~\ref{teorema reconstruccion 1} can be
trivially modified to handle this situation. In the Poisson reduction case the
theorem works without modification and when working with a solution of the symplectic
reduced space it suffices to change the principal fiber bundle $\pi: M \rightarrow M
/G $  by $\pi_\mu: \mathbf{J}^{-1}(\mu) \rightarrow \mathbf{J}^{-1}(\mu)/ G_{\mu} $
all over.

\section{Examples}
\label{examples we see}

\subsection{Stochastic collective Hamiltonian motion}
\label{Stochastic collective Hamiltonian motion}

Our first example shows a situation in which the symplectic reduction of a
symmetric  stochastic Hamiltonian system offers, not only the advantage of cutting
its dimension, but also of \emph{making it into a deterministic system}. From the
point of view of obtaining the solutions of the system, the procedures introduced in
the previous section allow in this case the splitting of the problem into two parts:
first, the solution of a standard ordinary differential equation for the reduced
system and second, the solution of a stochastic differential equation in the group
at the time of the  reconstruction.

Let $\left( M,\omega \right) $ be a symplectic manifold, $G$ a Lie group
and
$
\Phi :G\times M\rightarrow M$ a free, proper, and canonical action.
Additionally, suppose that this action has  a coadjoint equivariant momentum map
$\mathbf{J}:M\rightarrow 
\mathfrak{g}^{\ast }$ associated. Let $h_{0}\in C^\infty(M)^{G}$ be a $G$-invariant
function and consider the deterministic Hamiltonian system with Hamiltonian function
$h_{0}$.

A function of the form $f\circ\mathbf{J}\in  C^\infty(M)$, for some $f \in 
C^\infty(\mathfrak{g}^\ast) $, is called 
{\bfi  collective}. We recall that by the Collective Hamiltonian Theorem (see for
instance~\cite{symmetry and mechanics})
\begin{equation}
\label{collective Hamiltonian vector field}
X_{f \circ \mathbf{J}}(z)= \left(\frac{ \delta f }{\delta \mu}  \right)_M (z), \quad
z \in M, \, \mu= \mathbf{J} (z),
\end{equation}
where the functional derivative $\frac{ \delta f }{\delta \mu} \in  \mathfrak{g}$ is
the unique element such that for any $\nu\in \mathfrak{g}^\ast$, $Df(\mu) \cdot
\nu=\langle \nu, \frac{ \delta f }{\delta \mu}\rangle $. A straightforward
consequence of~(\ref{collective Hamiltonian vector field}) is that the $G$-invariant
functions, in particular $h _0 $,  commute with the collective functions. Indeed, if
$h \in C^\infty(M)^{G}$, then for any $z \in M $,
\begin{equation*}
\{h, f \circ \mathbf{J}\} (z)= \mathbf{d}h (z) \cdot X_{f \circ \mathbf{J}}(z)=
\mathbf{d}h (z) \cdot \left(\frac{ \delta f }{\delta \mu}  \right)_M (z)=0.
\end{equation*}
Collective functions play an important role to prove the complete
integrability of certain dynamical systems (see \cite{guillemin
integrability}). Moreover, some relevant physical systems may be described
using collective Hamiltonian functions. In that case, the (deterministic)
equations of motion exhibit special features and, in some favorable cases,
may be partially integrated using geometrical arguments (see \cite{guillemin
collective motion}). The aim of this example is to study stochastic
perturbations of deterministic symmetric mechanical systems introduced by means of
collective Hamiltonians.

Let $Y:\mathbb{R}_+\times\Omega\rightarrow\mathbb{R}^{r}$ be a $\mathbb{R}^{r}$%
-valued continuous semimartingale and $\left\{ f_{1},...,f_{r}\right\}
\subset C^{\infty}\left( \mathfrak{g}^{\ast}\right) $ a finite family of $
{\rm Ad}^{\ast}_G$-invariant functions on $\mathfrak{g}^{\ast}$. The coadjoint equivariance of the momentum map and the $
{\rm Ad}^{\ast}_G$-invariance of the functions allows us to construct the following
$G$-invariant Hamiltonian function%
\begin{equation*}
\begin{array}{rcl}
h:M & \longrightarrow & \mathbb{R\times R}^{r} \\ 
m & \longmapsto & \left( h_{0}\left( m\right) ,\left( f_{1}\left( \mathbf{J}%
\left( m\right) \right) ,...,f_{r}\left( \mathbf{J}\left( m\right) \right)
\right) \right).
\end{array}
\end{equation*}
Let $X$ be the continuous semimartingale 
\begin{equation*}
\begin{array}{rcl}
X:\mathbb{R}_+\times\Omega & \longrightarrow & \mathbb{R_+\times R}^{r} \\ 
\left( t,\omega\right) & \longmapsto & \left( t,Y_{t}\left( \omega\right)
\right) .%
\end{array}
\end{equation*}
Consider the stochastic Hamiltonian system $\left( M,\omega,h,X\right)$ which is, by
construction,
$G$-invariant. Noether's theorem (Proposition~\ref{several noether} {\bf
(ii)}) guarantees that the level sets of $\mathbf{J}  $  are left 
invariant by the solution semimartingales of $\left( M,\omega,h,X\right)$. As to the
reduction of this system, its main feature is that if we apply to it the reduction
scheme introduced in Theorem~\ref{reduction theorem statement} {\bf (ii)}, for any
$\mu \in  \mathfrak{g}^\ast$, the reduced stochastic Hamiltonian system
$\left( M_{\mu},\omega _\mu, h_{\mu},X\right) $ is such that
\begin{equation*}
h_{\mu}\circ\pi_{\mu}=h_{0}\circ i_{\mu}, 
\end{equation*}
since $\mathbf{J}$, and hence the functions $f_{i}\circ\mathbf{J}$, are constant on
the level sets
$\mathbf{ J}^{-1}\left( \mu\right) $, for any $i=1,...,r$. Consequently, the
reduced system 
$\left( M_{\mu},\omega _\mu, h_{\mu},X\right) $ is equivalent to the deterministic 
Hamiltonian system $
\left( M_{\mu},\omega _\mu,h_{\mu}\right) $.
In other words, the reduced system obtained from $\left( M,\omega,h,X\right) $
coincides with the one obtained in deterministic mechanics by symplectic reduction
of $\left( M,h_{0},t,\mathbb{R}_+\right) $. Thus, \emph{we have perturbed
stochastically a symmetric mechanical system preserving its symmetries and without
changing the deterministic behavior of its corresponding reduced system}.

\begin{remark}
\normalfont%
If we want to perturb the deterministic Hamiltonian system associated to $h_{0}$ 
with the only prescription that the level set $\mathbf{J}^{-1}\left( \mu \right)
$ is left invariant, for a given value
$
\mu \in \mathfrak{g}^{\ast }$, we can weaken the requirement on the  ${\rm Ad}^{\ast
}_G$-invariance of the functions $f_{i}\in C^{\infty }\left( \mathfrak{g}^{\ast
}\right) 
$, $i=1,...,r$. Indeed, if we just ask that $\delta f_{i}/ \delta\mu \in
\mathfrak{g}_\mu $, we then have that $X _{h _0}(z),X_{f _1 \circ \mathbf{J}}(z),
\ldots ,X_{f _r \circ \mathbf{J}}(z) \in T _z \mathbf{J}^{-1}(\mu)$, for any $z \in
\mathbf{J}^{-1}(\mu)$. The required invariance property follows then
from~(\ref{collective Hamiltonian vector field}) and Proposition~\ref{invariant
manifolds proposition}.  
\end{remark}

\begin{remark}
\normalfont
In this example, the reduction-reconstruction scheme provides a global decomposition of the system $\left( M,\omega,h,X\right) $ into its deterministic and stochastic parts. If one is willing to work only locally, this splitting could be carried out without reduction in the neighborhood of any point in phase space, given that as $\{h _0, f _i\circ  \mathbf{J}\}=0 $, for any $i \in \{ 1, \ldots , r \} $, then $[X_{h _0}, X_{f _i \circ \mathbf{J}}]=0 $.
\end{remark}

\subsection{Stochastic mechanics on Lie groups} 
\label{Stochastic mechanics on Lie groups} 

The presence of mechanical systems whose phase space is the cotangent bundle of a
Lie group is widespread. Besides the importance that this general case has
in specific applications it is also very useful at the time of illustrating
some of the theoretical developments in this paper since most of the
constructions that we presented admit very explicit characterizations.
We start by recalling the main features of
(deterministic) Hamiltonian systems over Lie groups. The reader interested in
further details is encouraged to check with \cite {fom,symmetry and mechanics} and
references therein.

Let $G$ be a Lie group. The tangent bundle $TG$ of $G$ is trivial since it is
isomorphic to the product $G\times \mathfrak{g}$, where $\mathfrak{g=}T_{e}G$ is the
Lie algebra of $G$ and $e\in G$ is the identity element. The identification $%
TG=G\times \mathfrak{g}$ is usually carried out by means of two
isomorphisms, denoted by $\lambda $ and $\rho $ and induced by left and
right translations on $G$, respectively. More specifically, let
$\lambda :TG\rightarrow G\times \mathfrak{g}$ be the map given by $\lambda \left(
v\right) =\left( g,T_{g}L_{g^{-1}}\left( v\right) \right) $, where $g=\tau
_{G}\left( v\right) $ with $\tau _{G}:TG\rightarrow G$  the natural
projection. On the other hand, $\rho :TG\rightarrow G\times \mathfrak{g}$ is
defined by $\rho \left( v\right) =\left( g,T_{g}R_{g^{-1}}\left( v\right)
\right) $. We refer to the image of $\lambda $ as 
{\bfi body coordinates}
and to the image of $\rho $ as 
{\bfi space coordinates}%
. The cotangent bundle $T^{\ast }G$ is also trivial and
isomorphic to $G\times \mathfrak{g}^{\ast }$. We can introduce 
{\bfi body coordinates}
and 
{\bfi space coordinates}
on $T^{\ast }G$ by $\bar{\lambda}\left( \alpha\right) =\left( g,T_{e}^{\ast
}L_{g}\left( \alpha\right) \right) \in G\mathfrak{\times g}^{\ast }$ and $\bar{%
\rho}\left( \alpha\right) =\left( g,T_{e}^{\ast }R_{g}\left( \alpha\right) \right) $
respectively, where $g=\pi _{G}\left( \alpha\right) $ and $\pi _{G}:T^{\ast
}G\rightarrow G$ is the canonical projection. The transition from body to
space coordinates is as follows:%
\begin{align*}
\left( \rho \circ \lambda ^{-1}\right) \left( g,\xi \right) & =\rho \left(
g,T_{e}L_{g}\left( \xi \right) \right) =\left( g,T_{g}R_{g^{-1}}\circ
T_{e}L_{g}\left( \xi \right) \right) =\left( g,{\rm Ad}%
_{g}\left( \xi \right) \right)  \\
\left( \bar{\rho}\circ \bar{\lambda}^{-1}\right) \left( g,\mu \right) &
=\rho \left( g,T_{g}^{\ast }L_{g^{-1}}\left( \mu \right) \right) =\left(
g,T_{e}^{\ast }R_{g}\circ T_{g}^{\ast }L_{g^{-1}}\left( \mu \right) \right)
=\left( g,{\rm Ad}_{g^{-1}}^{\ast }\left( \mu \right) \right) ,
\end{align*}%
for any $\left( g,\xi \right) \in G\times \mathfrak{g}$ and any $\left(
g,\mu \right) \in G\times \mathfrak{g}^{\ast }$.
The group action of $G$ by left
or right translations can be lifted to both $TG$ and $T^{\ast }G$. We will
denote by $\Phi _{L}:G\times TG\rightarrow TG$ and $\bar{\Phi}_{L}:G\times T^{\ast
}G\rightarrow T^{\ast }G$ the lifted action of left translations on the tangent and
cotangent bundle respectively, and by $\Phi _{R}:G\times TG\rightarrow TG$ and
$\bar{\Phi} _{R}:G\times T^{\ast }G\rightarrow T^{\ast }G$ the lifted actions of
right translations. The lifted actions have particularly simple expressions in
suitable body or space coordinates. Indeed, it is more convenient to express $\Phi
_{L}$ and $%
\bar{\Phi}_{L}$ in body coordinates, where for any $g,h \in G $, $\xi\in
\mathfrak{g}$, and $\mu\in \mathfrak{g}^\ast$, 
\begin{align*}
\left( \Phi _{L}\right) _{g}\left( h,\xi \right) & =\left( \lambda \circ
TL_{g}\circ \lambda ^{-1}\right) \left( h,\xi \right) =\left( gh,\xi \right), 
\\
\left( \bar{\Phi}_{L}\right) _{g}\left( h,\mu \right) & =\left( \bar{\lambda}%
\circ T^{\ast }L_{g^{-1}}\circ \bar{\lambda}^{-1}\right) \left( h,\mu
\right) =\left( g^{-1}h,\mu \right).
\end{align*}%
As to  $\Phi _{R}$ and $\bar{\Phi}_{R}$,  space coordinates are particularly
convenient; for any $g,h \in G $, $\zeta \in \mathfrak{g}$, and $\alpha\in
\mathfrak{g}^\ast$, 
\begin{align*}
\left( \Phi _{R}\right) _{g}\left( h,\zeta \right) & =\left( \rho \circ
TR_{g}\circ \rho ^{-1}\right) \left( h,\zeta \right) =\left( hg,\zeta
\right)  \\
\left( \bar{\Phi}_{R}\right) _{g}\left( h,\alpha \right) & =\left( \bar{\rho}%
\circ T^{\ast }R_{g^{-1}}\circ \bar{\rho}^{-1}\right) \left( h,\alpha
\right) =\left( hg^{-1},\alpha \right) .
\end{align*}

The actions $\bar{\Phi}_{L}$ and $\bar{\Phi}_{R}$, being the cotangent lifted actions
to $T^{\ast }G$ of an action on $G$, have canonical momentum maps $\mathbf{J}%
_{L}:T^{\ast }G\rightarrow \mathfrak{g}^{\ast }$ and $\mathbf{J}_{R}:T^{\ast
}G\rightarrow \mathfrak{g}^{\ast }$, respectively, when we endow $T ^\ast G $ with
its canonical symplectic form. Let
$\theta
\in
\Omega ^{1}\left( T^{\ast }G\right) $ be the Liouville canonical one-form on $
T^{\ast }G$. Then, $\mathbf{J}_{L}$ and $\mathbf{J}_{R}$ are given by 
\[
\left\langle \mathbf{J}_{L}\left( z_{g}\right) ,\xi \right\rangle 
=\langle z _g, \left( \xi \right) _{G}^{L} (g)\rangle,  \qquad
\left\langle \mathbf{J}_{R}\left( z_{g}\right) ,\xi \right\rangle 
=\langle z _g, \left( \xi \right) _{G}^{R} (g)\rangle, 
\]
for any $z_{g}\in T_{g}^{\ast }G$ and any $\xi \in \mathfrak{g}$. Here $
\left( \xi \right) _{G}^{L}\in \mathfrak{X}\left( G\right) $ (respectively $\left(
\xi \right) _{G}^{R}\in \mathfrak{X}\left( G\right) $) denotes the
infinitesimal generator associated to $\xi \in \mathfrak{g}$ by the left (respectively right)action of $G$  on itself. This expression clearly shows that $\mathbf{%
J}_{L}$ is right-invariant and $\mathbf{J}_{R}$ left-invariant. Observe that 
$\mathbf{J}_{L}={\rm Ad}_{g^{-1}}^{\ast }\circ \mathbf{J}_{R}$.
For example, in body coordinates, these momentum maps have the following
expressions (\cite[Theorem 4.4.3]{fom})%
\begin{equation}
\left( \mathbf{J}_{L}\right) _{B}\left( \left( g,\mu \right) \right) =%
{\rm Ad}_{g^{-1}}^{\ast }\left( \mu \right) \text{ \ \ and \ \ }%
\left( \mathbf{J}_{R}\right) _{B}\left( \left( g,\mu \right) \right) =\mu .
\label{momentum cotangent body}
\end{equation}

In this context, the classical results on symplectic and Poisson reduction that we
have described in the previous section  admit a particularly explicit
formulation. In all that follows we will suppose that the action with respect to
which we are reducing is lifted left translations. Using body coordinates, it is easy
to see that in this case the Poisson reduced space $T ^\ast G/G $  coincides with
the dual of the Lie algebra $\mathfrak{g}^\ast$ endowed with the {\bfi 
Lie-Poisson} structure given by 
\begin{equation*}
\left\{ f_{1},f_{2}\right\} _{\mathfrak{-}}^{\mathfrak{g}^{\ast }}\left(
\mu \right) =-\left\langle \mu , \left[ \frac{\delta f_{1}}{\delta \mu },\frac{
\delta f_{2}}{\delta \mu }\right] \right\rangle,
\end{equation*}
for any $\mu\in \mathfrak{g}^\ast$ and $f _1, f _2 \in C^\infty(\mathfrak{g}^\ast)$.
The symplectic reduced spaces $\mathbf{J}_L^{-1}(\mu)/G _\mu $ are naturally
symplectomorphic to the symplectic leaves of the Lie-Poisson structure on
$\mathfrak{g}^\ast$, that is, the coadjoint orbits endowed with the so-called {\bfi 
Kostant-Kirillov-Souriau} symplectic form $\omega_\mu^- $:
\begin{equation*}
\omega_\mu^-(\mu)(\xi_{\mathfrak{g}^\ast}(\mu), \eta_{\mathfrak{g}^\ast}(\mu))
=\omega_\mu^-(\mu)(- \mbox{\rm ad} ^\ast _\xi \mu,- \mbox{\rm ad} ^\ast _\eta \mu
)=-\langle
\mu,[\xi, \eta]\rangle.
\end{equation*}
Let now $V$ be a vector space, $X: \mathbb{R}_+ \times \Omega \rightarrow V $ a
continuous semimartingale, and $h :T ^\ast G \rightarrow V ^\ast  $ a smooth map
invariant under the lifted left translations of $G $ on $T ^\ast G $. If we use body
coordinates and we visualize $T ^\ast G $ as the product $G \times \mathfrak{g}^\ast
$, the invariance of $h:G \times \mathfrak{g}^\ast \rightarrow V ^\ast $ allows us to
write it as
$h=\sum_{i=1}^{r} h _i
\epsilon^i$, where $\{ \epsilon ^1, \ldots, \epsilon^r\} $ is a basis of $V ^\ast  $
and $h _1, \ldots, h _r \in C^\infty(\mathfrak{g}^\ast)$. Let $\{ e _1, \ldots, e
_r\} $ be the dual basis of $V$ and write $X=\sum _{i=1}^r X ^i e _i $. Using the
left trivialized expression of the Hamiltonian vector fields in the deterministic
case (see~\cite[Theorem 6.2.5]{hsr}) it is easy to see that the stochastic
Hamiltonian equations in this setup are
\begin{equation}
\label{stochastic trivialized vector field}
\delta\Gamma^h=\sum _{i=1}^r \left(T _eL_{ \Gamma^G}\left(\frac{\delta h _i}{\delta
\Gamma ^{\mathfrak{g}^\ast}}\right), \mbox{\rm ad} ^\ast _{\frac{\delta h _i}{\delta
\Gamma ^{\mathfrak{g}^\ast}}}\Gamma ^{\mathfrak{g}^\ast}
\right) \delta X ^i
\end{equation}
where $\Gamma^G $ and $\Gamma^{\mathfrak{g}^\ast} $ are the $G$ and
$\mathfrak{g}^\ast$ components of $\Gamma ^h $, respectively, that is, $
\Gamma^h:=\left(\Gamma^G,
\Gamma^{\mathfrak{g}^\ast}
\right)
$. In the left trivialized representation, the reduced Poisson and symplectic
Hamiltonians are simply the restrictions $h^{\mathfrak{g}^\ast}$ and
$h^{\mathcal{O}_{\mu}}$ of
$h$ to $\mathfrak{g}^\ast$ and to the coadjoint orbits $\mathcal{O}_{\mu}\subset
\mathfrak{g}^\ast$, respectively. Additionally, the  reduced stochastic Hamilton
equations on $\mathfrak{g}^\ast $ and $\mathcal{O}_{\mu} $ are given by 
\begin{equation}
\label{reduced trivialized equations}
\delta\Gamma^{\mathfrak{g}^\ast}=\sum _{i=1}^r \mbox{\rm ad} ^\ast _{\frac{\delta
h^{\mathfrak{g}^\ast} _i}{\delta
\Gamma ^{\mathfrak{g}^\ast}}}\Gamma ^{\mathfrak{g}^\ast}
\delta X ^i\text{\quad and \quad}
\delta\Gamma^{\mathcal{O}_{\mu}}=\sum _{i=1}^r \mbox{\rm ad} ^\ast _{\frac{\delta
h^{\mathcal{O}_{\mu}} _i}{\delta
\Gamma ^{\mathcal{O}_{\mu}}}}\Gamma ^{\mathcal{O}_{\mu}}
\delta X ^i
\end{equation}
where $h^{\mathfrak{g}^\ast}=\sum_{i=1}^{r} h _i^{\mathfrak{g}^\ast}
\epsilon^i$ and $h^{\mathcal{O}_{\mu}}=\sum_{i=1}^{r} h _i^{\mathcal{O}_{\mu}}
\epsilon^i$.

The combination of expressions~(\ref{stochastic trivialized vector field})
and~(\ref{reduced trivialized equations}) shows that in this setup, the dynamical
reconstruction of reduced solutions is particularly simple to write down. Indeed,
suppose that  we are given a solution $\Gamma^{\mathfrak{g}^\ast} $ of, say,  the
Poisson reduced system. In order to obtain the solution $\Gamma^h $ of the original
system such that $\Gamma ^h_0=(\Gamma ^G _0, \Gamma^{\mathfrak{g}^\ast}_0)
$ and
$\pi(\Gamma^h)=
\Gamma^{\mathfrak{g}^\ast}$, with $\pi: T ^\ast G\simeq G \times  \mathfrak{g}^\ast
\rightarrow T ^\ast G/G \simeq \mathfrak{g}^\ast$ the Poisson reduction projection,
it suffices to solve the stochastic differential equation in $G$
\begin{equation}
\label{stochastic trivialized vector field reconstruction}
\delta\Gamma^G=\sum _{i=1}^r  T _eL_{ \Gamma^G}\left(\frac{\delta h _i}{\delta
\Gamma ^{\mathfrak{g}^\ast}}\right)\delta X ^i,
\end{equation}
with the initial condition $\Gamma^G_0$. The reconstructed solution that we are
looking for is then $
\Gamma^h=\left(\Gamma^G,
\Gamma^{\mathfrak{g}^\ast}
\right).$

\subsection{Stochastic perturbations of the free rigid body}
\label{Stochastic perturbations of the free rigid body}

The free rigid body, also referred to as Euler top, is a particular case of systems
introduced in the previous section where the group $G$ is $SO (3, \mathbb{R})$. We
recall that in the context of mechanical systems on groups, a Hamiltonian system is
called \emph{free} 
when the energy of the system is purely kinetic and there is no
potential term. Let  $\left( \cdot ,\cdot
\right) $ be a left  invariant Riemannian metric on $G$; the kinetic energy $E$
associated to $\left( \cdot ,\cdot \right) $  is
$ E\left( v\right) =\frac{1}{2}\left( v,v\right) $, $v\in TG$. Then, using the left
invariance of the metric, we can write in body coordinates
\begin{equation*}
E   \left( g,\xi  \right) =\frac{1}{2}\left( \xi ,\xi \right)
_{e}=\frac{1}{2}\left\langle I\left( \xi \right) ,\xi \right\rangle, 
\end{equation*}%
for any $\left( g,\xi \right) \in G\times \mathfrak{g}$, where $\left\langle
\cdot ,\cdot \right\rangle $ is the natural pairing between elements of $%
\mathfrak{g}^{\ast }$ and $\mathfrak{g}$, and $I:\mathfrak{g\rightarrow g}^{\ast }$
is the map given by $\xi   \longmapsto   \left( \xi ,\cdot \right) _{e}$ and usually
known as the 
{\bfi inertia tensor} associated to the metric $(\cdot , \cdot )$.
The Legendre transformation associated to $E$ can be used to define a
Hamiltonian function $h:T^{\ast }G\rightarrow \mathbb{R}$ that, in body coordinates,
can be written as
\begin{equation}
h \left( g,\mu  \right) =\frac{1}{2}\left\langle \mu
,\Lambda \left( \mu \right) \right\rangle, 
\label{hamiltoniana en dual algebra}
\end{equation}%
where $\Lambda =I^{-1}$. Notice that as the kinetic
energy is left invariant (invariant with respect to the lifted $G$-action  to $T ^\ast  G  $ of  the action of $G$ on itself by left translations), then the components of $\mathbf{J}_{L}$ are
conserved quantities of the corresponding Hamiltonian system.
In order to connect with example in Section~\ref{Stochastic collective Hamiltonian
motion}, let $f\in C^{\infty }\left( \mathfrak{g}^{\ast }\right) $ be the function
$f:\mathfrak{g}^{\ast }  \rightarrow   \mathbb{R} $ given by $
\mu   \mapsto   \frac{1}{2}\left\langle \mu ,\Lambda \left( \mu
\right) \right\rangle$.
By~(\ref{momentum cotangent body}), the Hamiltonian function $h$ may be
expressed as $h=f\circ \mathbf{J}_{R}$. Therefore $h$ is collective with
respect to $\mathbf{J}_{R}$.

We now go back to the free rigid body case, that is, $G=SO\left( 3,\mathbb{R}\right)
$. We recall that the Lie algebra $\mathfrak{so}\left( 3,\mathbb{R}\right) $ is the
vector space of three dimensional skew-symmetric real matrices whose
bracket is just the commutator of two matrices. As a Lie algebra, $\left( 
\mathfrak{so}\left( 3\right) ,\left[ \cdot ,\cdot \right] \right) $ is naturally
isomorphic to $\left( \mathbb{R}^{3},\times \right) $, where $\times $
denotes the cross product of vectors in $\mathbb{R}^{3}$. Under this
isomorphism, the adjoint representation of $SO\left( 3,\mathbb{R}\right) $
on its Lie algebra is simply the action of matrices on vectors of $\mathbb{%
R}^{3}$ and the Lie-Poisson structure on $\frak{so}(3)^\ast \simeq \mathbb{R} ^3$ is
given by
$\{f,g\}(v)=-v \cdot \left(\nabla f \times \nabla g \right)$, for any $f , g  \in
C^\infty(\mathbb{R} ^3) $, where $\nabla $ is the usual Euclidean gradient and $\cdot $ denotes the Euclidean inner product.

Given a free rigid body with inertia tensor $I: \mathbb{R}^3 \rightarrow
\mathbb{R}^3$, since $\delta h _B/ \delta\mu= \Lambda (\mu) $, for any $\mu \in 
\mathbb{R}^3$, the left-trivialized equations of motion of the system are
\begin{equation}
\label{euler equation fulls}
(\dot A, \dot \mu)= \left(A \cdot  \widehat{\Lambda (\mu)}, \mu\times \Lambda
(\mu)\right),
\end{equation}
where the dot in the right hand side of~(\ref{euler equation fulls}) stands for matrix 
multiplication and
$\widehat{\Lambda (\mu)} $ is the skew-symmetric matrix associated to
$\Lambda (\mu)\in \mathbb{R}^3$ via the mapping that implements the Lie algebra
isomorphism between $\left( 
\mathfrak{so}\left( 3\right) ,\left[ \cdot ,\cdot \right] \right) $ and $\left(
\mathbb{R}^{3},\times \right) $.  In the context of the free rigid body motion
 the momentum map
$\mathbf{J}_{L}$ (respectively, $\mathbf{J}_R $) is called 
{\bfi spatial angular momentum} (respectively, {\bfi body angular momentum}). The
second component of~(\ref{euler equation fulls}), that is, 
\begin{equation}
\label{euler equations before}
\dot \mu= \mu\times
\Lambda (\mu) 
\end{equation}
are the  well-known {\bfi   Euler equations } for the free rigid
body.

\medskip

\noindent \textbf{Random perturbations of the body angular momentum}. We now
introduce  stochastic
perturbations of the free rigid body by using some of the geometrical tools that we
have introduced above. Later on we will compare this example with the
model of the randomly perturbed rigid body studied in \cite{Liao 1997} and
\cite{Liao 2005}, whose physical justification, as we will briefly discuss, involves
the same ideas as ours.

Let $V=\mathbb{R} \times \mathfrak{so}(3)\simeq \mathbb{R}^+ \times \mathbb{R} ^3$
and let
$h$ be the Hamiltonian function $h:T^{\ast }SO(3)\rightarrow V^{\ast
}=\mathbb{R} \times
\mathfrak{so}(3) ^{\ast }$ defined as
$
h  =\left( h_{0}  ,\mathbf{J}_{R} 
\right)$,
where $h_{0}$ is the Hamiltonian function of the free
(deterministic) rigid body. Observe
that $h$ is a left-invariant function because so is $\mathbf{J}_{R}$. Let $Y:%
\mathbb{R}_{+}\times \Omega \rightarrow \mathfrak{g}$ be a continuous
semimartingale which we may suppose, for the sake of simplicity, that it is a $
\mathfrak{g}$-valued Brownian motion and let $X:\mathbb{R}_{+}\times
\Omega \rightarrow \mathbb{R}^{\ast }\times \mathfrak{g}$ be the semimartingale
defined as $X_{t}\left( \omega \right) =\left( t,Y_{t}\left( \omega \right)
\right) $ for any $\left( t,\omega \right) \in \mathbb{R}\times \Omega $.
Consider the stochastic Hamiltonian system on $T ^\ast G $ associated to $h $ and
$X $. Since $h$ is left
invariant, the momentum  map $\mathbf{J}_{L}$ is preserved by the solution
semimartingales of this system and moreover, we
can apply the reduction scheme introduced in the previous sections. For example, if
we carry out Poisson reduction we have a reduced Hamiltonian function
$h^{\mathfrak{g}^{\ast }}:\mathfrak{g}^{\ast }\rightarrow V^{\ast }$ given by 
$h^{\mathfrak{g}^{\ast }} (\mu)=  \left( \frac{1}{2}\left\langle \mu ,\Lambda \left( \mu
\right) \right\rangle ,\mu \right)$.
Let $\left\{ \xi _{1}, \xi_2,\xi
_{3}\right\}$ a basis of the Lie algebra $\mathfrak{g}$ and
$\left\{ \epsilon ^{1},\epsilon^2,\epsilon ^{3}\right\} \subset \mathfrak{g}%
^{\ast }$ its dual basis. Observe that if we write $\mathbf{J}
_{R}\left( \mu \right)
=\sum_{i=1}^{3}\left\langle \mu ,\xi _{i}\right\rangle \epsilon ^{i}$ and
$Y=\sum_{i=1}^3 Y ^i \xi_i $, then the reduced stochastic Lie-Poisson equations can
be expressed as
\begin{equation}
\delta \mu _{t}= \mu _{t}\times \Lambda  \left(\mu _{t} \right) \delta
t+\sum_{i=1}^{3}\left(\mu _{t} \times \xi _{i} \right) \delta Y_{t}^{i}.
\label{dual algebra reduced model 1}
\end{equation}%
Regarding the reconstruction of the reduced dynamics, one has to solve the
stochastic differential equation on the rotations group $SO (3) $ given
by~(\ref{stochastic trivialized vector field reconstruction})  that, in this
particular case, is given by
\begin{equation}
\delta A_{t}=A _t \cdot \widehat{ \Lambda  \left(\mu _{t}\right)}
\delta t+\sum_{i=1}^{3}A _t \cdot \widehat{ \xi _{i}}
\delta Y_{t}^{i}.  \label{equation on group model 1}
\end{equation}%

A physical model whose description fits well in a stochastic Hamiltonian differential
equation like the one associated to $h$ and $X$ is that of a free
rigid body subjected to small random
impacts. Each impact  causes a small and instantaneous change in the \textit{body}
angular momenta
$\mu _{t}$ at time 
$t$ that justifies the extra term in~(\ref{dual algebra reduced model 1}), when
compared to the Euler equations~(\ref{euler equations before}). 

Our model is very
similar to the one proposed in~\cite{Liao 1997} where, instead of introducing the
random perturbation by means of a Hamiltonian function, a stochastic differential
equation on the group $G$ is introduced. This equation, also studied in detail in
\cite{Liao 2005}, is
\begin{equation}
\delta A_{t}=A_{t}\cdot \Lambda  \cdot {\rm Ad}
_{A_{t} }^{\ast }\left( \alpha \right) \delta t+\sum_{i=1}^{3}\left(
A_{t}\cdot  \Lambda \cdot  {\rm Ad}_{A_{t} }^{\ast
}\left( \epsilon ^{i}\right) \right) \delta Y^{i},
\label{liao starting equation}
\end{equation}
where $\alpha \in \mathfrak{g}^{\ast }$ is a constant vector. It important
to note that the drift terms of  equations (\ref{equation on group model 1})
and (\ref{liao starting equation}) coincide. Indeed,  for any $(g, \mu)\in G\times
\mathfrak{g}^\ast$ we can write
\begin{equation*}
\mu ={\rm Ad}_{g}^{\ast }\circ {\rm Ad}_{g^{-1}}^{%
\ast }\left( \mu \right) ={\rm Ad}_{g}^{\ast }\left(  
\mathbf{J}_{L} \left( g,\mu \right) \right). 
\end{equation*}%
Since in our model the
spatial angular momentum is conserved, $\Lambda \left( \mu _{t}\right) =\Lambda
\left( {\rm Ad}_{A_{t}}^{\ast }\left( {\rm Ad}_{A_{t}^{-1} }^{\ast } \mu_t\right)
\right)=\Lambda
\left( {\rm Ad}_{A_{t}}^{\ast }\left( \mathbf{J}_L(A _t, \mu _t)\right)
\right)=\Lambda
\left( {\rm Ad}_{A_{t}}^{\ast }\left(\alpha\right)
\right)
$, where $\alpha =
\mathbf{J}_{L}(A _t, \mu _t) $ is the preserved value of the spatial angular
momentum of a solution $(A _t, \mu _t)$ of~(\ref{dual algebra reduced model 1})
and~(\ref{equation on group model 1}). The difference between (\ref{equation on
group model 1}) and (\ref{liao starting equation}) lies in the stochastic
terms. The justification given by the author in~\cite{Liao 1997} for 
the equation~(\ref{liao starting equation}) is the following: since in  the
(deterministic) rigid body the spatial angular momentum  $\mathbf{J}_{L}$ is
conserved,  once we have fixed the value of this conserved quantity, we can
simply study the dynamics of the free rigid body by looking at the first
component of the ordinary differential equation (\ref{euler equation fulls}), now
rewritten as%
\begin{equation}
\dot{A}=A\left( \Lambda \left( {\rm Ad}_{A}^{\ast
}\left( \alpha \right) \right) \right)   \label{eq aux 6}
\end{equation}%
where $\alpha \in \mathfrak{g}^{\ast }$ is the $\mathbf{J}_{L}$-value of the
solution. Under random impacts, the spatial angular
momentum
$\alpha $, which was preserved in the deterministic case, is now randomly
modified. The idea is then to replace $\alpha dt$ in (\ref{eq aux 6}) by $\alpha
\delta t+\sum_{i=1}^{3}\epsilon ^{i}\delta Y^{i}$. Unlike our model, where the
random perturbation is introduced in the cotangent bundle respecting the underlying
symmetries of the deterministic system, there is no preservation of
$\mathbf{J}_{L}$ in the stochastic model of \cite{Liao 1997}. 

One
advantage of working on $T^{\ast }G
$ is that, even in the stochastic context,
classical quantities such as the angular momentum, are still well defined.
These objects do not have a clear counterpart if one follows the configuration space
based approach in  \cite{Liao 1997} (see for instance \cite{Liao 2005} for a non-trivial definition of angular
velocity in the stochastic context).

\bigskip

\noindent \textbf{Not so rigid rigid bodies. Random perturbation of the inertia
tensor}. In this example we want to write the equations that describe a rigid
body some of whose parts are slightly loose, that is,
the body is not a true rigid body and hence its mass distribution is constantly
changing in a random way. This will be modelled by stochastically perturbing the
tensor of inertia. 

For the sake of
simplicity, we will write $ G=SO\left( 3,\mathbb{R}\right) $ and
$\mathfrak{g}=\frak{so}(3)$. Let $\mathcal{L}\left( \mathfrak{g}^{\ast
},\mathfrak{g}\right) $ be the vector space of linear maps from $\mathfrak{g}^{\ast
}$ to $\mathfrak{g}$. As we know $\left(\mathfrak{so}\left( 3\right), [ \cdot ,
\cdot ] \right)
\simeq
\left(
\mathbb{R}^{3},\times
\right) $. Furthermore, we can establish an isomorphism $\mathbb{R}^{3}\simeq
\left( \mathbb{R}^{3}\right) ^{\ast }$ using the Euclidean inner product and hence
we can write  $\mathfrak{g}\simeq \mathfrak{g}^{\ast }$. Let $V=\mathcal{L}_{S}\left( \mathfrak{g}^{\ast },%
\mathfrak{g}\right) =\left\{ A\in \mathcal{L}\left( \mathfrak{g}^{\ast },%
\mathfrak{g}\right) ~|~A^{\ast }=A\right\} $ be the vector space of
selfadjoint linear maps from $\mathfrak{g}^{\ast }$ to $\mathfrak{g}$.
Define the Hamiltonian $h:T^{\ast }G\rightarrow V^{\ast }$ in \textit{body
coordinates} as%
\begin{equation*}
\begin{array}{rcc}
h:T^{\ast }G\simeq G\times \mathfrak{g}^{\ast } & \longrightarrow  & V^{\ast }
\\ 
\left( g,\mu \right)  & \longmapsto  & \bar{\mu},%
\end{array}%
\end{equation*}%
where $\bar{\mu}$ is such that 
\begin{equation*}
\begin{array}{rcl}
\bar{\mu}:\mathcal{L}_{S}\left( \mathfrak{g}^{\ast },\mathfrak{g}\right)  & 
\longrightarrow  & \mathbb{R} \\ 
A & \longmapsto  & \frac{1}{2}\left\langle \mu ,A\left( \mu \right)
\right\rangle .%
\end{array}%
\end{equation*}%
Observe that in body coordinates the Hamiltonian $h$ does not depend on $G$,
so the Hamiltonian is $G$-invariant by the action $\bar{\Phi}_{L}$ on $%
T^{\ast }G$. On the other hand, consider some filtered probability space $%
\left( \Omega ,\mathcal{F},\left\{ \mathcal{F}_{t}\right\} _{t\in \mathbb{R}%
},P\right) $ and introduce a stochastic component $X:\mathbb{R}_{+}\times \Omega
\rightarrow V$ in the following way:%
\begin{equation*}
\begin{array}{rcc}
X:\mathbb{R}_{+}\times \Omega  & \longrightarrow  & \mathcal{L}_{S}\left( \mathfrak{g}%
^{\ast },\mathfrak{g}\right)  \\ 
\left( t,\omega \right)  & \longmapsto  & \Lambda t+\varepsilon A_{t}\left(
\omega \right) ,%
\end{array}%
\end{equation*}%
where $\Lambda \in \mathcal{L}_{S}\left( \mathfrak{g}^{\ast },\mathfrak{g}%
\right) $ plays the role of the inverse of the tensor of inertia given
by the deterministic (rigid) description of the body, $\varepsilon $ is a
small parameter, and $A$ is an arbitrary $\mathcal{L}_{S}\left( \mathfrak{g}%
^{\ast },\mathfrak{g}\right) $-valued semimartingale. In order to show 
how the stochastic Hamiltonian system on  $T ^\ast G $ associated to
$h$ and $X$ models a free rigid body whose inertia tensor undergoes random
perturbations, 
we write down the associated stochastic reduced Lie-Poisson equations in
Stratonovich form
\begin{equation*}
\delta \mu _{t}=\mu _{t}\times \Lambda \left( \mu _{t}\right) \delta
t+\varepsilon \mu _{t}\times \delta A_{t}\left( \mu _{t}\right). 
\end{equation*}
Thus we see that these Lie-Poisson equations consist in changing $\Lambda
\left( \mu _{t}\right) d t$ in the Euler equations~(\ref{euler equations before}) by
$\Lambda
\left(
\mu _{t}\right) \delta t+\varepsilon \delta A_{t}\left( \mu _{t}\right) $, which
accounts for the stochastic perturbation of the inertia tensor.

\bigskip

\noindent\textbf{Acknowledgments} 
We thank an anonymous referee for comments that have much improved this paper.
The authors acknowledge partial support
from the French Agence National de la Recherche, contract number JC05-41465.
J.-A. L.-C. acknowledges support from the Spanish Ministerio de Educaci\'on
y Ciencia grant number BES-2004-4914. He also acknowledges partial support
from MEC grant BFM2006-10531 and Gobierno de Arag\'on grant DGA-grupos
consolidados 225-206. J.-P. O. has been partially supported by a ``Bonus
Qualit\'e Recherche" contract from the Universit\'e de Franche-Comt\'e.

\end{document}